\documentclass{article}
\usepackage[numbers,sort&compress,super,comma]{natbib} % this must be before neurips_2025
\usepackage[preprint]{neurips_2025}
\usepackage[utf8]{inputenc}
\usepackage[T1]{fontenc}
\usepackage{lmodern}
\usepackage[english]{babel}
\usepackage{amssymb, amsmath, amsthm, amsfonts}
\usepackage{thmtools, mathtools, mathrsfs, dsfont}
\usepackage{bbm}
\usepackage{bm} % bold greek
\usepackage{forloop}
\usepackage[pdftex]{graphicx}  
\usepackage{subcaption}
\usepackage{placeins}      %forcing figure placement
\usepackage{booktabs}
\usepackage{multirow}
\usepackage[shortlabels]{enumitem} %easy customization of itemize/enumerate lists
\PassOptionsToPackage{dvipsnames}{xcolor}
\usepackage{xcolor}
\definecolor{ForestGreen}{rgb}{0.13, 0.55, 0.13}
\definecolor{MidnightBlue}{rgb}{0.1, 0.1, 0.44}
\definecolor{BurntOrange}{rgb}{0.8, 0.33, 0.0}
\definecolor{Plum}{rgb}{0.56, 0.27, 0.52}
\usepackage[colorlinks=true,linkcolor=MidnightBlue,citecolor=ForestGreen,filecolor=TealBlue,urlcolor=Plum]{hyperref}
\hypersetup{breaklinks=true}
\usepackage{tikz}
\usetikzlibrary{positioning}
\usepackage{tikz-cd}
\usepackage{algorithm}
\usepackage{algpseudocode}
\definecolor{mpcolor}{rgb}{1, 0.1, 0.59}
\definecolor{ascolor}{rgb}{1, 0.5, 0.}

\newtheorem{proposition}{Proposition}
\newtheorem{lemma}{Lemma}

\theoremstyle{definition} \newtheorem{definition}{Definition}  
\theoremstyle{remark} 

\usepackage{letltxmacro}
\LetLtxMacro{\originaleqref}{\eqref}
\renewcommand{\eqref}{Eq.~\originaleqref}

\usepackage{forloop}
\newcommand{\defvec}[1]{\expandafter\newcommand\csname v#1\endcsname{{\mathbf{#1}}}}
\newcounter{ct}
\forLoop{1}{26}{ct}{
    \edef\letter{\alph{ct}}
    \expandafter\defvec\letter
}
\forLoop{1}{26}{ct}{
    \edef\letter{\Alph{ct}}
    \expandafter\defvec\letter
}

\usepackage{expl3}

\ExplSyntaxOn
\cs_new:Nn \define_bold_greek:
 {
  \seq_map_inline:Nn \g_greek_letters_seq
   {
    \cs_new:cpn { v##1 } { \bm { \csname ##1 \endcsname } }
   }
 }

\seq_new:N \g_greek_letters_seq
\seq_set_split:Nnn \g_greek_letters_seq { , }
 {
  alpha,beta,gamma,delta,epsilon,zeta,eta,theta,iota,kappa,lambda,
  mu,nu,xi,pi,rho,sigma,tau,upsilon,phi,chi,psi,omega,
  Delta
 }

\define_bold_greek:
\ExplSyntaxOff
\DeclareMathOperator*{\tr}{tr} % trace

\DeclareMathOperator*{\argmin}{\rm argmin}

\DeclarePairedDelimiter{\norm}{\lVert}{\rVert}

\newcommand{\tconju}{\mathbin{\overset{\text{t.c.}}{\sim}}}
\DeclareMathOperator*{\E}{\rm E} % expectation
\newcommand{\inv}{^{-1}}
\newcommand{\dm}[1]{\ensuremath{\mathrm{d}{#1}}} % dx dy dz dmu
\newcommand{\field}[1]{\ensuremath{\mathbb{#1}}}
\newcommand{\reals}{\field{R}}

\newcommand{\trp}{{^\top}} % transpose
\newcommand{\identity}{\ensuremath{\mathbb{I}}}

\newcommand{\posReals}{\reals^{+}}
\newcommand{\inputSpace}{\mathbb{X}}
\newcommand{\dcomplexity}{\ensuremath{d_\text{cxty}}}
\newcommand{\initDist}{\mu} % natural probability distribution over the points on the phase space

\newcommand{\homeo}{\Psi} % NOT \Phi: \phi^t is the flow, so \Phi vs \phi clash on the same letter (matches the talk)
\newcommand{\invhomeo}{\homeo\inv}

\AtBeginDocument{\setlength\abovedisplayskip{4pt}}
\AtBeginDocument{\setlength\belowdisplayskip{4pt}}
\AtBeginDocument{\setlength{\intextsep}{9pt}} % Vertical space above & below [h] floats
\AtBeginDocument{\setlength{\textfloatsep}{9pt}} % Vertical space below (above) [t] ([b]) floats
\usepackage[compact]{titlesec}
\titlespacing*{\section}{0pt}{2pt plus 2pt minus 2pt}{2pt plus 2pt minus 2pt}
\titlespacing*{\subsection}{0pt}{2pt plus 2pt minus 2pt}{0pt}

\title{Dynamical Archetype Analysis:\\Autonomous Computation}
\author{
\'Abel S\'agodi, Il Memming Park\\
    Champalimaud Centre for the Unknown\\
    Champalimaud Foundation, Lisbon, Portugal\\
    \texttt{\{abel.sagodi, memming.park\}@research.fchampalimaud.org} \\
}

\begin{document}
\maketitle

% keywords: interpretability, dynamical dissimilarity, dynamical motifs, neural computation, working memory, computation through dynamics

\begin{abstract}
The study of neural computation aims to understand the function of a neural system as an information processing machine.
Neural systems are undoubtedly complex, necessitating principled and automated tools to abstract away details to organize and incrementally build intuition.
We argue that systems with the same effective behavior should be abstracted by their ideal representative, i.e., \emph{archetype}, defined by its asymptotic dynamical structure. 
We propose a library of archetypical computations and a new measure of dissimilarity that allows us to group systems based on their effective behavior by explicitly considering both deformations that break topological conjugacy and diffeomorphisms that preserve it.
Crucially, because this dissimilarity can be estimated directly from observed trajectories, it provides a practical tool for analyzing neural dynamics.
Numerical experiments demonstrate our method's ability to overcome the previously reported fragility of existing (dis)similarity measures for approximate continuous attractors and high-dimensional recurrent neural networks.
Although our experiments focus on working memory systems, our theoretical approach naturally extends to the general mechanistic interpretation of recurrent dynamics in both biological and artificial neural systems.
We argue that abstract dynamical archetypes, rather than detailed dynamical systems, offer a more useful vocabulary for describing neural computation.
%
% Understanding neural computation in dynamical neural systems often requires interpreting the behavior of large, black-box models.
% This poses challenges for both interpretability and cross-model comparison--particularly when comparing models from different sources, such as two animals.
% To address this, we propose a framework for interpreting neural computation for working memory through \emph{dynamical motifs}--canonical, human-understandable building blocks of dynamics.
% Within this framework, we introduce a model reduction method that approximates the essential computational mechanism in a system by identifying small perturbations that render two systems topologically conjugate.
% Our method fits the homeomorphism that establishes the topological conjugacy as a Neural ODE.
% We derive measures of dissimilarity based on our homeomorphism-conjugacy mapping method.
% Our method uses a curated library of low-dimensional canonical motifs (e.g., limit cycles and various continuous attractors) to explain the core dynamics.
% This allows us to define and quantify similarity between systems in terms of shared dynamical structure.
~% The extracted motifs are robust to parametric variation, suggesting they reflect stable computational strategies rather than model idiosyncrasies.
% We demonstrate the approach on various toy examples, most importantly approximate continuous attractors.
% Furthermore, we extract an approximate ring attractors from a high-dimensional RNNs trained on a angular velocity integration task, and show how the found homeomorphism captures the complexity of the geometry of the embedding of the ring.
% Our results highlight the utility of this method for robust interpretation and cross-system comparison within the computation-through-dynamics framework.
\end{abstract}

\section{Introduction}\label{sec:intro}
\begin{figure}[htbp]
\centering
\includegraphics[width=\textwidth]{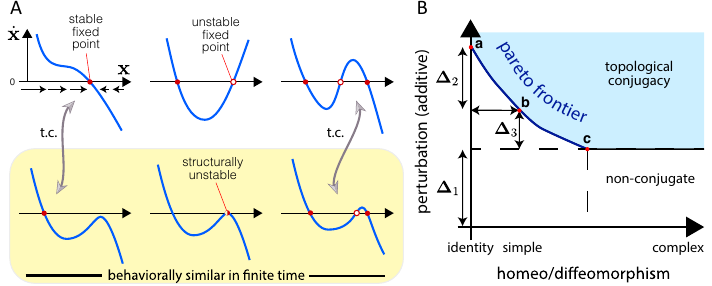}
\caption{
(\textbf{A})
Fixed points (red circles) characterize the asymptotic behavior and topological conjugacy (t.c.) in 1-D dynamical systems.
%They are topologically equivalent (t.c.) if their fixed points correspond to each other in the same order (2 pairs in the figure).
Although the systems at the bottom are not t.c.\ to each other, they behave similarly in the short time scale.
(\textbf{B})
Behavioral dissimilarity consists of two measures of distortion to map a source system to a target system.
Their trade-off forms a boundary of minimum amount of distortion needed.
See~\eqref{eq:decomp:triple} for the detailed explanation of the symbols.
}
\label{fig:pareto}
\end{figure}

Scientific understanding often advances not by improving predictions, but by changing the tools we use to describe systems.
Ordinary differential equations (ODEs) are a foundational modeling tool of choice for neural computation, used to describe both empirical signals and theoretical constructs.
For example, state-space models can be used to train a black-box that captures the spatiotemporal structure of time series data as an ODE~\citep{nair2023approximate,pei2neural,Dowling2024b}.
%Through abstraction and idealization, we can distill the essential features of neural trajectories, enabling meaningful understanding rather than merely presenting a detailed but opaque description.
A further abstraction, now widely adopted, is to focus on their long-term or asymptotic behavior---discretizing the continuous state space into attractors, their basins, and their relation to the lower-level implementation (e.g. neurons) and higher-level information processing (e.g. input and context).
These abstractions have yielded influential models of memory\citep{machens2008ca}, perceptual decision-making\citep{mante2013context}, and task trained recurrent neural networks (RNNs)\citep{driscoll2024flexible}, offering simple and intuitive explanations for otherwise high-dimensional dynamics~\citep{vyas2020ctd,versteeg2025computation}.
Yet, this abstraction of the mechanism is often based on confusing the effective and asymptotic behaviors of the ODE implementation, neglecting that systems with different asymptotic structures can exhibit indistinguishable behavior over task-relevant time scales.
%While canonical models aim to distill these computations into minimal motifs \citep{chirimuuta2014minimal}, actual biological neural implementations may vary significantly, motivating the need for flexible yet interpretable abstractions.
We need to conceptually separate the ODEs and their interpretable discrete abstraction, while making their connection precise through quantification of similarity in terms of their effective behavior at finite time scales (Fig.~\ref{fig:pareto}A).

%Understanding the relationship between neural activity and behavior is a fundamental question in neuroscience.
%Understanding neural computation requires more than capturing the details of neural dynamics; it necessitates the development of conceptual models that simplify and abstract the complexity of these processes into forms that are humanly interpretable.

%why no 20D attractors?
%Recent research has shown that, despite the inherent complexity of neuronal networks, their activity often resides on a low-dimensional manifold, representing only a small fraction of the system's full potential complexity \citep{duncker2021dynamics}.
%This empirical finding suggests that neural computation is governed by underlying constraints--potentially driven by task-relevant computations or external environmental factors--that allow for such a reduction in complexity.
%Such constraints not only make the idealization of neural systems feasible but also provide a path toward more interpretable models of neural activity.

%%%%
%\subsection{Memory maintenance in neural networks}\label{sec:ctd} 
%Discovering how the brain contributes to behavior is the central aim of systems neuroscience.
%We need a language that can describe how neural populations transform inputs into goal-directed behavior.
%Here we focus on the memory maintenance part of this transformation.

\section{Background}
\subsection{Equivalence of Qualitative Behavior}
Dynamical systems are widely used for studying neural computation, offering insights into cognitive functions such as working memory, short-term memory, and motor control~\citep{beer1995ctrnn, beer2006parameterspace, sussillo2014neural, vyas2020ctd,versteeg2025computation}.
Discrete and continuous forms of working memory, for example, can be modeled as attractor dynamics that maintain persistent neural activity \citep{wang2001synaptic,wimmer2014ca,masse2019circuit,panichello2021shared,zhang2022translation, hoeller2024bridging,ritter2025efficient}.
Short-term memory, involving transient yet structured activity, has been explored in dynamical models capturing task-dependent variability \citep{kurtkaya2025dynamical}.
% Motor control, too, exhibits low-dimensional dynamics shaped by recurrent structure \citep{wang2022representation}.

Since the time of Poincar\'e, dynamical system theory has heavily relied on the analysis of their \emph{qualitative behavior}---the eventual fate of the system, rather than the details of their trajectory.
For qualitative questions such as ``will the neural system eventually settle down to the same resting state regardless of the initial condition in the absence of external drive or will it oscillate forever?'', interestingly, often a simple and intuitive picture (technically a \emph{phase portrait}) can give a definite answer~\citep{Strogatz2000-tf}.

The space of all dynamical systems can be partitioned by the equivalence of their qualitative behavior~\citep{Arnold1994}, abstracting away the finite time transient behaviors and the details of the vector field.
Specifically, consider two autonomous dynamical systems defined on the same phase space $\inputSpace$:
\begin{align}
    \dot{\vx} = \vf(\vx),
    &\qquad
    \dot{\vy} = \vg(\vy).
%    \begin{dcases}
%	\dot{\vx} = \vf(\vx)
%	\\
%	\dot{\vy} = \vg(\vy)
%    \end{dcases}
\end{align}
%Poincar\'e introduced the notion of homeomorphism 
If there exists a homeomorphism $\homeo: \inputSpace \to \inputSpace$, a continuous invertible map from the trajectories (solutions) of one system to the other, then the two systems are \emph{topologically conjugate} (Def.~\ref{def:top_conj}) and we denote as $\vf(\cdot) \tconju \vg(\cdot)$.
Topological conjugacy defines an equivalence relation on the set of flows.
%Furthermore, if the map is differentiable, we have a diffeomorphism.

There are two major advantages of this abstraction:
(1) it is implementation agnostic, e.g., an insect neural circuit with 8 neurons can be topologically conjugate to an artificial neural network model where the time scale and the biophysical details may be inaccurate and (2) it promotes interpretability, i.e., it enables us to understand the qualitative behavior of a system in terms of simpler or more familiar models. 
It naturally emphasizes the use of attractor dynamics since they are asymptotically defined and preserved.
For neural computation, this means that only a finite collection of attractors remains to be analyzed.
Among those structures, two have a prominent status in the vernacular of neuroscience: fixed points and line attractors.
Mathematically, (hyperbolic) fixed points are locally conjugate to a linear system\citep{hartman1960lemma} which are well understood thanks to the easy decomposition and superposition of solutions.
Moreover, multiple isolated fixed points can serve as a substrate for discrete computational frameworks such as finite state machines or Markov chains.
Line attractors on the other hand can store an analog-valued representation (for indefinite time).
Theoretical neuroscience  has extensively used these as building blocks to help us expand our understanding of neural computation\citep{dayan2005theoretical,vyas2020ctd}.
While mathematically elegant, however, there are a few critical issues with this program.

\subsection{Structural Stability and Bifurcations}
If for any smooth vector field $\vDelta(\vx)$, there exists an $\epsilon > 0$, such that
\begin{align}\label{eq:vfp:epsilon}
\dot{\vx} &= \vf(\vx) + \epsilon\vDelta(\vx)
\end{align}
does not break the topological conjugacy, then the unperturbed system is said to be \emph{structurally stable}.
% In a parametric case, a recurrent neural network with weight vector $W$ perturbed by $\Delta W$ gives,
% \begin{align}
% 	\vf(\vx) = -\vx + W^\top\homeo(\vx),
% 	&\qquad
% 	\vDelta(\vx) = \nabla_W \vf\vert_\vx \Delta W + O(\Delta W^2) =
% \end{align}
% where $\homeo(\cdot)$ is a point nonlinearity and $\nabla_W$ denotes the Jacobian respect to $W$.
Since physical and numerical experiments are not perfect, the structurally unstable dynamical systems have been a big source of unease. % since Andronov and Pontryagin 1937) who advised to only work with structurally stable systems.
Unfortunately, several of the key asymptotic attractor structures that are very useful in expressing interpretable neural computation are not robust to vector field or synaptic weight perturbations, thus structurally unstable.
A notorious example is the continuous attractor (e.g., line attractor) dynamics~\citep{Sagodi2024a}.
Given that the biological neural systems are constantly fluctuating and both biological and artificial networks learn from noisy learning signals, structurally unstable neural computations are bound to bifurcate away~\citep{Park2023a}.
%Even structurally stable dynamics are known to go through bifurcations--breaking of topological conjugacy--during learning and necessarily so if the final dynamical system is not conjugate to the initial system~\citep{Doya1993-tz}.

%
\subsection{Methods to Compare and Reconstruct Dynamics}\label{sec:comparison}
%A foundational concept for the comparison of dynamics is \emph{topological conjugacy}, which formalizes the notion of two systems being dynamically equivalent under a homeomorphism that maps trajectories of one system to the other while preserving the temporal ordering (see Sec.~\ref{def:top_conj}).
%Another key approach focuses on asymptotic behavior, particularly through the study of $\omega$-limit sets, which describe the long-term evolution of trajectories.
%These sets help characterize the stability and attractor structure of a system and are often used to classify systems based on the geometry and topology of their invariant sets\citep{Jordan2019a}.
%Together, these methods provide a basis for principled comparisons of dynamical systems that go beyond superficial similarity of trajectories.

%geometry
A variety of methods compare the geometry of neural representations in latent spaces, including Representational Similarity Analysis (RSA) \citep{kriegeskorte2008representational,taylor2025framed}, Singular Vector Canonical Correlation Analysis (SVCCA) \citep{raghu2017svcca}, Procrustes Shape Distance \citep{williams2021generalized}, and Stochastic Shape Distance (SSD) \citep{lipshutz2024disentangling}, a probabilistic extension of Procrustes analysis \citep{duong2022representational, barbosa2025quantifying}; see also \citep{williams2024equivalence} for a broader comparison of such techniques.
A recently proposed method, Dynamical Similarity Analysis (DSA)\citep{ostrow2023beyond}, quantifies similarity between dynamical systems by comparing their local linearizations obtained through Dynamic Mode Decomposition (DMD) \citep{schmid2010dynamic}, thereby extending representational comparison beyond static geometry to include dynamics.
DSA constructs a high-dimensional embedding in which a nonlinear system is approximated by a linear vector field, capturing spatiotemporally coherent structure through so-called (truncated) Koopman modes. 
%To compare systems, DSA learns a linear coordinate transformation that maximizes the cosine similarity between the vector fields of their respective linear time-invariant approximations.
%comments:
When the linearization faithfully captures the underlying Koopman eigenspectrum, this provides an estimate of how close the systems are to being topologically conjugate. 
However, the reliance on spectral truncation can obscure critical nonlinear features---such as limit cycles, bifurcations, or multistability---that are often essential to understanding neural dynamics and computation (see Supp.Sec.~\ref{sec:dsa_sanity}).
%But we can capture some slow manifold behavior!
Complementary approaches extend geometric comparison to settings where neural population dynamics are noisy or stochastic.
Optimal transport–based metrics have been proposed to compare such systems by generalizing Procrustes shape distance to average trajectories and Stochastic Shape Distance (SSD) to marginal statistics or noise correlations \citep{nejatbakhsh2024comparing, lipshutz2024disentangling}.

%VF
Recent work has introduced several approaches for comparing or extracting structure from dynamical systems.
Time-Warp-Attend (TWA) was proposed to extract topologically invariant features by augmenting data with time-warped versions of the underlying vector field \citep{moriel2024timewarpattend}. DFORM uses invertible residual networks to learn orbit-matching coordinate transformations between systems, allowing comparisons across learned dynamical models \citep{chen2024dform, chen2026comparing}. % However, its orbital similarity loss (Eq.~5) serves only as a noisy proxy for topological equivalence.
Smooth Prototype Equivalences (SPE) aims to recover a diffeomorphism from observed data to known normal forms \citep{friedman2025characterizing}, but it has only been demonstrated on a narrow class of dynamics---specifically, comparisons between limit cycles and fixed points.
A key limitation of all these methods is their reliance on \emph{explicit vector field representations}, rather than working directly with observed trajectories.
This restricts their applicability in data-driven settings where only trajectory data is available and the vector field is unknown or intractable to estimate.

%realization
A central concept in comparing dynamical systems is that of \emph{realization} -- constructing a system that replicates the input-output behavior of another system\citep{grigoryeva2020dimension, gonon2023approximation}.
%A \emph{canonical} realization is the simplest system that produces the same effective behavior as a target system.
%Closely related is the notion of \emph{bisimilarity}, which formalizes behavioral equivalence through relations that match trajectories or observable outputs across systems \citep{vanderschaft2004bisimulation, vanderschaft2004equivalence, pola2004bisimulation, pola2006equivalence, tabuada2004bisimilar, girard2011approximate}.
%model identification of classes of cognitive models\citep{rmus2024artificial}
%
%reconstruction: 
Recent advances have sought to recover and interpret the full dynamics of complex systems from data, rather than relying solely on reduced or task-driven summaries \citep{durstewitz2023reconstructing, brenner2024almost, brenner2025learning}.
Several approaches focus on learning continuous-time models of latent stochastic dynamical systems, enabling structured representations of uncertainty and variability in neural activity \citep{duncker2019learning}.
For example, SSMLearn \citep{cenedese2022data} combines symbolic regression with state-space modeling to extract interpretable dynamical descriptions.
Other methods, such as Phase2vec \citep{ricci2022phase2vec}, learn embeddings that reflect the qualitative structure of dynamical systems---e.g., fixed points, limit cycles, or phase relationships---by organizing systems in a geometry shaped by their behavior.

%OTHER methods / miscellaneous
Other methods include MARBLE\citep{gosztolai2025marble}, which defines a similarity measure based on local vector field features (vectors and higher-order derivatives) sampled from across the phase space, where the features are embedded into lower-dimensional space through contrastive learning, and the embedded distributions from different systems can be compared using an optimal transport distance.
A more specialized approach, applicable only to limit-cycle dynamics, is the data-driven method for estimating the asymptotic phase of nonlinear oscillators using the Temporal 1-Form \citep{wilshin2025estimating}.

%\ascomment{Comparing noisy neural population dynamics with optimal transport: metric is an extension of Procrustes shape distance which compares average trajectories and SSD which compares marginal statistics or noise correlations\citep{nejatbakhsh2024comparing, lipshutz2024disentangling}.}

%Conley-Morse decomposition: ... \citep{arai2009database} \citep{kalies2017computational} \citep{chen2007vector}: periodic orbit extraction

%Koopman / linear  \citep{mezic2004comparison}

\section{Effective Behavior Similarity}\label{sec:effective_ds}
\subsection{Motivation and Measures of Dissimilarity}
Our guiding principle is to group dynamical systems by their similarity in what actually matters to neural computation---their \emph{effective behavior, not the qualitative behavior}.
%At the same time, to represent some groups by their most interpretable member.
%If two systems are equally useful for the same purpose, i.e. implementing a desired neural computation, then we will consider them ``the same''.
%What do we mean by effective behavior?
We develop criteria and measures of (dis)similarity that capture the concept of effective behavior.
Both artificial and biological agents perform neural computation in a finite amount of time.
For example, typical visuomotor behavior or perceptual-decision making tasks are performed within a few hundred milliseconds.
So if our goal is to understand the information processing capacity at the behaviorally relevant time scales, we should limit our analysis to the interval rather than the mathematically convenient infinite time scale.
We capture the behavioral time scale as a reasonable duration of time $[0, T_\text{max}]$.
This motivates the following time integrated point-wise distance measure:
For all $\vx \in \inputSpace$ and $t \in [0, T_\text{max}]$,
\begin{align}
    d_0(\vf(\vx), \vg(\vx))
	&\coloneq \int_0^{T_{\text{max}}} \norm{\phi_\vf^t(\vx) - \phi_\vg^t(\vx)}\, \dm{t},
\end{align}
where $\phi_\vf^t(\vx)$ is the solution (trajectory) of $\dot\vx = \vf(\vx)$ with initial condition $\vx$ at time $t$.

We must also restrict our analysis to the relevant portion of the phase space.
The dynamical system implemented defines a vector field over the entire space but typically most of the space is not utilized during normal operations.
In fact, neural population data analysis suggests that the trajectories are restricted to a low-dimensional manifold in the context of a ``task'' or  ``computation''.
Therefore, the distinction between the systems should only be compared at the utilized regions of the phase space.
We denote this natural probability distribution over the points on the phase space as $\initDist(\vx)$.
These two notions together motivate the space integrated measure (where $\E$ denotes the expectation operator):
\begin{align}
%     d_{0,\infty}(\vf, \vg; \initDist)
% 	&\coloneq \sup_{\vx \in \initDist(\vx)}
% 	d_0(\vf(\vx), \vg(\vx))
%     \\
    \bar{d}_{0}(\vf, \vg; \initDist)
	&\coloneq
	    \E_{\vx \sim \initDist(\vx)}
		d_0(\vf(\vx), \vg(\vx)).
	\label{eq:ed0}
\end{align}
If the axes of the space $\inputSpace$ are arbitrary, as in the case of state-space models with linear observation models, or theoretical models without a fixed implementation, we would want the measure to be invariant to isometric scaling, rotations, and/or affine transformations.
More general continuous coordinate transformation shall be allowed.
Let $\homeo: \inputSpace \to \inputSpace$ be a diffeomorphism mapping the $\vg$ system to the $\vf$ system; we refer to $\vg$ as the \textbf{source} system and $\vf$ as the \textbf{target} system.
Although topological conjugacy requires only a homeomorphism, we work in the smooth category and take $\homeo$ to be a diffeomorphism so that its deformation complexity (the deviation of $\nabla\homeo$ from the identity, introduced below) is well-defined; parametrizing $\homeo$ as a Neural ODE yields a diffeomorphism by construction.
We define a dissimilarity
\begin{align}\label{eq:ed1}
    {d}_{\text{traj}}(\vf(\vx), \vg(\vx); \homeo)
	&\coloneq 
	    \E_{\vx \sim \initDist(\vx)}
		\int
		    \norm{\phi_{\vf}^t(\vx) - \homeo(\phi_{\vg}^t(\homeo^{-1}(\vx)))}
		\,\dm{t},
\end{align}
which is non-negative, $0$ when topologically conjugate on the support of $\mu$, but not symmetric in general.
If a more complicated deformation is needed to minimize \eqref{eq:ed1}, it is natural to consider its effective behavior less similar\footnote{
We are relaxing the dichotomic nature of the topological conjugacy.
}.
We define  $\dcomplexity: (\inputSpace \to \inputSpace) \to \posReals$ to be a measure of complexity of the diffeomorphism such that the simplest map, the identity, gives $\dcomplexity(\homeo\colon\vx \mapsto \vx) = 0$, by the integral
%For example, we can define,
\begin{align}\label{eq:dcomplexity}
    \dcomplexity(\homeo) &\coloneq
    \E_{\vx \sim \initDist(\vx)}
    \norm{\nabla\homeo(\vx) - \identity} % + c_b \norm{\homeo(\vx) - \vx}
\end{align}
%for some $c_a, c_b > 0$.
%
%\subsection{Trade-off between perturbation and deformation}
We have arrived at a multi-valued dissimilarity between systems:
two systems are similar if,
(a) there is a small distance between trajectories after a continuous coordinate transform~\eqref{eq:ed1},
and (b) that transform has small complexity~\eqref{eq:dcomplexity}.
However, for a source--target pair, the two values are \emph{not uniquely determined}.
Below, we analyze the trade-off as a function of complexity.

First, we show that the trajectory distance measure is bounded by an additive vector field perturbation. 
Let $\vDelta(\cdot)$ be a perturbation to the vector field as in \eqref{eq:vfp:epsilon} with $\vDelta$ absorbing $\epsilon$.
%\begin{align}
%    \dot{\vx} &= \vf(\vx) + \vDelta(\vx).
%\end{align}
For a finite time horizon $[0, T_\text{max}]$, assuming a Lipschitz bound $L$ for $\vf(\cdot)$ in the relevant domain and a uniform bound on $\norm{\vDelta(\cdot)}$, we can bound the trajectory deviation from the unperturbed system.
By Gr\"onwall's inequality\cite{Howard2025}:
\begin{align}\label{eq:gronwall}
    \norm{
	\phi^t_\vf(\vx)
	-
	\phi^t_{\vf + \vDelta}(\vx)
    }
    &\leq
	e^{Lt} \int_0^{t} e^{-Ls} \norm{\vDelta(s)} \,\mathrm{d}s
    \leq
    \left(\sup_{\vs \in (0, t)} \norm{\vDelta(\vs)}\right) \frac{1}{L} (e^{Lt} - 1),
\end{align}
for any norm $\norm{\cdot}$ and $\vDelta(s) \coloneq \vDelta(\phi^s_\vf(\vx))$.
Thus, for small perturbation, the deviation of trajectories is also small for a short time:
\begin{align}\label{eq:d0bound}
    \bar{d}_0(\vf, \vf + \vDelta)
    &\leq
    \left(\sup_\vs \norm{\vDelta(\vs)}\right) \frac{1}{L^2} e^{LT_{\text{max}}}.
\end{align}

Now, let us consider the burden of distortion by using the best diffeomorphism that minimizes \eqref{eq:ed1}.
%Mapping the asymptotic structure of the two systems to align is a preferable distortion, however, 
Since diffeomorphisms cannot change the topological conjugacy class, we need the perturbation to take into account the slack if they are not conjugate.
Therefore, we first consider the minimum perturbation to achieve topological conjugacy (on the support of $\mu$).
\begin{align}
    \vDelta_1(\vx) = \argmin \norm{\vDelta}
    %\label{eq:d1:min}
    & \text{s.t.} \quad
    \vf(\vx) + \vDelta(\vx)
    \tconju
    \vg(\vx).
    \label{eq:d1:tc}
\end{align}
We identify this gap $\vDelta_1$ as the \emph{imperfection}.
By \eqref{eq:d0bound}, we have
\begin{align}
    \bar{d}_1(\vf, \vg; \homeo^\ast)
	&\leq
	\left(\sup_\vs \norm{\vDelta_1(\vs)}\right) \frac{1}{L^2} e^{LT_{\text{max}}},
\end{align}
where $\homeo^\ast$ achieves the desired topological conjugacy.
Note that $\vf \tconju \vg$ implies $\norm{\vDelta_1} = 0$.
Now, we want to consider the complexity of $\homeo$ and trade-off with (additive) perturbation.
To link the complexity of the diffeomorphism to the magnitude of the perturbation, we seek to find the equivalence up to first order of the vector field corresponding to the diffeomorphism.
Let
$
\vf(\vx)
    \tconju
\vg(\vx)
$.
That is, applying the diffeomorphism $\homeo$ to the source $\vg$ (as in \eqref{eq:ed1}) perturbs its vector field; we seek the equivalent additive perturbation $\epsilon\vDelta_2$,
\begin{align}
    \homeo_\ast\vg &= \vg + \epsilon\vDelta_2 + O(\epsilon^2),
    \label{eq:pareto:tc}
\end{align}
where $\epsilon > 0$.
We assume the following Ansatz for a diffeomorphism close to the identity, with first-order generator $\tilde\homeo$ (the $O(\epsilon)$ displacement field), together with the first-order expansions of its Jacobian and inverse:
\begin{align}\label{eq:h:ansatz}
    \homeo(\vx) &= \vx + \epsilon \tilde\homeo(\vx) + O(\epsilon^2),
    &
    \nabla\homeo(\vx) &= \identity + \epsilon\nabla\tilde\homeo(\vx) + O(\epsilon^2),
    &
    \homeo^{-1}(\vx) &= \vx - \epsilon\tilde\homeo(\vx) + O(\epsilon^2).
\end{align}
Expanding the pushforward $(\homeo_\ast\vg)(\vx)=\nabla\homeo(\homeo^{-1}(\vx))\,\vg(\homeo^{-1}(\vx))$ to first order (see Supp.~Sec.~\ref{sec:fojac}) identifies the induced perturbation with the Jacobian Lie bracket of $\vg$ and $\tilde\homeo$,
\begin{align}
    \vDelta_2 &= (\nabla\tilde\homeo)\,\vg - (\nabla\vg)\,\tilde\homeo = [\vg,\tilde\homeo] + O(\epsilon).
\end{align}
Adding the irreducible conjugacy slack $\vDelta_1$ of \eqref{eq:d1:tc}, which no diffeomorphism can remove, the total perturbation $\vDelta=\vDelta_1+\vDelta_2$ obeys the point-wise bound
\begin{align}
    \norm{\vDelta(\vx)}
    &\leq
    \norm{\vDelta_1(\vx)}
    + \norm{\vg(\vx)}\cdot\underbrace{\norm{\nabla\tilde\homeo(\vx)}}_{\substack{\text{scale c.}}}
    + \norm{\nabla\vg(\vx)}\cdot\underbrace{\norm{\tilde\homeo(\vx)}}_{\substack{\text{translation c.}}}
    + O(\epsilon).
    \label{eq:h:bound}
\end{align}
for appropriate vector and corresponding matrix norms.
We recognize two kinds of complexities of the diffeomorphism: (1) scale complexity which is a function of the Jacobian of diffeomorphism $\nabla\homeo(\vx)$, and (2) translation complexity.
Note that this form matches \eqref{eq:dcomplexity} if translation is considered invariant.

Under the assumption that $\norm{\vg(\vx)}$ and $\norm{\nabla\vg(\vx)}$ are bounded within the region of interest, we argue that a linear combination of the two complexity measures of the diffeomorphism can be traded off with the perturbation.
Therefore, we define the smallest perturbation needed to compensate for the suboptimal diffeomorphism with a given complexity $c$:
\begin{align}
    \vDelta_2(\vx ; c) \coloneqq \inf_{\dcomplexity{\homeo} = c}
	\nabla\homeo(\invhomeo(\vx)) \vg(\invhomeo(\vx))
	- \vg(\vx) - \vDelta_1(\vx).
\end{align}

Thus, for a given diffeomorphism complexity, the minimal perturbation can be expressed via the following additive decomposition of the total vector field perturbation (see Fig.~\ref{fig:pareto}B):
\begin{align}\label{eq:decomp:triple}
    \vg(\vx) &= 
	\vf(\vx)
	&\qquad\text{target system}
	\\
	&\quad+\vDelta_1(\vx)
	%&\qquad\text{min. required for topological conjugacy}
	&\qquad\text{imperfection}
	\\
	&\quad+\vDelta_2(\vx, c)
	&\qquad\text{covered by sub-optimal diffeomorphism}
	\\
	&\quad+\vDelta_3(\vx, c)
	&\qquad\text{remaining additive gap}
\end{align}

Figure~\ref{fig:pareto}B illustrates the interaction between the two distortion measures---we want the simplest diffeomorphism and the smallest perturbation, forming a Pareto frontier for the optimal deformation.
Point \textbf{a} represents the case explaining the full distortion solely using perturbation.
Point \textbf{c} represents only using the perturbation to explain the imperfection that cannot be captured by any diffeomorphism.
Point \textbf{b} represents a trade-off where the best diffeomorphism for the given complexity constraint creating additional gap for the perturbation.

Given the best diffeomorphism with the complexity constraint $c$, the trajectory loss $\bar{d}_1$ lower bounds the size of perturbation needed.
\begin{align}
    \bar{d}_1(\vf, \vg; c)
	&\leq
	    \left(\sup_\vs \norm{\vDelta_1(\vs) + \vDelta_3(\vs, c)}\right)
	    \frac{e^{LT_{\text{max}}}}{L^2}
	\leq
	    \left(\sup_\vs\norm{\vDelta_1(\vs)} + \sup_\vs\norm{\vDelta_3(\vs, c)}\right)
	    \frac{e^{LT_{\text{max}}}}{L^2}.
\end{align}
Since the $L$, $\norm{\vg}$, $\norm{\nabla\vg}$ are unknown in general and the Pareto frontier is difficult to compute, we propose to report the pair $(\bar{d}_1(\vf, \vg; c), c)$ in practice, to allow diffeomorphisms with $\vDelta_3 \neq 0$.
Reporting the pair, rather than a distance alone, is forced rather than merely convenient: on finite-horizon trajectory data a diffeomorphism of unbounded complexity can always achieve exact conjugacy, so the distance is informative only relative to a complexity budget (Supp.Sec.~\ref{sec:finite_time_conj}).

\begin{figure}[t!bhp]
    \centering
    \includegraphics[width=\linewidth]{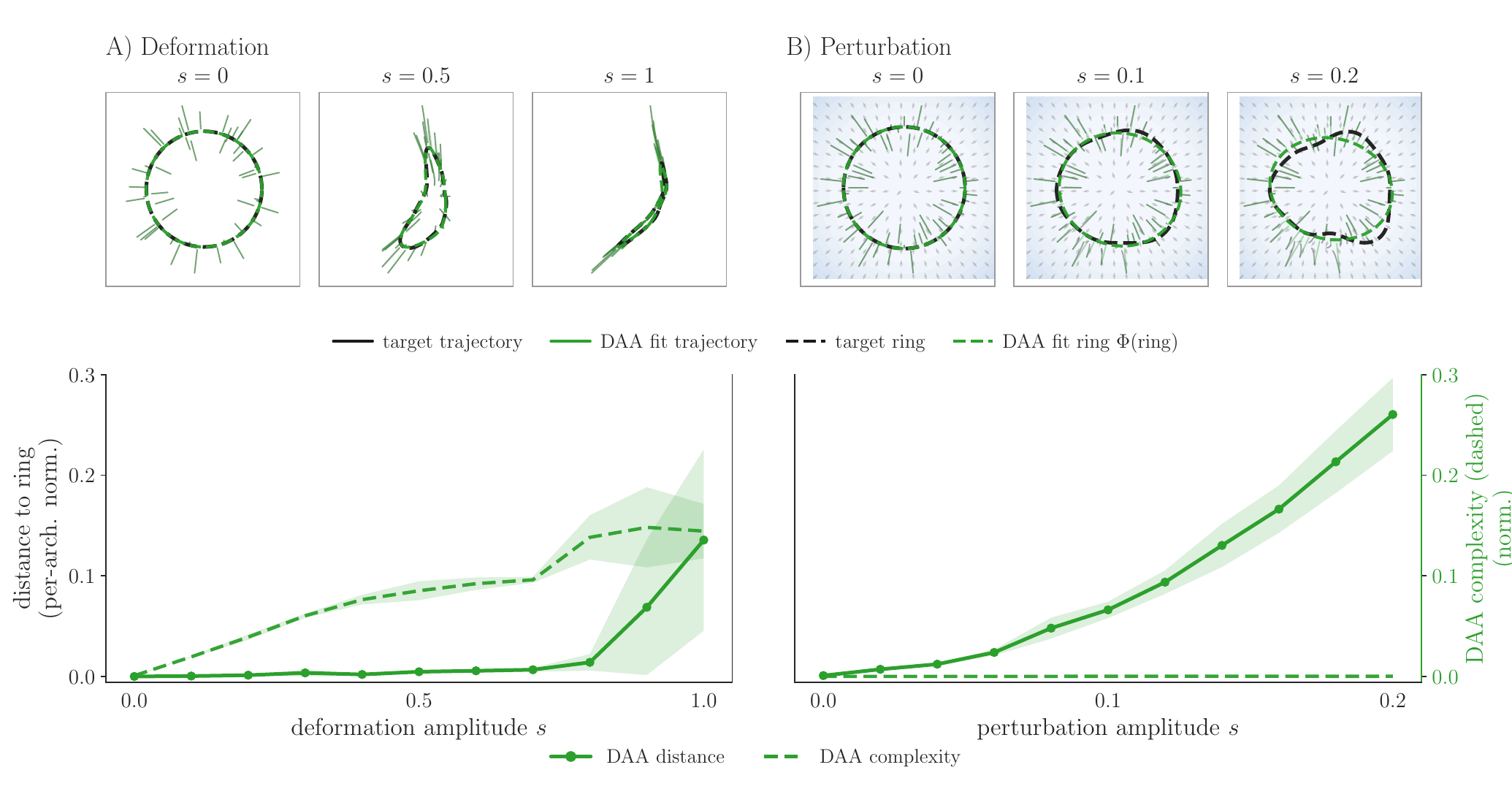}
    \caption{Dissimilarity of the \emph{ring attractor} archetype to an increasingly distorted ring, under DAA. In each column the top strip shows phase portraits at three amplitudes and the bottom panel the dissimilarity as $s$ increases. In the phase portraits the invariant ring is coloured by its angular coordinate on a perceptually uniform CIELAB wheel (so colour order tracks the ring's angular code); the target ring is drawn solid and the DAA fit ring $\homeo(\mathrm{ring})$ dashed, with trajectories in grey for context. In the bottom panels the target is DAA-fit and the dissimilarity is decomposed into trajectory distance (solid, left axis) and diffeomorphism complexity (dashed, right axis).
    (\textbf{A}) \emph{Deformation.} A diffeomorphic warp that leaves the continuum of fixed points intact, so the perturbed system remains topologically conjugate to the ring. DAA recovers the (inverse) diffeomorphism: the fitted ring tracks the target up to a large amplitude ($s\!\approx\!0.8$), the trajectory distance (solid, left axis) stays near zero, and the diffeomorphism complexity \eqref{eq:dcomplexity} (dashed, right axis) grows with the deformation.
    (\textbf{B}) \emph{Perturbation.} An additive vector-field perturbation that pins fixed points on the ring and immediately breaks topological conjugacy. The ring archetype still approximates the persisting slow invariant manifold --- the fitted ring tracks the target at small $s$ and departs from it as the pinned fixed points strengthen --- and the dissimilarity grows with $s$ because of the irreducible slack $\vDelta_1$ in \eqref{eq:decomp:triple}, while the diffeomorphism complexity stays low.
    }
    \label{fig:ring_pert_fig}
\end{figure}

% \subsubsection{The complexity of the homeomorphism}
% To assess how much a learned homeomorphism $\homeo$ deforms space, we distinguish between a set of canonical transformation types, such as translations, permutations, rotations, uniform and non-uniform scalings (of space), and rescaling of the vector field (which affects global speed).
%  While these operations are useful for interpreting specific deformation modes, they are not directly included in our complexity measure.
% Instead, we capture the full, potentially nonlinear deformation induced by $\homeo$ by evaluating how much it deviates locally from the identity map.
%  This is quantified by measuring how much the Jacobian of the homeomorphism deviates from the identity's:
% \begin{equation}\label{sec:jacobian_minid}
% \left\|\frac{\partial \homeo - \operatorname{id}}{\partial x}\right\|_{2,F} = \left\|\frac{\partial \homeo}{\partial x} - \mathbb{I}\right\|_{2,F},
% \end{equation}
% which captures how much the local linear behavior of $\homeo$ differs from that of the identity transformation.
% The norm of this deviation is computed either using the Frobenius norm or the spectral norm. %\ascomment{Frobenius or spectral?}
% Finally, the \( L^p \) norm of the batch of Jacobians is computed, providing a scalar value that quantifies the overall local deformation across the space. A norm close to zero indicates that the transformation is nearly identity, while higher values suggest more significant deformations.

\subsection{Empirical estimator: a trajectory-based loss function}\label{sec:loss}
Given a source system $\vg(\vx)$, in the scenario where we are analyzing neural trajectories from a target system without access to the underlying vector field $\vf(\cdot)$ nor its estimate, we need to estimate the dissimilarity vector and parameters based solely on a collection of trajectories.
Suppose the trajectories are uniformly sampled at $\Delta t$ intervals, that is, we are given data of the form:
$\mathcal{D} = \{ \vx_0, \vx_{\Delta t}, \vx_{2\Delta t}, \ldots, \vx_{n \Delta t}\}$.
The (mean squared error) loss function for the full trajectory is:
\begin{align}\label{eq:loss}
    L(\theta, \beta, \mathcal{D}) &\coloneq \frac{1}{n} \sum_{i=1}^n
	\norm{\vx_{i\Delta t} - \homeo_\theta(\phi_\vg^{i\Delta t}(\invhomeo_\theta(\vx_0); \beta))}^2
\end{align}
where $\theta$ represents the parameters for the diffeomorphism class $\homeo_\theta(\cdot)$ from $\vg$ space to $\vf$ space,
and $\beta$ are the parameters for a continuous family of $\vg(\vx; \beta)$, for example, the frequency for a standard limit cycle family.
Let $\theta^\ast, \beta^\ast$ be the minimizers of \eqref{eq:loss}.
$L(\theta^\ast, \beta^\ast, \mathcal{D}_{\text{train}})$ is 
an estimator for \eqref{eq:ed1}, and thus $\bar{d}_1(\vf, \vg; c = \dcomplexity(\homeo_{\theta^\ast}))$ using the squared 2-norm.
When the archetype's intrinsic dimension is smaller than the data dimension $D$, the source $\vg$ acts on $\reals^D$ as the archetype on its invariant manifold together with trivial contracting dynamics $\dot{\vx}_{\text{res}}=-\vx_{\text{res}}$ on the surplus directions; the choice of this residual block, and its dependence on the assumed regularity, is discussed in Supp.~Sec.~\ref{sec:trivial_dynamics}.

\section{Dynamical Archetype Analysis}\label{sec:dmm}\label{sec:daa}
Unlike typical clustering, where prototypes represent central or average examples, archetypes capture extreme, idealized forms that define the outer boundaries of a category. 
For example, the Platonic solids serve as idealized geometric forms---perfectly symmetrical, regular polyhedra that have long been regarded as the purest spatial shapes. 
While no physical object ever achieves this perfection, many natural and artificial forms resemble or approximate one of them.
In this sense, Platonic solids function as archetypes: extreme, abstract ideals that help us categorize and interpret the messy, imperfect forms observed in the real world.  In the same way, dynamical archetypes are discrete, canonical systems into which real neural dynamics can be categorized.

In neural computation, archetypes include theoretical constructs such as ``standardized'' fixed points, limit cycles, ring attractors, and continuous attractors\citep{sussillo2013blackbox, katz2017fibers, golub2018fixedpointfinder,townley2000existence, pals2024inferring,tian2025differential,liang2025symmetry}.
These idealized systems encapsulate essential computational roles in their purest and often most analytically tractable forms. For instance, an ideal ring attractor exhibits perfect rotational symmetry, free from distortion or drift. 
While biological and artificial neural circuits inevitably deviate from such perfection\footnote{They are non-repeating over time even for the same individual due to the ever changing brain\citep{chirimuuta2024brain}.}\citep{panichello2019error,gu2025attractor}, treating them as deformations of these archetypes offers a powerful interpretive framework.

Constructing a library of dynamical archetypes involves selecting a discrete set of theoretically grounded dynamical structures that cover a space of simple computations.
Rather than inferring archetypes directly from data, we define them a priori based on theoretical insights into neural computation.
Once established, this library serves as a reference vocabulary for analyzing observed dynamics by measuring their dissimilarity or deformation relative to these archetypes.
We call this approach \textbf{Dynamical Archetype Analysis} (DAA).

To make the archetype library interpretable and composable, we organize its elements by the topological dimension of their invariant manifolds.
This perspective reflects recent findings that task-trained RNNs often implement behavior through low-dimensional motifs---such as fixed points, limit cycles, and transient flows---that serve as functional subroutines \citep{driscoll2024flexible}. These motifs can be recombined across task conditions to support generalization \citep{tafazoli2024building}, and neural activity frequently lies on low-dimensional task manifolds embedded in state space\citep{langdon2023unifying,can2021emergence,cueva2021continuous,gort2024emergence,mishra2021continual,chaudhuri2019attractor,ghazizadeh2021slowmanifold,duncker2021dynamics, pezon2024linking,fortunato2024nonlinear}.
%Attractors can be further characterized by both their intrinsic and extrinsic dimensionality \citep{jazayeri2021interpreting,chaudhuri2019attractor,Sagodi2024a}. 

Based on these observations, we organize archetypes in a hierarchy based on the topological dimension of their invariant manifold and dynamics on it.
This classification provides a natural vocabulary for describing the dynamics we observe in neural systems (see Tab.~\ref{tab:archetypes} and Supp.~Sec.~\ref{sec:library}). 
The lowest-dimensional attractors in our framework are \emph{fixed points}, which correspond to isolated stable states --- typical of decision-making or memory models that settle into discrete outcomes~\citep{wong2006recurrent,bruno2017spiral}. 
Here, we consider two archetypes, the \emph{single fixed point} and the \emph{bistable system} (two stable fixed points, see Supp.~Sec.~\ref{sec:bistable}).
One-dimensional attractors in our library include the \emph{bounded line attractor} for integration tasks, the \emph{ring attractor} for circular variables, and the \emph{limit cycle} for periodic behaviors.
%At dimension 2, we find attractors like the plane and cylinder, representing combinations of continuous variables (e.g., two independent integrators, or a position and a phase).
 %More structured 2D attractors include the torus, formed by coupling two ring-like dimensions, and the sphere, a closed surface that supports more global integration or constrained exploration.
These archetypes serve as building blocks for constructing higher-dimensional attractors via Cartesian products, allowing for complex dynamics to arise from simple, interpretable components (see Supp. Sec.~\ref{sec:composite}).

%\mpcomment{Let's put a table of archetype and corresponding theoretical model or example memory manifold structure.}

\begin{table}[h]
\centering
\caption{Representative dynamical archetypes organized by dimension of the invariant manifold.}
\label{tab:archetypes}
\begin{tabular}{lll}
\toprule
\textbf{Dim.} & \textbf{Archetype} & \textbf{Memory manifold structure} \\
\midrule
0 & Single fixed point      & Discrete memory, decay to rest state \\
0 & Bistable system         & Binary decisions, winner-take-all dynamics \\
\midrule
1 & Bounded line attractor  & Evidence integration, continuous memory \\
1 & Ring attractor          & Head direction, orientation, phase coding \\
1 & Limit cycle             & Rhythmic activity, locomotion, oscillatory behavior \\
%\midrule
%2 & Torus                   & Coupled oscillations, grid cells, two phase-like variables \\
%2 & Sphere                  & 3D head direction representation \\
\bottomrule
\end{tabular}
\end{table}

\begin{figure}[t!b!hp]
    \centering
    \includegraphics[width=\linewidth]{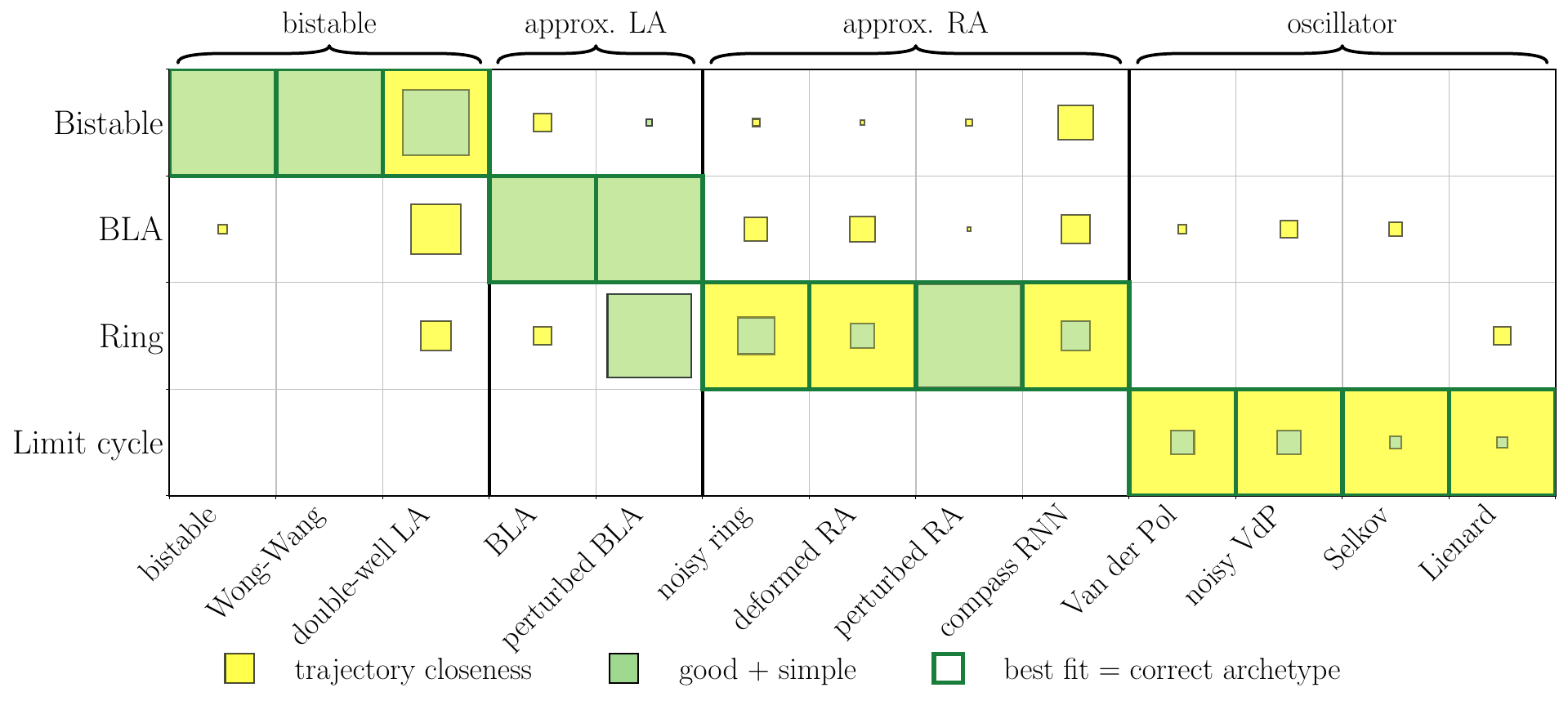}
    \caption{DAA identifies the archetype behind each target. Rows are the five archetypes (Stable FP, Ring, Limit cycle, BLA, Bistable), columns the target systems grouped by dynamical class (black rules; the three multi-member families are labelled above: approximate ring attractors, oscillators, bistable systems). Squares are scaled as described in Supp.Sec.~\ref{sec:squares_scaling}. The yellow square encodes trajectory closeness and is \emph{column-normalized}: within each target the sizes scale to that target's own range of archetype fits, so a target is compared only against its alternative explanations rather than across targets of different ambient dimension or fit difficulty. The green overlay encodes \emph{simplicity} of the required diffeomorphism (a near-identity map, small complexity \eqref{eq:dcomplexity}) on a fixed \emph{absolute} scale, so a green square of a given size means the same complexity in every column. A large green square is thus a good fit obtained with little deformation, whereas a yellow square with only a small green core is a good fit that demanded a substantial diffeomorphism.
Complexity is reported only where the fit is close enough for it to be meaningful, namely within a factor of three of the best distance in that column; elsewhere the square is left yellow. This matters because a mismatched archetype can reach a middling distance by letting the diffeomorphism run away, and the resulting complexity describes the optimiser rather than the target: fitting the BLA archetype to Selkov, for instance, reaches a distance $3.6\times$ the best available while paying a complexity of $33$, two orders of magnitude above any accepted fit. For every target the best-fitting archetype (smallest trajectory distance) is outlined in green when it matches the target's class. DAA assigns each target to its own archetype (all outlines correct), including the trained $64$-dimensional \emph{compass RNN} (an angular-velocity integrator, Fig.~\ref{fig:avi_rnn_recttanh}), which it identifies as a \emph{ring attractor}; the two-bounded-line-attractor target (``2 BLAs'') has no single matching archetype and is left unmarked.
    }
    \label{fig:archetype2target}
\end{figure}

%We consider parameterized versions of motifs that involve adding to on- or rescaling the vector field both on- and off-manifold.
Although attractors are categorized by their topological structure, the dynamics on those attractors can vary continuously within each category. For example, for the stable limit cycle, its dynamics are fully specified after choosing a scalar velocity $v\in \reals$ such that the dynamics are given by:
\begin{equation}
\dot{r} =  -r (r - 1), %r_0?\alpha
\quad \dot{\theta} = v,
\end{equation}
where \( r \) represents the radial coordinate and \( \theta \) the angular coordinate and the system exhibits stable behavior at \( r = 1 \).

\section{Results}
We apply DAA to a range of target systems to assess its ability to recover and classify canonical dynamical archetypes, highlighting how it captures underlying structure, distinguishes them, and quantifies deviations from idealized dynamics.

%justification of this  experiment
% How do deformed attractors look through our motif-based analysis?
%This experiment should shed light on the claim: We propose a library of canonical computations and a new measure of dissimilarity that allows us to group systems based on their effective behavior by explicitly considering both deformations that break the topological conjugacy as well as diffeomorphisms that preserve it.

%\subsection{Quantifying the deformation to a ring attractor}\label{sec:imp_ring}
\subsection{Faithful recovery of deformed ring attractors}\label{sec:imp_ring}
To evaluate the effectiveness of DAA, we investigate how deformed ring attractors deviate from the archetypal ring attractor.
The dissimilarity measure can account for two types of variation:
(i) topological deformations that preserve the qualitative structure via diffeomorphisms, and
(ii) perturbations that break topological conjugacy, altering the underlying computation.
%
%description
We apply two different deformations to the ring attractor:
(i) warped by diffeomorphic deformations (see Supp.~Sec.~\ref{sec:homeopert_exp_description}) and 
(ii) adding a small random perturbation to the ring attractor vector field to induce dynamics that increasingly break topological equivalence (see Supp.~Sec.~\ref{sec:vfpert_exp_description}).

\begin{figure}[tbhp]
    \centering
    \includegraphics[width=\linewidth]{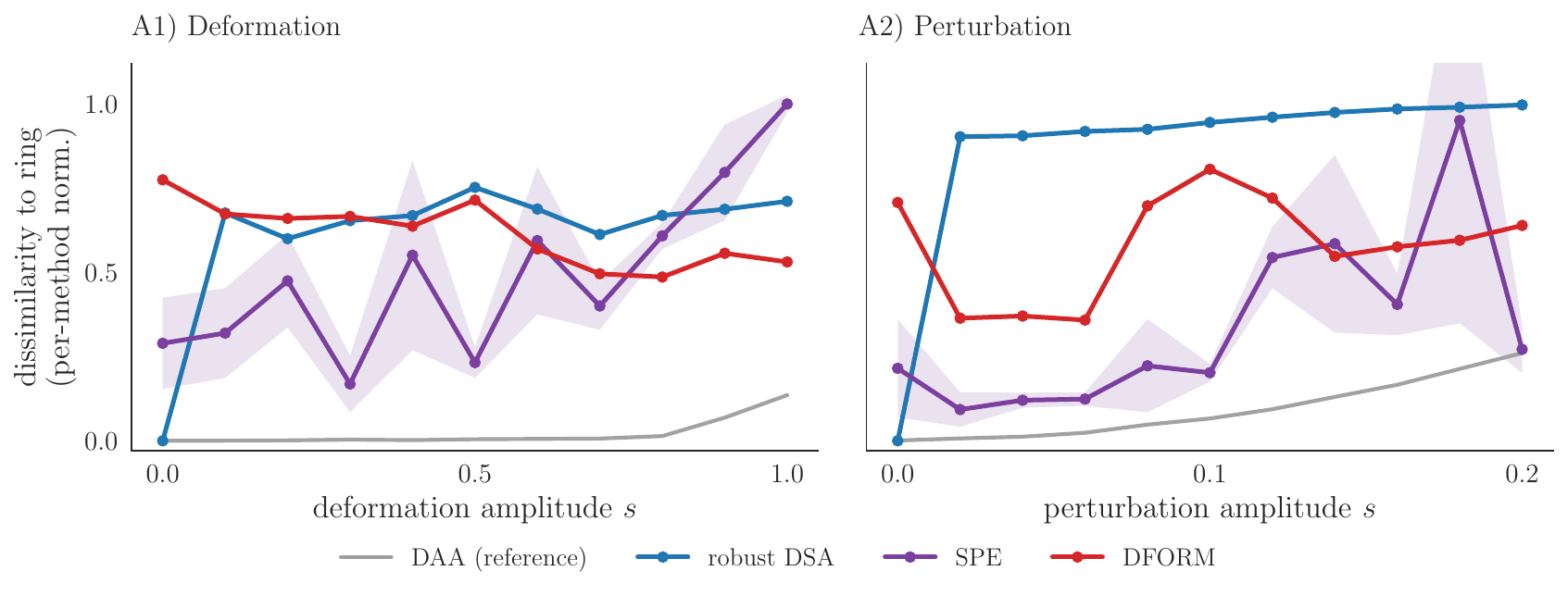}\\[4pt]
    \includegraphics[width=0.82\linewidth]{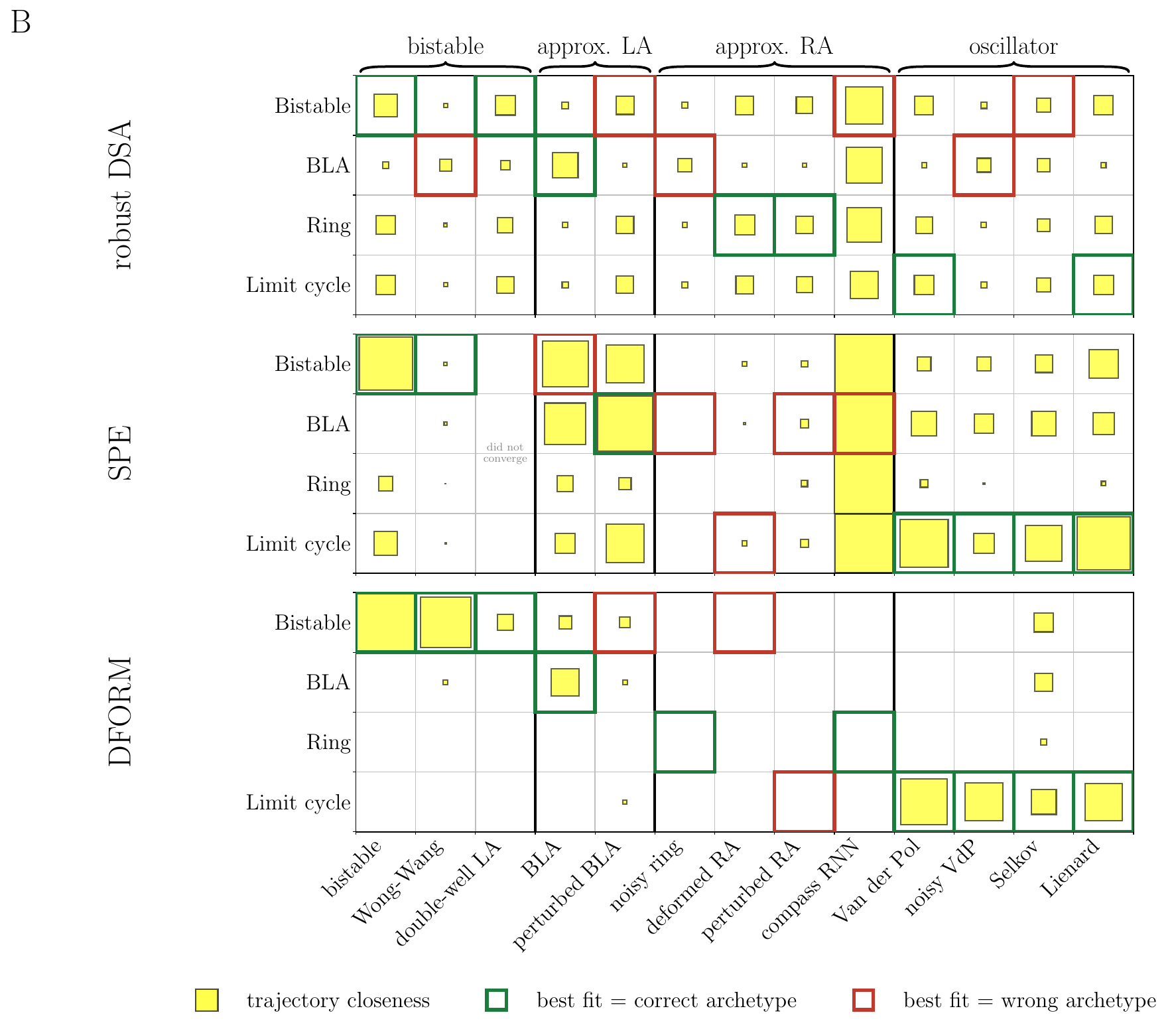}
    \caption{The three comparison methods (robust DSA, SPE, DFORM) on the same two analyses that DAA is applied to in Figs.~\ref{fig:ring_pert_fig} and~\ref{fig:archetype2target}.
    (\textbf{A}) Dissimilarity reported by each method as the ring is diffeomorphically deformed (\emph{Deformation}, left) or its vector field perturbed (\emph{Perturbation}, right), each normalized to its own range; DAA is drawn faint (grey) for reference. Unlike DAA, none of the three methods is invariant to the topology-preserving deformation: the reported dissimilarity is already large at small amplitude. Under the vector-field perturbation robust DSA saturates almost immediately while SPE and DFORM respond erratically, whereas DAA grows smoothly.
    (\textbf{B}) Archetype--target trajectory closeness for robust DSA, SPE and DFORM (rows), on the same absolute scale as Fig.~\ref{fig:archetype2target} (yellow only, since these methods carry no diffeomorphism-complexity measure); for each target that has a single matching archetype the best-fitting archetype is outlined green when it matches the target's class and red otherwise (the ambiguous ``2 BLAs'' target is left unmarked). Robust DSA and SPE misclassify about half of the targets, and DFORM, which is the strongest of the three here, still misclassifies a quarter of them and is, like the other two, not invariant to the diffeomorphic deformation in (\textbf{A}).
    }
    \label{fig:comparison3}
\end{figure}

We learn to map the trajectories of the \emph{ring attractor} archetype using a diffeomorphism \( \hat{\homeo}_\theta \), parameterized by a Neural ODE (NODE) \citep{chen2018neural} (for details see Supp.~Sec.~\ref{sec:architectures}) to the two approximate ring attractors (see Fig.~\ref{fig:ring_pert_fig}).
%interpretation of results
For the diffeomorphism-type deformation, our method approximates the inverse of the diffeomorphism, up to a maximal perturbation size after which the fit quality drops (Fig.~\ref{fig:ring_pert_fig}A).
For the vector field-type deformation, even without an exact topological conjugacy, the ring attractor archetype can be accurately mapped onto the approximate ring attractor (Fig.~\ref{fig:ring_pert_fig}B), while the dissimilarity grows because of $\vDelta_1$ in \eqref{eq:decomp:triple} and the complexity (\eqref{eq:dcomplexity}) grows with the deformation of the slow invariant manifold that persists after the ring attractor breaks apart.
For both deformation types, the location of the invariant manifold is accurately approximated (Fig.~\ref{fig:ring_pert_fig}A and B).
Running the three comparison methods on the same two analyses confirms that, unlike DAA, none of them is invariant to the topology-preserving deformation, and that robust DSA and SPE in particular misclassify many targets that DAA identifies correctly (Fig.~\ref{fig:comparison3}).

\subsection{Classifying the effective behavior of target systems}\label{sec:target_demonstrations}
To probe the ability of DAA to capture and distinguish canonical dynamical archetypes, we consider a range of target systems that vary in topology, noise, and deformation.
We now consider a library consisting of the \emph{ring attractor}, the \emph{limit cycle}, the \emph{single stable fixed point}, the \emph{bistable system} (two stable fixed points), and the \emph{bounded line attractor} (BLA).
% Can we find the Platonic ideal?
For a list of target systems, we try to fit the mold of each of the archetypes by fitting a diffeomorphism.
By comparing the different fits we can infer which of the archetypes best captures the effective behavior.
We test this by measuring the dissimilarity of our library of archetypes to a list of target systems:
1) Ring attractors (and its diffeomorphism- and additive vector field-deformed and noisy versions),
2) the \emph{compass RNN} trained on an angular-velocity integration task (with a slow ring invariant manifold, see Fig.~\ref{fig:avi_rnn_recttanh}A),
3) limit cycles (Van der Pol oscillator, Li\'enard sigmoid oscillator and Sel'kov oscillator)
4) system with two bounded line attractors (2 BLAs, to which no archetype aligns well).
See for the details of the target systems Supp.~Sec.~\ref{sec:target_systems} and the experimental setup Supp.~Sec.~\ref{sec:manytargets_exp_details}.

Normalized similarity pairs every target with the single archetype to which it is almost topologically conjugate via a small deformation (Fig.~\ref{fig:archetype2target}). 
The ring attractor-related targets all have high similarity to the \emph{ring attractor} archetype, while the oscillators have the highest similarity to the \emph{limit cycle} archetype (Fig.~\ref{fig:archetype2target}), while this is not captured by (robust) DSA, SPE or DFORM (Fig.~\ref{fig:comparison3}). 
The system with two bounded line attractors has the highest similarity to the \emph{bistable system} archetype; even though they are not topologically conjugate the bistability is shared.
%other way?If there is no canonical motif that is able to capture the effective behavior of a target, then all dissimilarities
Finally, if a target has a large dissimilarity to all the canonical archetypes in our library, that means that our library is not expressive enough to cover that computation with a corresponding Platonic archetype.
Alternatively, it is also possible that our parametrization of the diffeomorphism is not expressive enough ($\vDelta_2$ in \eqref{eq:decomp:triple}).

%\subsection{Comparison to other methods}\label{sec:compare_method}

%\paragraph{DSA}%https://github.com/mitchellostrow/DSA
%metric: DSA score
%%\paragraph{Phase2Vec} %https://github.com/nomoriel/phase2vec
%%\paragraph{DFORM} %https://github.com/rq-Chen/DFORM_stable
%\paragraph{SPE} %https://github.com/nitzanlab/prototype-equivalences
%metric: cycle-error

%\paragraph{SSMLearn} %https://github.com/haller-group/SSMLearnPy/tree/main

% \paragraph{MARBLE}
%MARBLE \citep{gosztolai2025marble} proposes a similarity measure based on local vector field features (vectors and higher-order derivatives) sampled from across the phase space.
% The features are embedded into lower-dimensional space through contrastive learning, and the embedded distributions from different systems can be compared using an optimal transport distance.

%%%%%%%%%%%%%%%%%
\section{Discussion}
As Popper put it, ``Science may be described as the art of systematic over-simplification -- the art of discerning what we may with advantage omit''~\citep{Popper1988}.
We propose to omit details between systems in the asymptotic time scale by characterizing the effective behavior of a neural system.
Although we have only analyzed autonomous computation (e.g., working memory) which forms the backbone for temporally extended computation, the approach can be extended to input-driven, or even embodied systems that interact with an environment.

The space of theories for neural computation corresponds to the space of possible dynamical systems the brain might instantiate---a vast landscape.
We linked the similarity of behavior in finite time with the deformation in the vector field (perturbation) and geometry (diffeomorphism).
DAA provides a means to chart regions of this space relevant for working memory maintenance by identifying and parametrizing recurrent computations in a human-interpretable way.

\textbf{Limitations:}
Compared to the state-of-the-art methods, our approach is computationally more demanding because it must learn the nonlinear diffeomorphism; however, this yields better identification accuracy and interpretability.
Another limitation is that the archetypes are not compositional in the present form, making it difficult to systematically expand the space of neural computation.

\makeatletter
\if@submission
\makeatother
\else
\section*{Acknowledgements}
% We would like to thank the anonymous reviewers from NeurIPS for their insightful comments and constructive feedback, which helped improve the quality of this work.
This project was supported by the Champalimaud Foundation and the Portuguese Recovery and Resilience Plan (PPR), through project number 62, Center for Responsible AI, and the Portuguese national funds, through FCT - Fundação para a Ciência e a Tecnologia - in the context of the project UIDB/04443/2020.
Part of this research was conducted while IMP was visiting the Okinawa Institute of Science and Technology (OIST) through the Theoretical Sciences Visiting Program (TSVP).
\fi

%This enables us to isolate families of dynamical motifs that serve as building blocks for memory-related neural dynamics.
%It accommodates both infinite-horizon systems and finite-horizon systems like echo state networks, highlighting that time-limited trajectories may obscure deeper structural distinctions only visible in the full vector field.

% though our method is computationally more demanding, yet brings better interpretability.

%\subsection{Limitations}
%Assumption: the invariant set of the target needs to be well-sampled.
%
%%Main shortcoming: difficulty with dealing with complex attractors.
%
%% When is decomposition of attractor into motifs possible?
%% The archetypes are not quite compositional in the present form.
%We focus on understandable attractors for now. 
%These are covered by our method to build higher-dimensional composite models from lower-dimensional dynamical motifs.
%
%% Input-driven computations would require an extension of the theory and method.
%
%%\subsection{}
%%Limits to neuroscientific understanding \citep{chirimuuta2024brain}
%%
%%simplest isn't always the best \citep{dyer2023simplest}
%%
%%planning: complex dynamical motifs that correspond to current and future actions \citep{vyas2020ctd}

\newpage
\bibliographystyle{unsrtnat_IMP_v1}
\bibliography{all_ref,catniplab}

\newpage
\appendix

%%%%%TOP
\section{Approximating diffeomorphisms}\label{sec:topology}
Our comparison represents each pair of systems by a learned diffeomorphism $\homeo$ that maps the phase space of a canonical archetype onto that of the target (Sec.~\ref{sec:daa}). This appendix collects the ingredients this construction relies on: the approximation-theoretic guarantees that a flow-based network can represent such a map (Supp.Sec.~\ref{sec:uap}), the concrete flow-based parameterizations we use for $\homeo$ (Supp.Sec.~\ref{sec:flowparam}), and the measure we use to quantify how far a fitted $\homeo$ departs from the identity, its complexity (Supp.Sec.~\ref{sec:complexity}).

\subsection{Universal approximation guarantees}\label{sec:uap}
Since a diffeomorphism is in particular a homeomorphism, the classical universal approximation results for continuous invertible maps apply directly to the maps we fit.
Any homeomorphism of $\mathbb{R}^n$ can be approximated by a Neural ODE or an i-ResNet operating in a lifted space of dimension $2n$; for the Neural ODE construction this follows from combining the Whitney embedding theorem for manifolds with the universal approximation theorem for neural networks, which gives an embedding dimension of $2n+1$ \citep{zhang2020approximation}.
Capping such a model with a single linear layer further turns it into a universal approximator for non-invertible continuous functions \citep{zhang2020approximation}.
These guarantees ensure that the flow-based parameterizations of the next section are, in principle, expressive enough to represent the archetype-to-target maps our method searches over.

%diffeomorphisms from gradient information of desired costs \citep{lai2021parallelised}

\subsection{Flow-based parameterizations}\label{sec:flowparam}
\paragraph{Neural ODEs}
Neural ODEs \citep{chen2018neural} represent the latest instance of continuous deep learning model for time series, first developed in the context of continuous recurrent networks \citep{cohen1983absolute}.
Neural ODEs model the continuous trajectory of the hidden state \( \mathbf{h}(t) \) as the solution to an ordinary differential equation parameterized by a neural network. Mathematically, the dynamics are defined by:
\[
\frac{d \mathbf{h}(t)}{dt} = f(\mathbf{h}(t), t, \theta)
\]
where \( \mathbf{h}(t) \) is the hidden state at time \( t \), \( f \) is a neural network with parameters \( \theta \), and \( t \) represents time.
Given an initial condition \( \mathbf{h}(t_0) = \mathbf{h}_0 \), the evolution of the hidden state over time is computed by solving the ODE using a numerical solver, allowing the model to capture complex, continuous transformations.

%example of f
\subsubsection{Feedforward Neural Network}\label{sec:fnn}
We will be using a Multi-layer Perceptron as our parametrized vector field. 
A single block is described by
\[
f(\mathbf{h}(t), t, \theta) = W_2 \sigma(W_1 \mathbf{h}(t) + \mathbf{b}_1) + \mathbf{b}_2,
\]
with \( \sigma \) is an activation function (e.g., ReLU, tanh, or sigmoid), \( W_1, W_2 \) are weight matrices, and \( \mathbf{b}_1, \mathbf{b}_2 \) are bias vectors.
This form models the ODE dynamics using a simple two-layer feedforward network.

\subsection{Measuring deformation via the Jacobian}\label{eq:jacobian_norm}
 For a smooth diffeomorphism \( \homeo \), the Jacobian \( J_{\homeo}(x) \in \mathbb{R}^{n \times n} \) is the matrix of partial derivatives:
\[
J_{\homeo}(x) = \left[ \frac{\partial \homeo_i}{\partial x_j}(x) \right]_{i,j}.
\]

\paragraph{Matrix norm of the Jacobian at a point \( x \).} Typically, we use either the operator norm (e.g., induced 2-norm) or the Frobenius norm:
\begin{equation}
\|J_{\homeo}(x)\| :=
\left\{
\begin{aligned}
\|J_{\homeo}(x)\|_{\mathrm{op}} &= \|J_{\homeo}(x)\|_2 = \sigma_{\max}(J_{\homeo}(x)) && \text{(spectral/operator norm)} \\
\|J_{\homeo}(x)\|_F &= \left( \sum_{i,j} \left| \frac{\partial \homeo_i}{\partial x_j}(x) \right|^2 \right)^{1/2} && \text{(Frobenius norm)}
\end{aligned}
\right.
\end{equation}

\subsubsection{Complexity: Measuring Deformation Relative to the Identity Map}\label{sec:complexity}
To quantify how much a diffeomorphism \( \homeo: \mathbb{R}^n \to \mathbb{R}^n \) deviates from the identity map, we analyze how much its Jacobian matrix \( J_{\homeo}(x) \) deviates from the identity matrix.
We capture the local deviation of the diffeomorphism from the identity map at point \( x \in \mathbb{R}^n \) by calculating:
\begin{equation}\label{eq:jac_id}
\|J_{\homeo}(x) - \mathbb{I}\|.
\end{equation}
When this norm is small, \( \homeo \) is locally close to the identity, meaning that it minimally stretches, rotates, or shears the space at that point.

To quantity measures the average deviation of the diffeomorphism from the identity map over a domain \( \Omega \subseteq \mathbb{R}^n \), with respect to a measure \( \mu \), we calculate: 
\[
\|J_{\homeo} - \mathbb{I}\|_{L^p(\mu)} := \left( \int_{\Omega} \|J_{\homeo}(x) - \mathbb{I}\|^p \, d\mu(x) \right)^{1/p}.
\]
 It reflects how much, on average, the transformation expands, contracts, or twists space away from the identity. 
 Unless specified otherwise, we use $p=2$.

Here, the domain \( \Omega \subseteq \mathbb{R}^n \) is taken to be the set of all trajectory points from a collection of dynamical systems. That is, \( \Omega \) consists of states 
\[
\{ \mathbf{x}_0^j, \mathbf{x}_{\Delta t}^J, \ldots, \mathbf{x}_{n \Delta t}^j \}
\]
for trajectories $j=1, \dots B$.
The measure \( \mu \)  corresponds to a uniform weighting over these points.
This norm captures the average deviation from the identity across the dynamic range of the system and is useful for regularizing global deformation induced by \( \homeo \).

%
%\newpage
%\section{Geometry}
%Geometric Learning on Manifolds\citep{mostowsky2024geometrickernels}
%
%
%\paragraph{Standard Mat\'ern Kernel (Based on Geodesic Distance)}
%Replace the Euclidean distance  $\| x - x' \|$  in the Mat\'ern kernel with the distance to the manifold  $ d_{\text{manifold}}(x, x')$. 
%This distance is now the measure of how far two points are from each other, but considering the geometry of the manifold.
%\begin{equation}
%k(x, x') = \frac{1}{\Gamma(\nu)} \left( \frac{\sqrt{2\nu} d_{\text{manifold}}(x, x')}{\ell} \right)^\nu K_\nu\left( \frac{\sqrt{2\nu} d_{\text{manifold}}(x, x')}{\ell} \right).
%\end{equation}
%
%
%\subsection{Implementations}
%\citep{miolane2020geomstats}
%%https://github.com/geomstats/geomstats?tab=readme-ov-file
%

%%%%%DS
\newpage
\section{Dynamical systems background}
This section recalls the equivalence relations on flows that our comparison rests on. We use two levels of strictness: topological conjugacy, a $C^0$ equivalence that preserves only the qualitative orbit structure, and smooth ($C^k$) equivalence, which additionally preserves the linearization and hence invariants such as eigenvalues. This distinction is what motivates reporting a match by two separate quantities, a topological distance and a smooth-deformation complexity.

%\subsection{Comparing dynamics}
%Binary/Discrete: equivalence
%
%Continuous comparison: metric/distance/dissimilarity

\subsection{Exact equivalence: Topological conjugacy}\label{sec:top_conj}

Topological conjugacy is an equivalence relation on the category of flows. %Behavioral space
\begin{definition}[Topological conjugacy]\label{def:top_conj}
let $\phi$ be a flow on $X$, and $\psi$ a flow on $Y$, with $X$, $Y$, and $h\colon Y \to X$ as above.

We say that $\phi$ is \emph{topologically semiconjugate} to $\psi$ if, by definition, $h$ is a surjection such that
\[
\phi(h(y), t) = h \circ \psi(y, t), \quad \text{for all } y \in Y, \; t \in \mathbb{R}.
\]
Furthermore, $\phi$ and $\psi$ are said to be \emph{topologically conjugate} if they are topologically semiconjugate and $h$ is a homeomorphism.
\end{definition}

Smooth equivalence is an equivalence relation in the category of smooth manifolds with vector fields. %model space
\begin{definition}[Smooth equivalence]\label{def:smooth_equivalence}
Two dynamical systems defined by the differential equations 
\[
\dot{x} = f(x) \quad \text{and} \quad \dot{y} = g(y)
\]
are said to be \emph{smoothly equivalent} if there exists a diffeomorphism \( h \colon X \to Y \) such that
\[
f(x) = M^{-1}(x) \, g(h(x)) \quad \text{where} \quad M(x) = \frac{d h(x)}{d x}.
\]
In that case, the dynamical systems can be transformed into each other by the coordinate transformation \( y = h(x) \).
\end{definition}

\begin{definition}[Orbital equivalence]\label{def:orb_eq}
Two dynamical systems on the same state space, defined by 
\[
\dot{x} = f(x) \quad \text{and} \quad \dot{x} = g(x),
\]
are said to be \emph{orbitally equivalent} if there exists a positive function \( \mu \colon X \to \mathbb{R} \) such that
\[
g(x) = \mu(x) f(x).
\]
\end{definition}
Orbitally equivalent systems differ only in their time parametrization.

%\paragraph{With scaling}
%Two dynamical systems on the same state space:
%\[
%\dot{x} = f(x) \quad \text{and} \quad \dot{x} = g(x).
%\]
%Let 
%\begin{itemize}
%\item $p\in C^1$
%\item $\homeo$ a homeomorphism
%\item $\mu\colon \reals^n\rightarrow \reals$ positive function
%\end{itemize}
%such that 
%\begin{enumerate}
%\item $\psi_{\mu \cdot (f+p)}$ is topologically conjugate to $\psi_{g}$
%\item Take $\inf_p \|p\|$, $\inf_\mu \|\mu\|$, $\inf_\homeo \|\homeo\|$
%%\item $f(x) = \mu(x) D\homeo(x)^{-1}(x) \, g(h(x))$ % is topologically conjugate
%\end{enumerate}

%\begin{definition}\label{def:top_conj}
%Let \( f: X \to X \) and \( g: Y \to Y \) be two continuous dynamical systems, where \( X \) and \( Y \) are topological spaces.
% The systems are said to be \emph{topologically conjugate} if there exists a homeomorphism \( h: X \to Y \) such that the following conjugacy condition holds:
%\[
%h \circ f = g \circ h.
%\]
%That is, the following diagram commutes:
%\begin{tikzcd}[row sep=large, column sep=large]
%X \arrow[r, "f"] \arrow[d, "h"'] & X \arrow[d, "h"] \\
%Y \arrow[r, "g"'] & Y
%\arrow[from=1-2, to=2-1, phantom, "\circlearrowleft", description, pos=0.5]
%\end{tikzcd}
%\end{definition}

%\paragraph{Topological conjugacy of time series}
%\citep{dlotko2024topconj}

\paragraph{Topological conjugacy on and off the invariant manifold}
When the homeomorphism from a topological conjugacy is restricted to the invariant set---such as an attractor or an $\omega$-limit set---it captures asymptotic dynamics faithfully, ensuring that trajectories on the archetype and the transformed system correspond one-to-one in their long-term behavior.
 This equivalence is robust under transformations that preserve topological structure, making it ideal for characterizing dynamics intrinsic to the invariant manifold.
 
 However, conjugacy on the invariant set does not guarantee equivalence in the transient dynamics that precede convergence. Outside the invariant set, trajectories may differ significantly even if the long-term behavior aligns. Such discrepancies in the off-manifold regions are crucial for biological systems, where transients often carry computational or functional relevance \citep{koch2024biological}. Understanding and controlling these off-manifold dynamics---while preserving topological conjugacy on the attractor---becomes essential when interpreting or constructing models that approximate natural archetypes.
 
% 
%\begin{itemize}
%\item On asymptotic / $\omega$-limit set / invariant set
%\item On transients \citep{koch2024biological}
%\end{itemize}

\subsection{Trivial off-manifold dynamics and the regularity of equivalence}
\label{sec:trivial_dynamics}
Every archetype in our library is low-dimensional (dimension $d$), whereas the observed data live in $\reals^D$ with $D \ge d$. To compare them we embed the archetype in $\reals^D$ by augmenting the $D-d$ \emph{residual} directions with trivial contracting dynamics,
\begin{align}\label{eq:residual}
  \dot{\vx}_{\text{res}} = -\vx_{\text{res}},
\end{align}
a global sink onto the $d$-dimensional archetype manifold. These residual directions are shared by all archetypes and carry no computation: their only role is to relax the state onto the computational manifold. It is therefore tempting to regard \eqref{eq:residual} as \emph{the} canonical trivial dynamics. The subtlety we address here is that which residual dynamics count as trivial, and hence what the single representative $\dot{\vx}_{\text{res}}=-\vx_{\text{res}}$ stands for, depends on the regularity of the equivalence used to compare systems.

Consider a hyperbolic contraction on $\reals^{D-d}$, i.e.\ a linear sink $\dot{\vy} = A\vy$ with $\mathrm{Re}\,\mathrm{spec}(A) < 0$.
\begin{itemize}[leftmargin=1.5em, noitemsep, topsep=2pt]
  \item \textbf{Topological ($C^0$).} By Hartman--Grobman \citep{hartman1960lemma,grobman1959homeomorphism} and the Kuiper--Ladis classification \citep{kuiper1975topology,ladis1973topological}, any two such sinks are topologically conjugate as soon as they share the same stable dimension. At this level a spiral sink, a fast node, a slow node, and $D-d$ independent $\dot x=-x$ modes are all the same: rates and rotations are invisible, and \eqref{eq:residual} is the unique trivial contraction.
  \item \textbf{Smooth ($C^k$, $k\ge 1$).} The full real Jordan form of $A$ is a conjugacy invariant \citep{sternberg1957local}. The representatives now split apart. In particular, a two-dimensional spiral sink with eigenvalues $-1 \pm bi$ is \emph{not} $C^k$-conjugate to two $\dot x=-x$ modes (eigenvalues $-1,-1$): the rotation rate $b$ is a smooth invariant that no diffeomorphism can remove. Differing decay rates ($-1$ versus $-2$) are likewise smoothly inequivalent, becoming equivalent only once a positive time rescaling is allowed ($C^k$-equivalence \citep{chen1963equivalence}), under which a spiral still never collapses to a real node.
  \item \textbf{Lipschitz.} An intermediate case: the imaginary parts at simple Jordan blocks are forgotten \citep{kawan2009lipschitz}, so a spiral and a node with matching real parts become equivalent again, while decay rates remain visible.
\end{itemize}
Thus \eqref{eq:residual} is canonically trivial only up to \emph{topological} equivalence; in the smooth category it is a genuine modeling commitment.

This matters for DAA because the dissimilarity \eqref{eq:loss} is parametrization-preserving (a conjugacy) and is evaluated in the smooth category, through the diffeomorphism $\homeo$. A target whose off-manifold relaxation spirals, or decays at a different rate, is therefore not matched exactly by the $-\identity$ padding, and the residual mismatch surfaces as trajectory loss and diffeomorphism complexity even though the off-manifold behavior is computationally irrelevant. This is the precise sense in which what counts as trivial dynamics depends on the assumed smoothness; it is the off-manifold counterpart of the on-manifold homeomorphism versus diffeomorphism distinction made in Sec.~\ref{sec:effective_ds}.

We adopt the following convention. The residual directions are compared only up to \emph{topological} equivalence: they encode fast collapse onto the computational manifold and are not part of the computation. Two facts make the specific representative \eqref{eq:residual} immaterial in practice. First, by separation of timescales the natural measure $\initDist$ concentrates near the manifold and the residual transients decay quickly relative to $T_{\text{max}}$, so their contribution to \eqref{eq:loss} is negligible provided the residual rate is fast enough. Second, when this fails the residual rate (and, if necessary, a rotation block) can be fit as nuisance parameters alongside the archetype parameters $\beta$, matching the residual block up to $C^k$-equivalence. We take $\dot{\vx}_{\text{res}}=-\vx_{\text{res}}$ as the canonical representative of the topologically unique trivial contraction class.

\paragraph{Diagonal versus block-diagonal residual generators.}
The residual generator in \eqref{eq:residual} is diagonal ($A_{\text{res}}=-\identity$, real spectrum). If the transverse relaxation of a target spirals (transverse eigenvalues $-a\pm i\omega$ with $\omega\neq 0$), no diffeomorphism conjugates the spiral to the non-rotating $-\identity$ block (smooth case above), so the rotation is charged either to the imperfection $\vDelta_1$ or to the complexity $c$, despite being computationally inert. A faithful alternative replaces $-\identity$ by a real-Jordan generator assembled from $2\times 2$ blocks $\left[\begin{smallmatrix}-a & -\omega\\ \omega & -a\end{smallmatrix}\right]$ with $(a,\omega)$ fit as nuisance parameters, matching the transverse dynamics up to $C^k$-conjugacy; equivalently, one may project the transverse directions out of \eqref{eq:loss} and \eqref{eq:dcomplexity}. Whether the diagonal block suffices in practice (it always recovers the correct archetype, and the quantitative readout is barely affected under sufficient timescale separation) or the $2\times 2$ block-diagonal generalization is required (when off-manifold rotation is non-negligible) remains open.

\subsection{Defective equivalence: Mostly Conjugate Dynamical Systems}
The framework of Mostly Conjugate Dynamical Systems relaxes strict topological conjugacy by allowing equivalence to hold only on a large-measure subset of the phase space.
Introduced by \citet{skufca2007relaxing, skufca2008mostlyconjugate}, this notion acknowledges that full conjugacy is often too rigid for practical applications, especially in complex or noisy systems. 
Instead, MCDS seeks a homeomorphism that maps trajectories accurately on a dense, dynamically relevant subset---often the global attractor of the system---while tolerating deviations in less critical regions.
%\citep{skufca2007relaxing, skufca2008mostlyconjugate, bollt2010comparing}
This relaxation provides a more flexible lens for comparing systems with similar qualitative dynamics but differing off-attractor behavior. 
\citet{bollt2010comparing} expanded on this idea by formalizing the statistical and geometric tools necessary to quantify ``closeness'' of systems beyond strict conjugacy. In the context of biological or engineered systems, this perspective is particularly useful: one can maintain fidelity to essential system dynamics while accommodating the variability inherent in real-world implementations.

\subsection{Conjugacy on finite-time data: topology is a cost, not an obstruction}\label{sec:finite_time_conj}
All of the equivalence relations above are properties of the \emph{asymptotic} dynamics: the invariants that obstruct a conjugacy (fixed points, periodic orbits, their eigenvalues, the topology of the attractor) are defined by the behavior as $t \to \infty$.
Any practical method, ours included, instead observes a finite family of trajectory segments over a finite horizon.
That distinction is not cosmetic, and it is what forces the two-axis formulation.

\begin{proposition}[Finite-time conjugacy is unobstructed]\label{prop:finite_time_conj}
Let $\vf$ and $\vg$ be Lipschitz vector fields on $\reals^n$ with $n \geq 2$, and let $\{\phi_\vf(\vx_i, t)\}_{i=1}^N$ and $\{\phi_\vg(\vy_i, t)\}_{i=1}^N$, $t \in [0,T]$, be two finite families of trajectory segments containing no equilibrium and no closed orbit.
Then there exists a diffeomorphism $\homeo$ of $\reals^n$ satisfying the conjugacy relation $\homeo(\phi_\vf(\vx_i,t)) = \phi_\vg(\homeo(\vx_i), t)$ for all $i$ and all $t \in [0,T]$; that is, $\vDelta_1 = 0$ \emph{on the observed data}.
\end{proposition}

\begin{proof}[Proof sketch]
By uniqueness of solutions distinct trajectories are disjoint, and a segment containing no equilibrium and no closed orbit is a compact embedded arc.
The two families are therefore embeddings of the same compact $1$-manifold (a disjoint union of $N$ arcs) into $\reals^n$.
Define $\homeo$ on the arcs by matching time-parameterized points, $\homeo(\phi_\vf(\vx_i,t)) := \phi_\vg(\vy_i,t)$, which is precisely the conjugacy relation restricted to the data.
Since the arcs are finitely many, compact and pairwise disjoint, they admit pairwise disjoint tubular neighborhoods; extend $\homeo$ over each tube as the arc map times a fiber map and glue to the identity outside with a partition of unity.
Equivalently, by the isotopy extension theorem any two embeddings of a compact $1$-manifold in $\reals^n$, $n \geq 2$, are ambient isotopic, and the ambient isotopy may be taken to realize the prescribed time-matching.
\end{proof}

Two remarks delimit the statement.
First, the hypothesis $n \geq 2$ is necessary: a diffeomorphism of $\reals$ is monotone, so in one dimension the \emph{order} of the segments along the line is an invariant that no diffeomorphism can permute.
Second, the excluded configurations are asymptotic, not generic.
A trajectory of a Lipschitz field cannot reach an equilibrium in finite time without violating backward uniqueness, and a segment closes into a periodic orbit only if its initial condition lies exactly on the cycle, a set of measure zero.
Generic finite-horizon data therefore consists of disjoint arcs, and the proposition applies.

The consequence is that the topological structure a method is asked to identify is never literally present in finite-horizon data; the data merely approaches it.
What survives is quantitative.
A trajectory that has spiralled close to a limit cycle is, mathematically, still an arc that may be unwound, but the diffeomorphism that unwinds it must shear by roughly the accumulated winding, and its cost diverges as the observation approaches the asymptotic regime.
For the rotation family of Supp.Sec.~\ref{sec:omega_matrix} this cost is explicit: the required map is the log-spiral twist, of complexity $\sim |\Delta\omega|\,\ln(1/r_{\min})$, which diverges exactly as the data approaches the fixed point.
We summarize this as: \emph{in finite time, topology is not an obstruction but a cost}.

This is why neither axis is informative alone, and why the complexity budget $c$ in $\bar{d}_1(\vf,\vg;c)$ is a necessity rather than a convenience.
On the distance side, Prop.~\ref{prop:finite_time_conj} shows that a diffeomorphism of unbounded complexity drives the trajectory distance to zero for \emph{every} archetype, so an unconstrained distance cannot discriminate.
On the complexity side, the conjugacy is not unique: if $\homeo$ conjugates $\vg$ to $\vf$, so does $\homeo \circ S$ for any $S$ commuting with the flow of $\vg$, so an unpenalized fit returns an arbitrary representative of that family and the complexity it reports is an upper bound of no particular significance (Supp.Sec.~\ref{sec:cplx_validation} measures this slack).
Constraining complexity resolves both defects at once, which is why we report the pair $(\bar{d}_1(\vf,\vg;c), c)$ and, in practice, fit with an explicit penalty on the complexity.
It also locates the identification claim correctly: an archetype is selected not because the alternatives are impossible, but because the correct one is \emph{cheap}.

\newpage
\section{Mathematical derivations}
This section collects the technical derivations supporting the main text. The central tool is a Grönwall bound on how far two flows diverge under a perturbation of the vector field; from it we obtain the estimates relating the reported dissimilarity to the size of the perturbation needed to align a target with an archetype.

\subsection{Bound on the flow divergence for a perturbation (\eqref{eq:gronwall})}\label{sec:flow_bound}
\begin{lemma}[Grönwall's inequality]
Let \( u(t) \) be a non-negative continuous function on \( [0, T] \), and suppose  
\[
u(t) \leq \alpha + L \int_0^t u(s)\, ds
\]
for constants \( \alpha \geq 0 \), \( L \geq 0 \). Then  
\[
u(t) \leq \alpha e^{Lt}.
\]
\end{lemma}

Let \( \phi^t_{\vf}(\vx) \) and \( \phi^t_{\vf + \vDelta}(\vx) \) denote the flows of the vector fields \( \vf \) and \( \vf + \vDelta \), respectively, starting from the same point \( \vx \). Assume \( \vf \) is Lipschitz with constant \( L \), and let \( \vDelta(t) \) be a time-dependent perturbation. Then the deviation between flows is bounded by
\[
\norm{ \phi^t_{\vf}(\vx) - \phi^t_{\vf + \vDelta}(\vx) }
\leq
e^{Lt} \int_0^{t} e^{-Ls} \norm{\vDelta(s)}\, \mathrm{d}s.
\]
This follows from Grönwall’s inequality applied to the differential inequality for the flow difference.

Using the fact that \( \norm{\vDelta(s)} \leq \sup_{s \in [0, t]} \norm{\vDelta(s)} \), we obtain the simplified upper bound:
\[
\norm{ \phi^t_{\vf}(\vx) - \phi^t_{\vf + \vDelta}(\vx) }
\leq
\left( \sup_{s \in [0, t]} \norm{\vDelta(s)} \right) \cdot \frac{1}{L}(e^{Lt} - 1).
\]
This shows that the trajectory deviation grows proportionally with the perturbation size and is amplified over time depending on the system's Lipschitz constant.

\subsubsection{A tighter bound for \eqref{eq:d0bound}}
The estimate \eqref{eq:d0bound} discards the $-1$ in $(e^{Lt}-1)$ before integrating over the horizon.
Retaining it and integrating the pointwise Gr\"onwall bound \eqref{eq:gronwall} directly yields a strictly sharper estimate.
With $\bar{d}_0(\vf,\vf+\vDelta)=\E_\vx\int_0^{T_\text{max}}\norm{\phi^t_\vf(\vx)-\phi^t_{\vf+\vDelta}(\vx)}\,\dm t$ and $\norm{\phi^t_\vf(\vx)-\phi^t_{\vf+\vDelta}(\vx)}\le\big(\sup_\vs\norm{\vDelta(\vs)}\big)\tfrac{1}{L}\big(e^{Lt}-1\big)$,
\begin{align}
    \bar{d}_0(\vf, \vf + \vDelta)
    &\leq
    \left(\sup_\vs \norm{\vDelta(\vs)}\right)\frac{1}{L}\int_0^{T_\text{max}}\!\big(e^{Lt}-1\big)\,\dm t
    =
    \left(\sup_\vs \norm{\vDelta(\vs)}\right)\left(\frac{1}{L^2}\big(e^{LT_\text{max}}-1\big) - \frac{1}{L}T_\text{max}\right)
    \notag\\
    &\leq
    \left(\sup_\vs \norm{\vDelta(\vs)}\right) \left(\frac{1}{L^2} e^{LT_{\text{max}}} - \frac{1}{L}T_{\text{max}}\right),
    \label{eq:d0bound_tight}
\end{align}
where the last step drops the negative constant $-1/L^2$.
The retained $-\tfrac1L T_\text{max}$ and the exact prefactor make \eqref{eq:d0bound_tight} smaller than \eqref{eq:d0bound}.

\subsubsection{Stochastic dynamics: a mean-square bound}\label{sec:stoch_bound}
For noisy single trials the target is better modelled as an It\^o SDE,
\begin{align}
    \dm{\vX_t} &= \big(\vg + \vDelta\big)(\vX_t)\,\dm t + \Sigma^{1/2}\,\dm{\vB_t},
    \qquad \vX_0=\vx,
\end{align}
with $\vg$ globally $L$-Lipschitz, $\vB_t$ a standard $d$-dimensional Brownian motion, and constant diffusion covariance $\Sigma=\Sigma^{1/2}(\Sigma^{1/2})\trp\succeq 0$.
Comparing the noisy trajectory to the deterministic archetype flow $\phi^t_\vg(\vx)$ and applying It\^o's lemma to $\norm{\vX_t-\phi^t_\vg(\vx)}^2$ (the martingale $2\langle\vX_t-\phi^t_\vg,\Sigma^{1/2}\dm{\vB_t}\rangle$ has zero mean; the second-order It\^o term contributes $\tr\Sigma$), then Cauchy--Schwarz, Young, and Gr\"onwall, gives
\begin{align}\label{eq:stoch_bound}
    \E\norm{\vX_t-\phi^t_\vg(\vx)}^2
    &\leq
    \Big(\underbrace{\sup_{s\in[0,t]}\E\norm{\vDelta(s)}^2}_{\text{model mismatch}}
    + \underbrace{\tr\Sigma}_{\text{noise floor}}\Big)\,\frac{e^{(2L+1)t}-1}{2L+1}.
\end{align}
The bound splits into a reducible mean-square drift mismatch and an \emph{irreducible noise floor} $\tr\Sigma$ that persists even for a perfect archetype ($\vDelta\equiv 0$): single-trial trajectory dissimilarity cannot fall below the diffusion floor, the theoretical counterpart of the empirically measured split-half noise ceiling.
The rate constant $2L+1$ arises from the cross term $2\langle\vX_t-\phi^t_\vg(\vx),\vDelta\rangle$ in $\tfrac{\dm{}}{\dm t}\E\norm{\vX_t-\phi^t_\vg(\vx)}^2$: bounding it with Young's inequality $2\langle\va,\vb\rangle\le\gamma\norm{\va}^2+\gamma^{-1}\norm{\vb}^2$, valid for any weight $\gamma>0$, contributes $\gamma$ to the growth rate and $\gamma^{-1}$ to the mismatch term. Taking $\gamma=1$ gives the displayed constant; in general the growth rate is $2L+\gamma$ and the mismatch prefactor $\gamma^{-1}$, so $\gamma$ trades exponential amplification against sensitivity to $\vDelta$. This is the standard mean-square (moment) stability estimate for It\^o SDEs~\citep[Ch.~4]{mao2007stochastic}.

\subsection{First-order bound and the Pareto frontier}\label{sec:fojac}
Let $\homeo$ be a diffeomorphism that can be used for achieving the smooth conjugacy in \eqref{eq:d1:tc}.
Now, we want to consider the complexity of $\homeo$ and trade-off with additive perturbation.
For this, we will do a first-order perturbation analysis.

To link the complexity of the diffeomorphism to the magnitude of an additive perturbation, we work out, to first order, the perturbation of the vector field that is equivalent to conjugating by a near-identity diffeomorphism. Consistent with \eqref{eq:ed1}, $\homeo$ carries the source $\vg$ to the target, so it is $\vg$ that is transformed; we seek the additive perturbation $\epsilon\vDelta_2$ of $\vg$ that its first-order action reproduces,
\begin{align}
    \homeo_\ast\vg &= \vg + \epsilon\vDelta_2 + O(\epsilon^2).
\end{align}

We assume the near-identity Ansatz with first-order generator $\tilde\homeo$ (the $O(\epsilon)$ displacement field), and record its Jacobian together with the inverse of the near-identity map:
\begin{align}\label{eq:appendix:h:ansatz}
    \homeo(\vx) &= \vx + \epsilon \tilde\homeo(\vx) + O(\epsilon^2),
    &
    \nabla\homeo(\vx) &= \identity + \epsilon\nabla\tilde\homeo(\vx) + O(\epsilon^2),
    &
    \homeo^{-1}(\vx) &= \vx - \epsilon\tilde\homeo(\vx) + O(\epsilon^2).
\end{align}
Transforming a vector field by a diffeomorphism is the \emph{pushforward}, which is defined with the inverse map, exactly the $\homeo^{-1}$ that \eqref{eq:ed1} applies to the state before flowing the source. At the point $\vx$,
\begin{align}
    (\homeo_\ast\vg)(\vx) &= \nabla\homeo\!\big(\homeo^{-1}(\vx)\big)\,\vg\!\big(\homeo^{-1}(\vx)\big).
\end{align}
Substituting the expansions (both factors evaluated at $\homeo^{-1}(\vx)=\vx-\epsilon\tilde\homeo(\vx)+O(\epsilon^2)$),
\begin{align}
    (\homeo_\ast\vg)(\vx)
    &= \big(\identity + \epsilon\nabla\tilde\homeo(\vx)\big)\Big[\vg(\vx) - \epsilon\,(\nabla\vg(\vx))\,\tilde\homeo(\vx)\Big] + O(\epsilon^2)
    \\
    &= \vg(\vx) + \epsilon\big[(\nabla\tilde\homeo(\vx))\,\vg(\vx) - (\nabla\vg(\vx))\,\tilde\homeo(\vx)\big] + O(\epsilon^2).
\end{align}
Matching identifies the diffeomorphism-induced perturbation with the (Jacobian) Lie bracket $[\va,\vb]\coloneqq(\nabla\vb)\,\va-(\nabla\va)\,\vb$ of $\vg$ and the generator $\tilde\homeo$,
\begin{align}\label{eq:pareto:delta}
    \vDelta_2 &= (\nabla\tilde\homeo)\,\vg - (\nabla\vg)\,\tilde\homeo \;=\; [\vg,\tilde\homeo] \;+\; O(\epsilon).
\end{align}
The two terms come from the two ways $\homeo$ enters the pushforward: the outer Jacobian $\nabla\tilde\homeo$ scaling $\vg$, and the argument shift $\homeo^{-1}$ (evaluating $\vg$ at the back-mapped point) through $\nabla\vg$.

The full perturbation carrying the source into the target is $\vDelta=\vDelta_1+\vDelta_2$, with $\vDelta_1$ the irreducible slack of \eqref{eq:d1:tc} (the minimal perturbation into the conjugacy class, which no diffeomorphism removes). With \eqref{eq:pareto:delta} and the triangle inequality, for all $\vx$,
\begin{align}
    \norm{\vDelta(\vx)}
    &\le \norm{\vDelta_1(\vx)} + \norm{(\nabla\tilde\homeo(\vx))\,\vg(\vx)} + \norm{(\nabla\vg(\vx))\,\tilde\homeo(\vx)} + O(\epsilon)
    \notag\\
    &\le \norm{\vDelta_1(\vx)} + \norm{\vg(\vx)}\cdot\underbrace{\norm{\nabla\tilde\homeo(\vx)}}_{\substack{\text{scale c.}}}
    + \norm{\nabla\vg(\vx)}\cdot\underbrace{\norm{\tilde\homeo(\vx)}}_{\substack{\text{translation c.}}} + O(\epsilon),
\end{align}
for appropriate vector and corresponding matrix norms.
We recognize two kinds of complexities of the diffeomorphism: (1) scale complexity which is a function of the Jacobian of diffeomorphism $\nabla\homeo(\vx)$, and (2) translation complexity.

Under the assumption that $\norm{\vg(\vx)}$ and $\norm{\nabla\vg(\vx)}$ are bounded within the region of interest, we argue that a linear combination of two complexity measures of the diffeomorphism can be traded off with the additive perturbation.

\subsubsection{Integrated bound}
The bound \eqref{eq:h:bound} holds point-wise.
Given a distribution over trajectories (orbits of finite duration),
$p(\phi(\vx_0,T_\text{max}))$, we can integrate the inequality over the marginal probability measure $p(\vx)$:
\begin{align}\label{eq:h:bound:integrated}
    \E_\vx \norm{\vDelta(\vx)}
    &\leq
	\E_\vx \norm{\vDelta_1(\vx)}
	+
	\E_\vx
	\left[
	    \norm{\vg(\vx)}
	    \cdot
	    \norm{\nabla\tilde\homeo(\vx)}
	\right]
	+
	\E_\vx
	\left[
	    \norm{\nabla\vg(\vx)}
	    \cdot
	    \norm{\tilde\homeo(\vx)}
	\right]
	+ O(\epsilon).
\end{align}
The imperfection term $\E_\vx\norm{\vDelta_1(\vx)}$ is carried over from the point-wise bound \eqref{eq:h:bound}: it is the irreducible slack that no diffeomorphism can remove, and it vanishes exactly when $\vf \tconju \vg$.

\section{What the complexity measures}\label{sec:cplx_validation}
Sec.~\ref{sec:trivial_dynamics} establishes what a diffeomorphism can and cannot remove: for hyperbolic sinks the stable dimension is the only topological ($C^0$) invariant, while the real Jordan form---in particular a rotation rate---is a smooth ($C^k$) invariant that no diffeomorphism can absorb, so a rotation must be charged \emph{either} to the imperfection $\vDelta_1$ \emph{or} to the complexity $c$.
This appendix settles that alternative empirically, and then asks what the complexity number is actually made of.

\subsection{The rotation matrix: complexity charges the deformation that is required}\label{sec:omega_matrix}
Consider the one-parameter family of planar linear systems
\begin{align}\label{eq:omega_family}
    \dot{\vx} = A(\omega)\,\vx,
    \qquad
    A(\omega) = \begin{bmatrix} -a & -\omega \\ \omega & -a\end{bmatrix},
    \qquad a = 1 .
\end{align}
Its eigenvalues are $-a \pm i\omega$, so $\operatorname{Re}\lambda = -a < 0$ for every $\omega$: the origin is a hyperbolic sink throughout and nothing crosses the imaginary axis.
By Sec.~\ref{sec:trivial_dynamics} the whole family is therefore a \emph{single} topological class, while $-a\pm i\omega$ is a $C^k$ invariant, so no two members with $\omega_1\neq\omega_2$ are smoothly conjugate.
(At $\omega=0$ the family crosses the node/focus boundary $\operatorname{tr}^2-4\det = -4\omega^2 = 0$; this changes the phase portrait's appearance but not its topological type, and is not a bifurcation.)
The family is thus an unusually clean test bed: \textbf{distance is blind by construction}, so if the complexity is meaningful it must carry the entire $\omega$-dependence.
The conjugacy between two foci is the log-spiral twist $\theta \mapsto \theta + \tfrac{\omega_2-\omega_1}{a}\ln r$, which predicts $\dcomplexity(\homeo) \sim |\omega_2-\omega_1|\cdot\ln(1/r_{\min})$: complexity should track the \emph{difference} of rotations, not either rotation itself.

We fit every archetype $A(\omega_a)$ to every target $A(\omega_t)$, $\omega_a,\omega_t\in\{0,0.2,\dots,1.0\}$.
The archetype's rotation must be \emph{frozen} for this to mean anything: a learnable full $A$ would simply become the target and the matrix would collapse to zero.
We therefore use the fixed-rotation archetype of Sec.~\ref{sec:lds_w} (zero learnable parameters).
The same observation is what makes the $\omega_a=0$ row---the node~$\to$~spiral sweep of Fig.~\ref{fig:ring_pert_fig}---a valid experiment: the default linear archetype is diagonal, has real eigenvalues, and \emph{cannot} rotate at all.

\begin{table}[htbp]
\centering
\caption{Complexity $\dcomplexity(\homeo)$ for every (archetype $\omega_a$, target $\omega_t$) pair; $50$ trajectories from an annulus, $T_{\text{max}}=2.5$. The trajectory distance is $\leq 1.6\times10^{-3}$ in \emph{all} $36$ cells and is omitted. Low along the diagonal, growing away from it in both directions.}
\label{tab:omega_matrix}
\begin{tabular}{c|cccccc}
\toprule
$\omega_a\backslash\omega_t$ & $0.0$ & $0.2$ & $0.4$ & $0.6$ & $0.8$ & $1.0$ \\
\midrule
$0.0$ & \textbf{0.074} & 1.408 & 1.993 & 2.238 & 2.627 & 1.432 \\
$0.2$ & 1.046 & \textbf{0.137} & 1.259 & 1.887 & 2.459 & 2.510 \\
$0.4$ & 1.936 & 1.138 & \textbf{0.134} & 0.970 & 1.794 & 2.242 \\
$0.6$ & 2.315 & 1.815 & 1.274 & \textbf{0.148} & 0.837 & 1.633 \\
$0.8$ & 1.767 & 2.141 & 1.790 & 1.285 & \textbf{0.193} & 0.772 \\
$1.0$ & 1.084 & 1.669 & 1.966 & 1.644 & 1.203 & \textbf{0.244} \\
\bottomrule
\end{tabular}
\end{table}

Three results (Tab.~\ref{tab:omega_matrix}).
First, the \textbf{trajectory distance is $\leq 1.6\times 10^{-3}$ in all $36$ cells}: topological conjugacy holds across the entire family, exactly as Sec.~\ref{sec:trivial_dynamics} predicts, and $\vDelta_1 \approx 0$.
This resolves the alternative left open there: the rotation is charged \emph{entirely to the complexity}, not to the imperfection.
Second, \textbf{complexity tracks the required deformation and not the target}: $\mathrm{corr}(\dcomplexity(\homeo),\,|\omega_t-\omega_a|) = +0.74$ whereas $\mathrm{corr}(\dcomplexity(\homeo),\,\omega_t) = +0.06$.
Had the second correlation been the large one, the measure would merely be reporting a property of the target rather than the cost of the map, and reading complexity as ``the deformation the archetype requires'' would be unsupported everywhere else in this work.
Third, the \textbf{diagonal is the control}: $\dcomplexity(\homeo) \in [0.07, 0.24]$ when $\omega_a=\omega_t$, an order of magnitude below the off-diagonal entries---$\homeo$ collapses to the identity precisely when the archetype already matches.

Two honest caveats follow from the same source.
The matrix is \emph{not} symmetric ($\mathrm{mean}|c_{ij}-c_{ji}| = 0.30$, i.e.\ $17\%$ of the typical entry) even though the twist $\pm(\Delta\omega/a)\ln r$ has the same magnitude in both directions, because \eqref{eq:dcomplexity} is evaluated on the \emph{target's} trajectory support, and a node and a focus populate the plane differently.
Complexity is therefore directional---a property of the triple (archetype, target, support)---and the matrix should not be read as a distance matrix on systems.
The same support-dependence appears in the $\ln(1/r_{\min})$ factor: the complexity of the node~$\to$~spiral conjugacy \emph{diverges} as the data approach the fixed point, and is finite here only because the initial conditions lie on an annulus ($r_{\min}\approx 0.07$).
Finally, the extreme corners ($|\Delta\omega| = 1$) carry both the largest distance and a \emph{depressed} complexity, the signature of an under-converged fit on the hardest twist; the reported $+0.74$ is therefore a lower bound on the true association.

\begin{figure}[htbp]
    \centering
    \includegraphics[width=\linewidth]{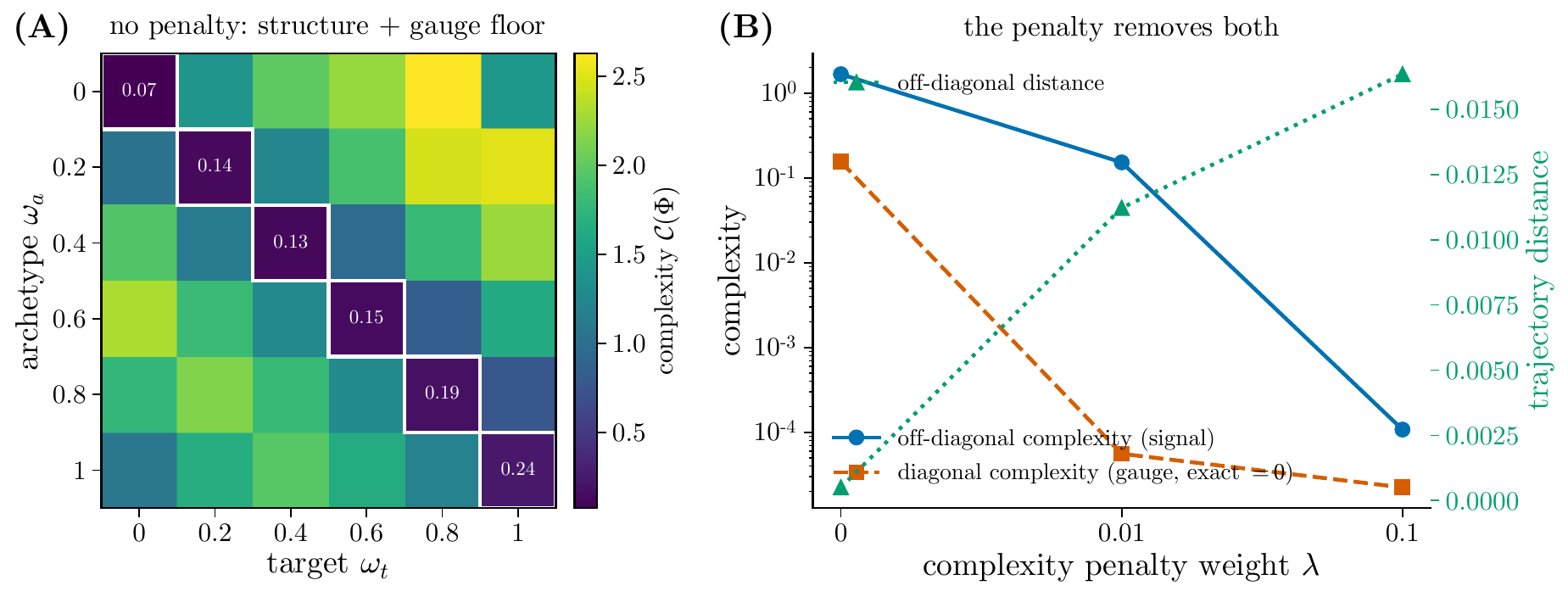}
    \caption{\textbf{The complexity penalty removes the gauge and the signal together.}
    \textbf{(A)} The unpenalized complexity matrix of Tab.~\ref{tab:omega_matrix}: entries grow away from the diagonal with $|\omega_t - \omega_a|$, and the diagonal (outlined, values inset) is the control, where the exact complexity is $0$ but the fit still reports $0.07$ to $0.24$.
    \textbf{(B)} Refitting at penalty weight $\lambda$. The diagonal (gauge, exact value $0$) is cleaned as intended, but the off-diagonal signal collapses with it, by four orders of magnitude at $\lambda = 0.1$, while the off-diagonal trajectory distance rises (right axis). Note the logarithmic complexity axis.
    }
    \label{fig:omega_lambda}
\end{figure}

\paragraph{Turning the penalty on.}
The diagonal of Tab.~\ref{tab:omega_matrix} is a pure gauge artifact: there the archetype \emph{is} the target, the identity is an exact conjugacy, and the correct complexity is $0$, yet the unpenalized fit reports $0.07$ to $0.24$, growing with $\omega$ because a larger rotation offers a larger commuting rotation to absorb.
Refitting the whole matrix with a complexity penalty of weight $\lambda$ removes it, and also removes everything else (Fig.~\ref{fig:omega_lambda}).
The diagonal mean falls from $0.155$ to $6\times10^{-5}$ and $2\times10^{-5}$ at $\lambda = 0.01$ and $0.1$, which is the intended effect.
But the off-diagonal mean, which is the signal, falls from $1.671$ to $0.152$ and then to $1.1\times10^{-4}$ over the same range, while the off-diagonal trajectory distance degrades by factors of $22$ and $33$.
At $\lambda = 0.1$ the entire matrix is numerically zero, so the correlation with $|\Delta\omega|$ that nominally survives there ($+0.72$) is a correlation among values of order $10^{-4}$: the penalty preserves the rank ordering while destroying the scale, which is of no use for a quantity whose purpose is to be compared against an analytic value.
Cold-started regularization also fails outright on the hardest cells.
At $\lambda = 0.01$ the fit for $\omega_a = 1 \to \omega_t = 0$ does not converge at all: the loss remains flat while the Jacobian grows monotonically until the integrator diverges, and reducing the solver step twice does not change the outcome (this cell is left blank; $35$ of $36$ converged).
The opposite corner converges only with a reduced step, and then to a distance five times the typical cell.
Both are the direction in which the map must \emph{un}-wind a rotation, which is precisely the growth the penalty opposes.
The conclusion is that a penalty applied from a cold start cannot separate removing the gauge from destroying the signal, because from a near-identity initialization the fit never builds the required twist in the first place.
Warm-starting from the unpenalized solution and annealing the penalty upward does separate them, recovering the analytic complexity to $0.94\times$ and $1.00\times$ on the two cases where it is known, with the trajectory distance \emph{improving} rather than degrading.

\subsection{Alternative complexity measures: $L^p$, matrix norm, and rank}\label{sec:cplx_alternatives}
Our complexity \eqref{eq:dcomplexity} aggregates a \emph{per-point matrix norm} of $\nabla\homeo - \identity$ over the trajectory support with an \emph{$L^p$ average}; we use the Frobenius norm and $p=2$.
Both are choices, and Sec.~\ref{sec:trivial_dynamics}'s lesson---that what counts as equivalent depends on the category one works in---has an exact counterpart here: what counts as \emph{complex} depends on the norm one picks.
The $p=2$/Frobenius pair is not arbitrary: it is the discrete counterpart of the $H^1$ seminorm, which is what links \eqref{eq:dcomplexity} to the LDDMM warp energy.
Alternatives depart from that correspondence and should be motivated, not defaulted into: $p\to\infty$ reports the \emph{worst} point (a Lipschitz reading, which would expose the $\ln(1/r_{\min})$ divergence above that $p=2$ dilutes), while $p=1$ reports the mean.

A more consequential choice is the per-point norm, because $\norm{\cdot}_F=(\sum_k\sigma_k^2)^{1/2}$ conflates deformations of different \emph{kind}: a rank-one warp of size $s$ and a full-rank warp of the same total size receive the same score.
The spectrum of $\nabla\homeo-\identity$ distinguishes them.
A convenient summary that requires no refitting is the \emph{stable rank} $\norm{\cdot}_F^2/\norm{\cdot}_2^2 = \sum_k\sigma_k^2/\sigma_{\max}^2 \in [1, D]$, since both norms are already reported.
On the family \eqref{eq:omega_family} the stable rank is $1.39$ on average (range $[1.16, 1.69]$ over the off-diagonal cells of Tab.~\ref{tab:omega_matrix}) and is essentially \emph{uncorrelated} with $|\Delta\omega|$ ($r = -0.16$), while the Frobenius size varies over a factor of $35$ across the same cells.
Size and shape are thus decoupled: increasing the demanded rotation makes the warp \emph{larger} but not different in kind.
The value itself is interpretable---close to $1$ in a space where the stable rank is at most $2$---and says the conjugacy is predominantly a single-direction shear, which is precisely the analytic log-spiral twist.
This suggests reporting complexity as at least a two-dimensional object (a size and a shape); we leave a systematic study, including the differentiable low-rank surrogates (nuclear norm, participation ratio) that would let one \emph{prefer} low-rank warps, to future work.

\subsection{Parameterized sweeps against the exact conjugacy}\label{sec:decomp_sweeps}
The rotation matrix of Supp.Sec.~\ref{sec:omega_matrix} varies one invariant at a time.
Here we do the same along three one-parameter families chosen so that the \emph{exact} conjugacy, and therefore the complexity it requires, is available in closed form.
This turns the complexity axis into a quantity with a known correct answer rather than a number to be interpreted.

\paragraph{The three families.}
In \textbf{(A) node $\to$ spiral} the archetype is a stable node and the target is the planar linear system with matrix $\begin{psmallmatrix} -a & -\omega \\ \omega & -a\end{psmallmatrix}$, so the swept parameter is the rotation rate $\omega$.
Every member is a hyperbolic sink, so all are topologically conjugate and the distance is blind to $\omega$ \emph{by construction}; the entire $\omega$-dependence must therefore be carried by the complexity.
The conjugacy is the log-spiral twist $\homeo(r,\theta) = (r, \theta - (\omega/a)\ln r)$, whose complexity we evaluate on the observed support.
In \textbf{(B) limit cycle with non-uniform speed} the on-cycle dynamics are $\dot\theta = -c(s)\,(1 + s\cos\theta)$, with $c(s)$ chosen to hold the period fixed (the period is a conjugacy invariant that a spatial diffeomorphism cannot change), and the conjugacy solves the circle equation $dg/d\theta = 1 + s\cos g$.
In \textbf{(C) ring attractor with process noise} the target \emph{is} the ring attractor, driven by additive noise of standard deviation $\sigma$.
No deformation whatsoever is required, so the exact complexity is $0$ at every $\sigma$ and any value the fit reports is noise-fitting.
All three sweeps are fitted without any penalty on the complexity, so the value reported is what an unregularized fit returns.

\begin{figure}[htbp]
    \centering
    \includegraphics[width=\linewidth]{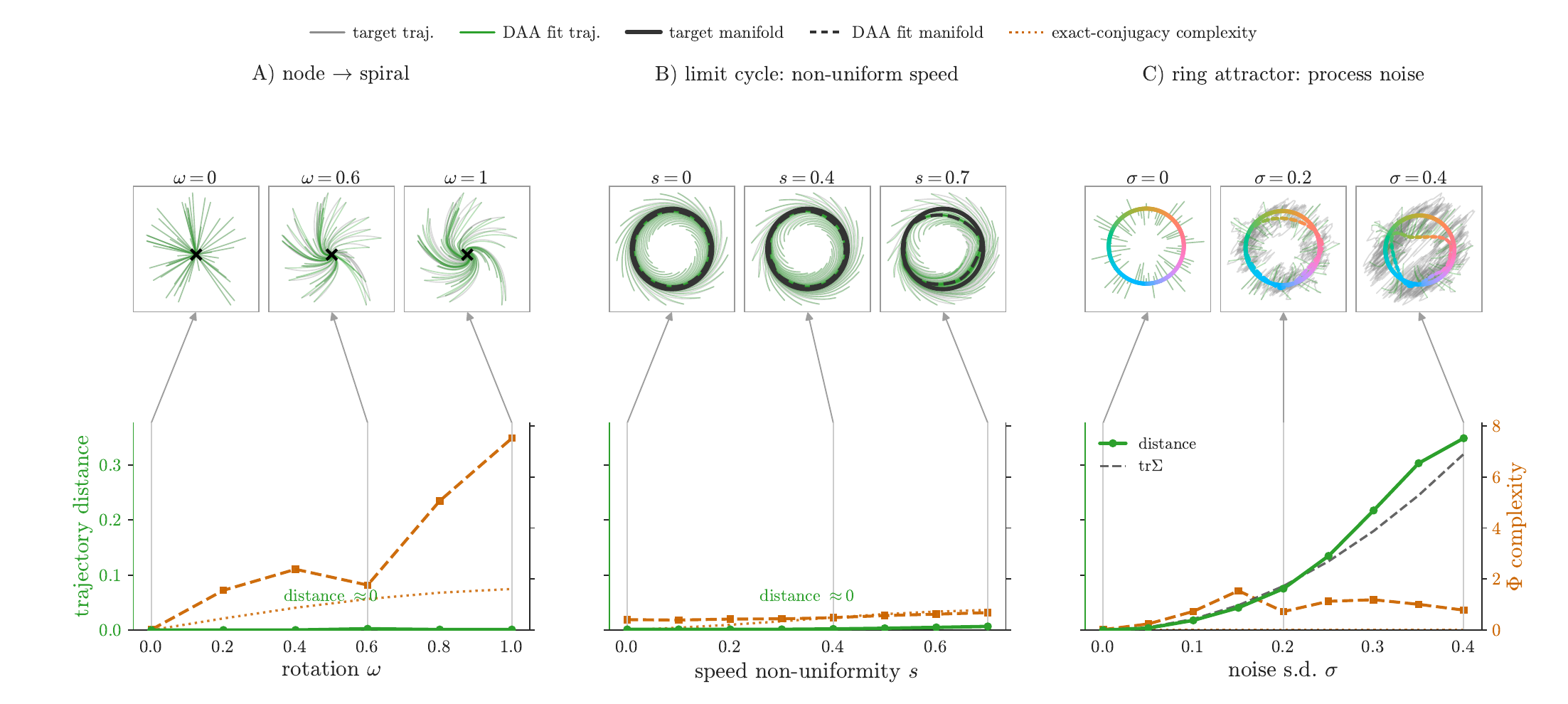}
    \caption{\textbf{Complexity measured against a known answer.}
    Three one-parameter families in which the exact conjugacy is known analytically (Supp.Sec.~\ref{sec:decomp_sweeps}).
    \emph{Top of each column:} target trajectories (grey) with the fitted archetype trajectories (green) and the target and fitted invariant manifolds (solid and dashed black), at three values of the swept parameter.
    \emph{Bottom:} trajectory distance (green, left axis) and diffeomorphism complexity \eqref{eq:dcomplexity} (orange, right axis) against the swept parameter, with the \emph{exact-conjugacy} complexity as the dotted orange curve.
    \textbf{(A)} Distance stays at $\approx 0$ for every $\omega$, as topological conjugacy requires, so complexity carries the whole effect. It rises with $\omega$ but overshoots the exact value by $1.5$ to $4.7\times$ and does so erratically.
    \textbf{(B)} At $s=0$ the target \emph{is} the archetype, the identity is an exact conjugacy and the exact complexity is $0$, yet the fit reports $0.41$. Once the required complexity exceeds this floor ($s \gtrsim 0.4$) the fit tracks the exact curve to within $0.79$--$0.99\times$.
    \textbf{(C)} The exact complexity is $0$ at every $\sigma$, yet the fit reports up to $1.54$, peaking near $\sigma = 0.15$ and falling away above it: at moderate noise the diffeomorphism spends capacity fitting noise, while at large noise it can no longer follow individual trajectories and the distance (which tracks the injected variance, grey) takes over instead.
    }
    \label{fig:decomp_sweeps}
\end{figure}

\paragraph{What the sweeps show.}
Complexity tracks the deformation that is genuinely required when that deformation is large, and is contaminated when it is small.
Two distinct contaminants are visible and they are separable.
The first is \emph{gauge slack}: the conjugacy is not unique, so an unpenalized fit returns an arbitrary member of the family $\homeo \circ S$ with $S$ commuting with the archetype's flow.
Panel (B) at $s=0$ isolates it in the cleanest possible way, since there the correct answer is exactly zero and the distance confirms the fit has found the identity, yet the reported complexity is $0.41$; the same floor appears on the diagonal of the rotation matrix (Supp.Sec.~\ref{sec:omega_matrix}).
The second is \emph{noise-fitting}, isolated by panel (C), where the required deformation is identically zero by construction and the reported complexity is therefore entirely attributable to capacity spent on the noise realization.
Neither contaminant is intrinsic to the measure.
Warm-starting the fit from the unpenalized solution and annealing the penalty upward recovers the analytic value in the cases where it is known, and removes the $s=0$ floor entirely, at no cost in distance.
This is the empirical counterpart of Supp.Sec.~\ref{sec:finite_time_conj}: on finite-horizon data the distance alone cannot discriminate, so the complexity must be both reported and gauge-fixed for the pair to be meaningful.

\newpage
\section{Archetype library}\label{sec:library}
This section catalogs the archetypes in our library. For each we give the defining vector field and the behavior it represents, organized by the topology of its attractor: fixed points, line and ring attractors, limit cycles, and their multistable and higher-dimensional compositions.

\subsection{Fixed point}\label{sec:fp}
\[
\begin{aligned}
\dot{x} &= \alpha x, \\
\end{aligned}
\]
with $\alpha<0$.
The learnable form is diagonal, $\dot{\vx} = \operatorname{diag}(\valpha)\,\vx$ with $\valpha<0$: its spectrum is real, so this archetype is a \emph{node} and cannot rotate.
That is a deliberate restriction rather than an omission---it is what makes the node~$\to$~spiral experiment of Sec.~\ref{sec:omega_matrix} meaningful, since the archetype cannot absorb the target's rotation into its own parameters and $\homeo$ is forced to account for it.

\subsection{Fixed-rotation linear system}\label{sec:lds_w}
For the rotation matrix of Sec.~\ref{sec:omega_matrix} we need a planar sink whose rotation is \emph{prescribed} rather than fitted,
\begin{align}\label{eq:lds_w}
    \dot{\vx} = A(\omega)\,\vx,
    \qquad
    A(\omega) = \begin{bmatrix} -a & -\omega \\ \omega & -a \end{bmatrix},
\end{align}
with $a=1$, eigenvalues $-a\pm i\omega$, and $A$ held \emph{fixed} (no learnable parameters).
Freezing $A$ is essential: were the full matrix learnable, the archetype would simply reproduce the target's own rotation, $\homeo$ would collapse to the identity, and every entry of the matrix would be zero by construction---the experiment would measure nothing.
$\omega=0$ recovers a (star) node and $\omega\neq 0$ a focus; by Sec.~\ref{sec:trivial_dynamics} all members are mutually topologically conjugate but pairwise smoothly inequivalent.

\subsection{Multi-stable}\label{sec:multistable}
The dynamics are governed by a polynomial vector field, where the time derivative $\dot{x}$ is proportional to the product
\begin{equation}
\dot{x} = \alpha \prod_{i=1}^{n} (x - r_i),
\end{equation}
 over all specified roots $r_i$, scaled by $\alpha\in \reals$.
When $\alpha < 0$, the system exhibits multistability, with stable fixed points at the specified roots.

\subsubsection{Bistable}\label{sec:bistable}
The \emph{Bistable} archetype used in the experiments is defined as
\begin{equation}
\dot x = - (x^3 - x).
\end{equation}

\begin{figure}[htbp]
    \centering
    \includegraphics[width=.5\linewidth]{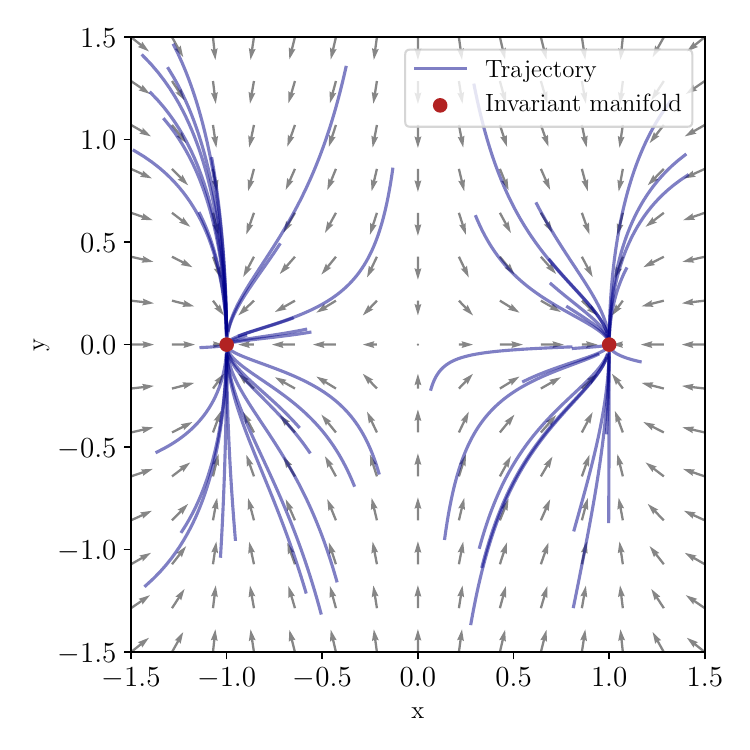}
    \caption{The \emph{Bistable} archetype's vector field and phase portrait.
     }
    \label{fig:bistable_2d}
\end{figure}

%PFTop
%\subsection{Fixed-point toplogy}\label{sec:fpt}
%We implementing any fixed-point topology dynamics in the following way.
%We consider the motif to be a directed graph. % (almost a Morse graph, but we allow for cyclicity).
%The nodes are the fixed points and the (lack of) outcoming vertices determines their stability.
%
%We embed the graph into $\reals^e$ with $e=2$ or $e=3$ based on planarity, through force-directed minimization.
%If graph is planar we embed the graph in $\reals^2$.
%If graph is not planar we embed the graph in  $\reals^3$.
%The remaining dimensions follow the dynamics $\dot x = - x$. % (or \dot x_i = - \alpha_i x_i)
%
%We define the dynamics in $\reals^e$ as follows.
%Move linearly to the closest point on graph.
%In the neighborhood of the graph we add a dynamics term that specifies speed on connecting orbits: a Gaussian bump of the connecting vector field around it.

%%LC eqs
\subsection{Limit cycle}\label{sec:lc}
The limit cycle system is defined in polar coordinates as follows:
\begin{align}%\label{eq:lc}
\dot{r} &= \alpha r(r - 1), \label{eq:lc_radial} \\
\dot{\theta} &= v, \label{eq:lc_angular} %more generally, we can allow for different speeds along the LC
\end{align}
where \( \alpha \) controls the rate of radial convergence to the unit circle, and \( v \) determines the angular velocity.
The limit cycle lies at \( r = 1 \).

Transforming to Cartesian coordinates using \( x = r \cos\theta \) and \( y = r \sin\theta \), the dynamics become:
\[
\begin{aligned}
\dot{x} &= \alpha \left( \sqrt{x^2 + y^2} - 1 \right) x - v y, \\
\dot{y} &= \alpha \left( \sqrt{x^2 + y^2} - 1 \right) y + v x.
\end{aligned}
\]
Higher dimensional: $\dot x_i = \alpha x_i$ for $i=3, \dots, D$ for $x\in \reals^D$.

For the canonical limit cycle  $\alpha=-1$ and $v=-1$.

\begin{figure}[htbp]
    \centering
    \includegraphics[width=.5\linewidth]{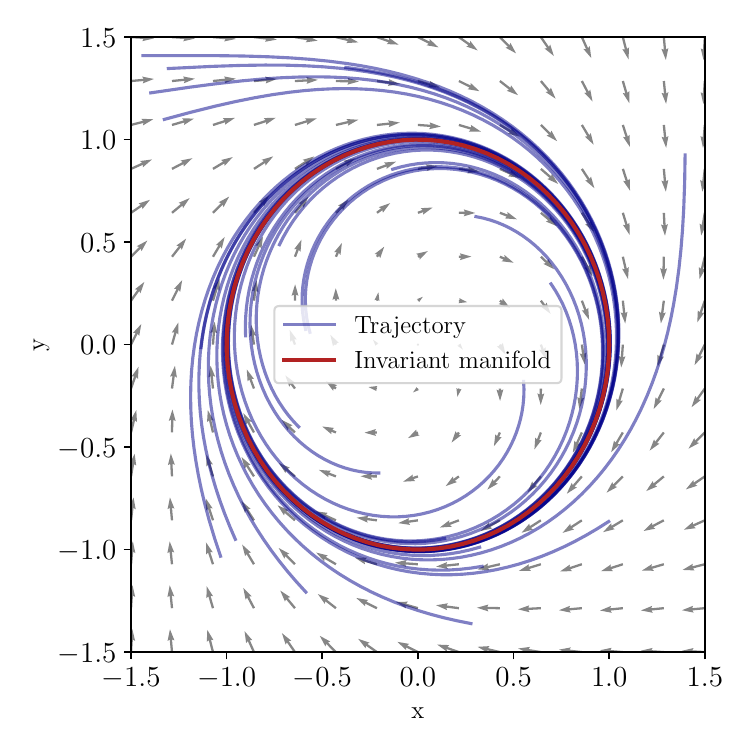}
    \caption{The \emph{limit cycle} archetype's vector field and phase portrait.
     }
    \label{fig:lc_2d}
\end{figure}

\subsubsection{Analytical limit cycle}\label{sec:analytical_lc}
Given the radial dynamics in \eqref{eq:lc_radial} the solution is:
\[
r(t) = \frac{1}{1 + C' e^{-\alpha t}}, \quad \text{where } C' = \frac{1 - r_0}{r_0}.
\]
%\[
%r(t) = \frac{1}{1 + C' e^{-t}}, \quad \text{where } C' = \frac{1 - r_0}{r_0}.
%\]

For the angular dynamics  in \eqref{eq:lc_angular} the solution is:
\[
\dot{\theta} = v
\quad \Rightarrow \quad 
\theta(t) = vt + \theta_0
\]

For the linear differential equation
\[
\dot{x} = \alpha x,
\]
the solution is:
\[
x(t) = x_0 e^{\alpha t}.
\]

%%RA eqs
\subsection{Ring attractor}\label{sec:ra}
The ring attractor system is defined in polar coordinates as follows:
\begin{equation}\label{eq:ra}
\begin{aligned}
\dot{r} &= \alpha r(r - 1), \\
\dot{\theta} &= 0.
\end{aligned}
\end{equation}

Transforming to Cartesian coordinates using \( x = r \cos\theta \) and \( y = r \sin\theta \), the dynamics become:
\[
\begin{aligned}
\dot{x} &= \alpha\left( \sqrt{x^2 + y^2} - 1 \right) x, \\
\dot{y} &= \alpha\left( \sqrt{x^2 + y^2} - 1 \right) y.
\end{aligned}
\]
Higher dimensional: $\dot x_i = \alpha x_i$ for $i=3, \dots, D$ for $x\in \reals^D$.

for canonical ring attractor $\alpha=-1$.

\begin{figure}[htbp]
    \centering
    \includegraphics[width=.5\linewidth]{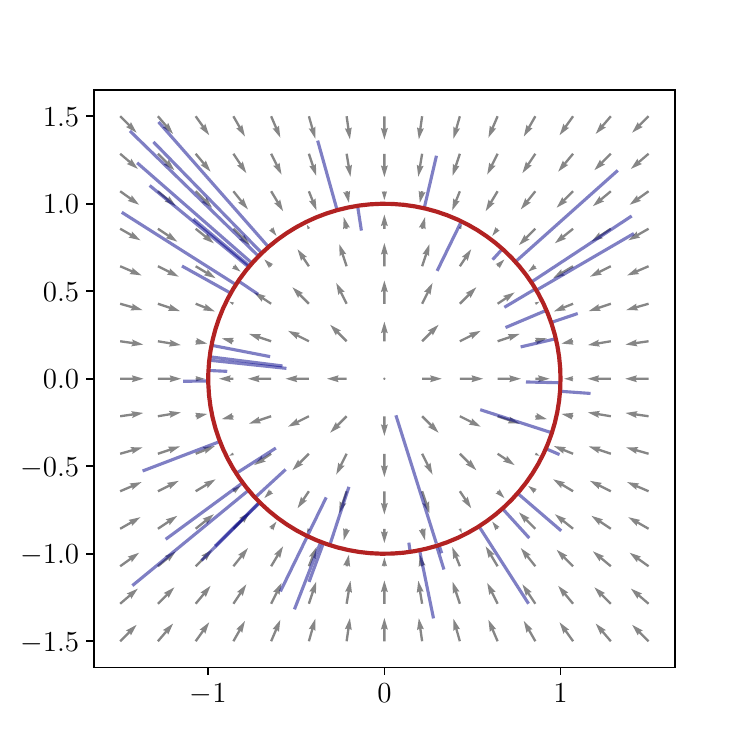}
    \caption{The \emph{ring attractor} archetype's vector field and phase portrait.
     }
    \label{fig:ra_2d}
\end{figure}

\subsubsection{Analytical}
See Supp.Sec.~\ref{sec:analytical_lc} with $v=0$.

\subsubsection{Bump-modulated perturbation on the ring}
modifying the angular dynamics from a static $\dot\theta=0$ to a \emph{bump-modulated perturbation} on the ring.
The dynamics are given by
\begin{align}
\dot{\theta} &= \exp\left(-\frac{(r - 1)^2}{2\sigma^2} \right) \cdot \psi(\theta), \label{eq:bumpmod_angular}
\end{align}
where \( \psi(\theta) \) is a learnable \( 2\pi \)-periodic function defining angular velocity perturbations on the ring, and the Gaussian bump localizes the effect to a narrow band around the unit circle.

%SPHERE
\subsection{Sphere attractor}\label{sec:sa}
The system's dynamics are described by the following ordinary differential equations (ODEs) for the sphere attractor:
\begin{equation}
 \mathbf{x}(t) = \left[ \mathbf{x}_{\text{sphere}}(t), \mathbf{x}_{\text{residual}}(t) \right],
\end{equation}
 where \( \mathbf{x}_{\text{sphere}}(t) \in \mathbb{R}^{d+1} \) represents the components of the system confined to the embedded sphere \( S^d \) and \( \mathbf{x}_{\text{residual}}(t) \in \mathbb{R}^{D - (d+1)} \) represents the residual components of the system in the higher dimensions.

The dynamics for the sphere part of the system are given by:
\[
\dot{\mathbf{x}}_{\text{sphere}}(t) = \alpha \cdot \left( \| \mathbf{x}_{\text{sphere}}(t) \| - \text{R} \right) \cdot \frac{\mathbf{x}_{\text{sphere}}(t)}{\|\mathbf{x}_{\text{sphere}}(t)\|}
\]
where \( \alpha \) is the constant controlling the strength of radial attraction. 
% \( \| \mathbf{x}_{\text{sphere}}(t) \| \) is the norm (or radial distance) of \( \mathbf{x}_{\text{sphere}}(t) \),
%and  \( \frac{\mathbf{x}_{\text{sphere}}(t)}{\|\mathbf{x}_{\text{sphere}}(t)\|} \) is the unit vector in the direction of \( \mathbf{x}_{\text{sphere}}(t) \).

The residual dynamics, representing the evolution in the residual dimensions, are given by:
\[
\dot{\mathbf{x}}_{\text{residual}}(t) = \beta \cdot \mathbf{x}_{\text{residual}}(t)
\]
where \( \beta \) is a constant that governs the attraction strength in the residual dimensions. We take $\beta < 0$; as for all archetypes this residual block is trivial, and its precise value is immaterial up to topological equivalence ($\beta=-1$ being the canonical choice, see Supp.~Sec.~\ref{sec:trivial_dynamics}).

Thus, the full system of ODEs is:
\[
\begin{aligned}
\dot{\mathbf{x}}_{\text{sphere}}(t) &= \alpha \cdot \left( \| \mathbf{x}_{\text{sphere}}(t) \| - \text{R} \right) \cdot \frac{\mathbf{x}_{\text{sphere}}(t)}{\|\mathbf{x}_{\text{sphere}}(t)\|} \\
\dot{\mathbf{x}}_{\text{residual}}(t) &= \beta \cdot \mathbf{x}_{\text{residual}}(t)
\end{aligned}
\]
where the first equation governs the radial attraction of \( \mathbf{x}_{\text{sphere}}(t) \) toward the sphere's surface, and the second equation describes the attraction of \( \mathbf{x}_{\text{residual}}(t) \) toward the origin.

\subsection{Bounded continuous attractor}\label{sec:bca}
Consider a continuous dynamical system described by the following equations. Let \( \mathbf{x}(t) \in \mathbb{R}^D \) represent the state of the system at time \( t \), where \( D \) is the dimensionality of the system.
The system exhibits the following dynamics.

\subsubsection*{1. Inside the Bounded Region}
When the state \( \mathbf{x} \) is within the hypercube \( [-B, B]^D \), the flow is zero:
\[
\frac{d \mathbf{x}}{dt} = 0, \quad \text{for } \mathbf{x} \in [-B, B]^D
\]

\subsubsection*{2. Outside the Bounded Region}
When the state \( \mathbf{x} \) is outside the hypercube \( [-B, B]^D \), the flow is directed towards the nearest boundary, with a scaling factor \( \alpha \) controlling the flow magnitude:
\[
\frac{d \mathbf{x}}{dt} = -\alpha (\mathbf{x} - \text{proj}(\mathbf{x})), \quad \text{for } \mathbf{x} \notin [-B, B]^D,
\]
where \( \text{proj}(\mathbf{x}) \) denotes the projection of \( \mathbf{x} \) onto the hypercube \( [-B, B]^D \). The projection \( \text{proj}(\mathbf{x}) \) can be written component-wise as:
\[
\text{proj}(\mathbf{x}_i) = \text{clip}(\mathbf{x}_i, -B, B) \quad \text{for each } i = 1, 2, \dots, D.
\]

Thus, for each dimension \( \mathbf{x}_i \), if \( \mathbf{x}_i \) exceeds the bounds, it is projected back onto the nearest boundary.

\subsubsection*{Residual Dynamics}
If the system has a subspace with dimension \( D - d_{\text{bca}} \), residual dynamics are added for the remaining dimensions \( \mathbf{x}_{d_{\text{bca}}+1}, \dots, \mathbf{x}_D \). These dimensions evolve according to the dynamics:
\[
\frac{d \mathbf{x}_{\text{residual}}}{dt} = \alpha \mathbf{x}_{\text{residual}}
\]

where \( \mathbf{x}_{\text{residual}} \) represents the residual part of the state vector. Here $\alpha < 0$; this is the same trivial residual block common to all archetypes, with rate immaterial up to topological equivalence and canonical value $\alpha=-1$ (Supp.~Sec.~\ref{sec:trivial_dynamics}).

\paragraph{Naming.}
The attracting set is a $d_{\text{bca}}$-dimensional box, so this one construction supplies several of the archetypes referred to elsewhere by name: $d_{\text{bca}}=1$ is the \emph{bounded line attractor} (\texttt{bla}, ``BLA''), and $d_{\text{bca}}=2$ is the \emph{plane attractor} (\texttt{plane}, ``plane''), which appears as an archetype row and a target column in the comparison matrices and as the $\vDelta_1$-free control in Supp.Sec.~\ref{sec:target_systems}.
We state this explicitly because both names are used throughout the figures while only the general $d_{\text{bca}}$ construction is defined here.

Combining the above, the system's ODE can be written as:
\[
\frac{d \mathbf{x}}{dt} =
\begin{cases}
0, & \text{if } |\mathbf{x}_i| \leq B, \, \forall i \in \{1, 2, \dots, D\} \\
-\alpha (\mathbf{x} - \text{proj}(\mathbf{x})), & \text{if } |\mathbf{x}_i| > B \, \text{for any } i. \\
\end{cases}
\]

Additionally, for the residual dynamics when \( D > d_{\text{bca}} \), the evolution is governed by:
\[
\frac{d \mathbf{x}_{\text{residual}}}{dt} = \alpha \mathbf{x}_{\text{residual}}, \quad \text{for the dimensions } d_{\text{bca}}+1, \dots, D
\]
where \( \alpha \) is a learnable parameter that scales the flow near the boundary, and \( \text{proj}(\mathbf{x}) \) is the projection of \( \mathbf{x} \) back to the hypercube \( [-B, B]^D \).

The invariant manifold of the system corresponds to the region inside the bounded region, which is the hypercube \( [-1, 1]^D \) when \( B = 1 \).
The trajectory is confined to this region as long as the state remains within the bounds.
The invariant manifold \( \mathcal{M} \) can be described as:
\[
\mathcal{M} = \{ \mathbf{x} \in \mathbb{R}^D \mid |\mathbf{x}_i| \leq 1, \, \forall i \in \{1, 2, \dots, D\} \}.
\]

\subsubsection{Bounded Line Attractor}
For $d_{bca}=1$ and $n=2$ we get a systems that is topologically conjugate to the Bounded Line Attractor from \citep{Sagodi2024a}.

\subsection{Composite systems}\label{sec:composite}
We can construct composite systems by combining two or more dynamical archetypes from the list above.
Each subsystem comes with an ODE
\begin{equation}
\dot x_i = f(x_i),
\end{equation}
with $x_i\in \reals^{D_i}$.
The composite systems dynamics is then defined as a concatenation of the $N$ subsystems dynamics
\begin{equation}
\dot x = [f(x_0), \dots f(x_N)],
\end{equation}
with $x\in \reals^{D}$ where $D=\sum_i^N D_i$.

%\subsubsection{Examples of composite systems}
%\subsubsection{Two Bounded Line Attractors}
%BLA+ Bistable
%The composite system we use for the experiments is the \emph{Two Bounded Line Attractors} (2 BLAs) system, which is the composite of the \emph{Bistable} and \emph{Bounded Line Attractor} archetypes.
%\begin{figure}[htbp]
%    \centering
%    \includegraphics[width=.5\linewidth]{twoblas_system.pdf}
%    \caption{The \emph{Two Bounded Line Attractors} system's vector field and used trajectories to train and evaluate.
%     }
%    \label{fig:twoblas_system}
%\end{figure}

%\paragraph{Torus}
%%RA+RA
%%LC+RA
%%LC+LC
%
%\paragraph{Cylinder}
%%BLA+RA
%%BLA+LC

%\subsection{Additive vector field to continuous attractor systems}
%\subsubsection{Periodic}
%\subsubsection{Non-periodic}

%%%%%%
\newpage
 \section{Target systems}\label{sec:target_systems}
This section defines the target systems we fit. These have known ground-truth structure, chosen so that the archetype a correct method should recover is unambiguous.
We also consider a stochastic version of the ring attractor and Van der Pol target systems with additive Gaussian noise:
\begin{equation}\label{eq:sde}
d\mathbf{x}(t) = f(\mathbf{x}(t), t)\,dt + \sigma\,d\mathbf{W}(t),
\end{equation}
where $\mathbf{W}(t)$ denotes a standard Wiener process and $\sigma$ is the diffusion coefficient, whose value is stated per target below ($\sigma = 0.25$ for the noisy Van der Pol system, $\sigma = 0.1$ for the noisy ring attractor).

\subsection{Limit cycles}
\subsubsection{Van der Pol oscillator}
The Van der Pol oscillator is a seminal nonlinear dynamical system that exhibits self-sustained oscillations through non-conservative damping\citep{vanderpol1926relaxation}.
The standard formulation is:
\begin{equation}
\begin{aligned}
\dot{x} &= y, \\
\dot{y} &= \mu (1 - x^2) y - x,
\end{aligned}
\end{equation}
where $\mu > 0$ controls the nonlinearity and damping. 
Unless otherwise specified, we set $\mu=0.3$.
%In Fig.~\ref{fig:all_archetypes_vdp}, we set $\mu = 0.3$.

%Also used in \citet{friedman2025characterizing}.

\begin{figure}[htbp]
    \centering
    \begin{subfigure}[b]{0.48\linewidth}
        \centering
        \includegraphics[width=\linewidth]{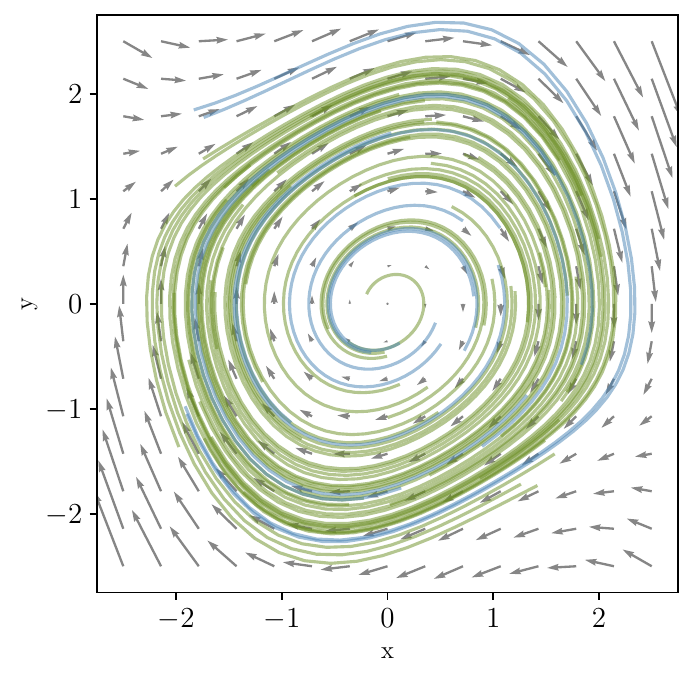}
        \caption{Deterministic}
        \label{fig:vdp_vf_trajs}
    \end{subfigure}
    \hfill
    \begin{subfigure}[b]{0.48\linewidth}
        \centering
        \includegraphics[width=\linewidth]{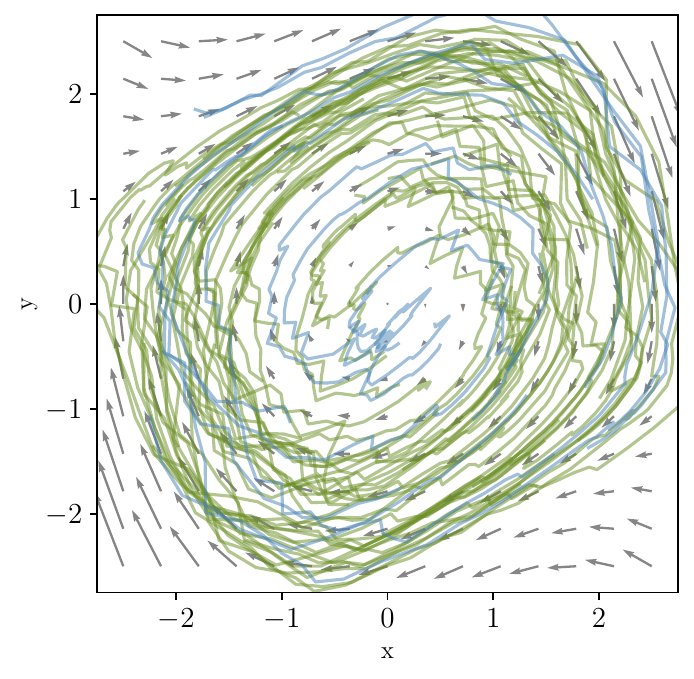}
        \caption{Noisy}
        \label{fig:vdpnoisy_vf_trajs}
    \end{subfigure}
    \caption{Deterministic and noisy Van der Pol systems.}
    \label{fig:vdp}
\end{figure}

\paragraph{Noisy Van der Pol}  
We set $\sigma = 0.25$ for the diffusion coefficient (\eqref{eq:sde}).

%MOVE TO RESULTS?
%\begin{figure}[htbp]
%    \centering
%    \includegraphics[width=.9\linewidth]{all_archetypes_vdp_noisy_var1_ntraj5_ai.pdf}
%    \caption{Comparison of fitting different archetypes to the noisy trajectories generated by a Van der Pol oscillator.
%     The archetypes considered in the fitting process are the ring attractor, limit cycle, and stable stable fixed point.
%     }
%    \label{fig:all_archetypes_vdp}
%\end{figure}

\subsubsection{Sel'kov}
The Selkov model, introduced in 1968 by E.E. Selkov, is a foundational mathematical framework for studying oscillatory dynamics in glycolysis \citep{selkov1968self}.
\begin{equation}
\begin{aligned}
\dot{x} &= - x + ay + x^2y, \\
\dot{y} &= b - ay - x^2y.
\end{aligned}
\end{equation}
We set $a= 0.05$ and $b=0.5$. 

%Used in \citep{friedman2025characterizing}.

%\begin{figure}[htbp]
%    \centering
%    \includegraphics[width=.5\linewidth]{selkov_vf_trajs.pdf}
%    \caption{
%     }
%    \label{fig:selkov_vf_trajs}
%\end{figure}

\subsubsection{Li\'enard sigmoid}
The Li\'enard sigmoid model is a nonlinear oscillator system incorporating a sigmoid function in its restoring force
\begin{equation}
\begin{aligned}
\dot{x} &= y, \\
\dot{y} &= - (1+e^{-ax})^{-1} + \tfrac{1}{2} - bx^2y.
\end{aligned}
\end{equation}
We set $a= 1.5$ and $b=-0.5$. 
% Used in \citep{friedman2025characterizing}. % + TWA

%\begin{figure}[htbp]
%    \centering
%    \includegraphics[width=.5\linewidth]{lienard_vf_trajs.pdf}
%    \caption{
%     }
%    \label{fig:lienard_vf_trajs}
%\end{figure}

\begin{figure}[htbp]
    \centering
    \begin{subfigure}[b]{0.48\linewidth}
        \centering
        \includegraphics[width=\linewidth]{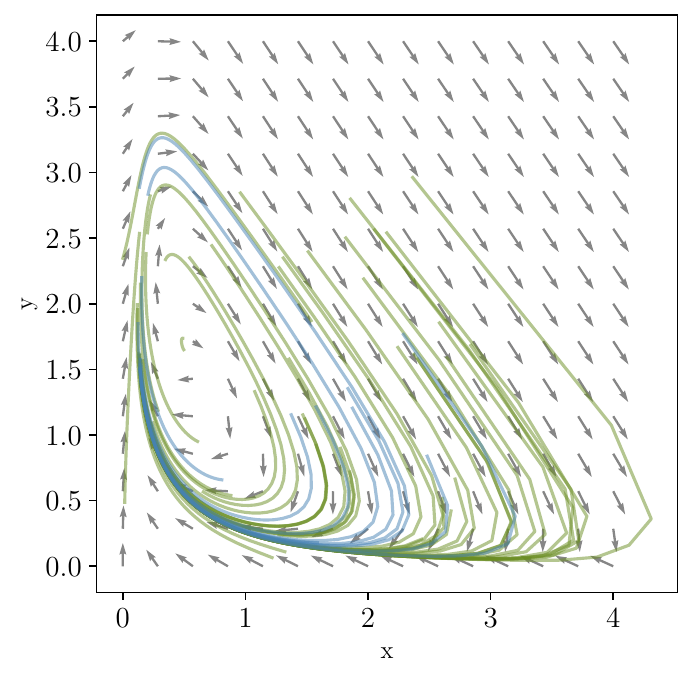}
        \caption{Selkov system}
        \label{fig:selkov_vf_trajs}
    \end{subfigure}
    \hfill
    \begin{subfigure}[b]{0.48\linewidth}
        \centering
        \includegraphics[width=\linewidth]{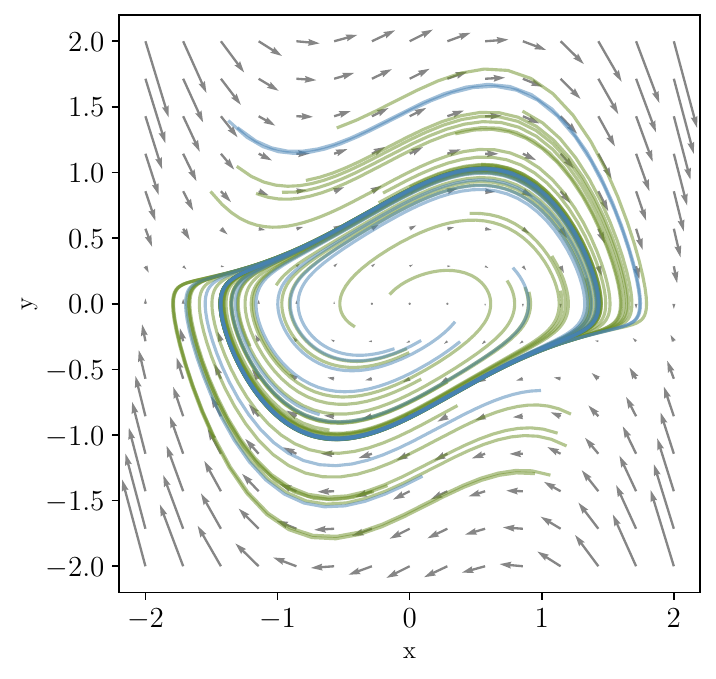}
        \caption{Li\'enard system}
        \label{fig:lienard_vf_trajs}
    \end{subfigure}
    \caption{Vector fields and trajectories of Sel'kov and Li\'enard systems.}
    \label{fig:selkov_lienard_combined}
\end{figure}

%\subsubsection{BZ reaction}
%\begin{equation}
%\begin{aligned}
%\dot{x} &= a - x - \frac{4xy}{1+x^2}, \\
%\dot{y} &= bx \frac{1-y}{1+x^2}.
%\end{aligned}
%\end{equation}

%\subsubsection{Lorenz stable limit cycle} %3D
%If $\rho>313$, then the global attractor in the Lorenz system\citep{lorenz1963deterministic} is a stable limit cycle\citep{gaiko2014global}.
%Therefore, we choose $\rho=313, \sigma=10, \beta=\tfrac{8}{3}$.

\FloatBarrier
\subsection{Approximate ring attractors}
%\subsubsection{Line attractor}

%\subsubsection{Ring attractor}
\begin{figure}[htbp]
    \centering
    \begin{subfigure}[b]{0.48\linewidth}
        \centering
        \includegraphics[width=\linewidth]{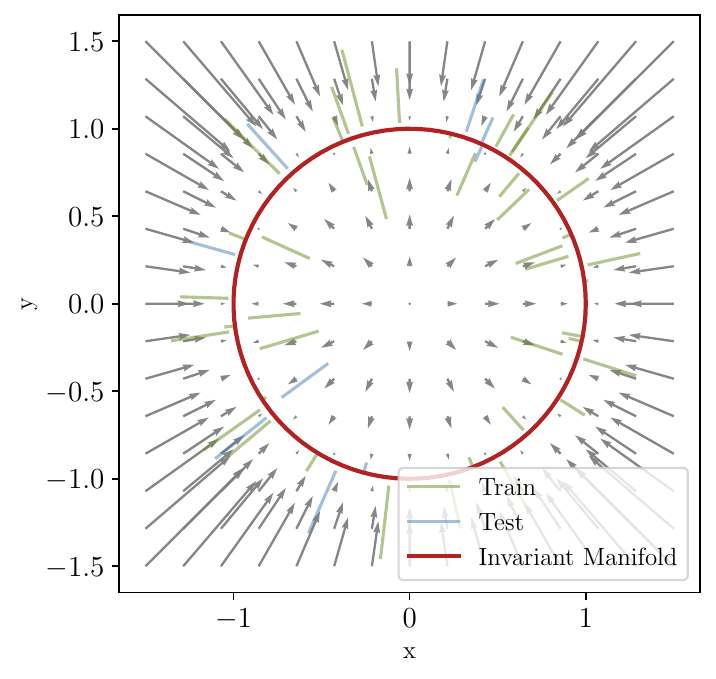}
        \caption{Deterministic}
        \label{fig:ra_vf_trajs}
    \end{subfigure}
    \hfill
    \begin{subfigure}[b]{0.48\linewidth}
        \centering
        \includegraphics[width=\linewidth]{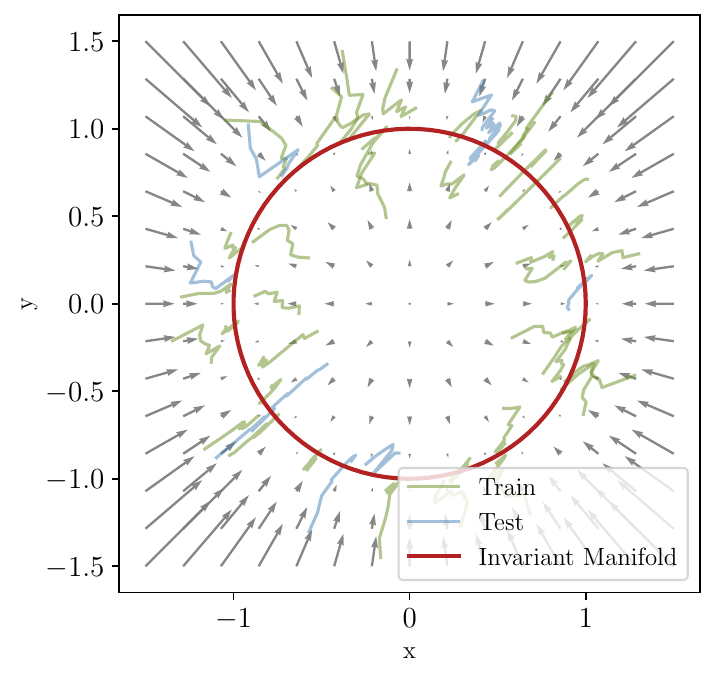}
        \caption{Noisy}
        \label{fig:ranoisy_vf_trajs}
    \end{subfigure}
    \caption{Deterministic and noisy Van der Pol systems.}
    \label{fig:ra}
\end{figure}

\paragraph{Noisy ring attractor}
We set $\sigma = 0.1$ for the diffusion coefficient (\eqref{eq:sde}).

%\subsubsection{Ring attractor}
\subsubsection{Homeomorphism based perturbation}\label{sec:homeopert_exp_description}
We sample trajectories from a ring attractor $\{x(t)\}_i^0$ by integrating the ring attractor, \eqref{eq:ra}.
The initial points are sampled randomly in the annulus 0.5 away from the ring attractor (which is the circle with radius 1).
We define a random diffeomorphism $\homeo_1$ by randomly setting the weights of a Neural ODE.
We interpolate between the identity and $\homeo_1$ by rescaling the vector field with a scalar $s \in \{0, 0.1, \dots, 1\}$ and get diffeomorphisms $\homeo_s$.
We map the trajectories with $\homeo_s$ for each $s\in \{0, 0.1, \dots, 1\}$, and get target trajectories $\{x(t)\}_i^s = \homeo_s(\{x(t)\}_i^0)$.
See for further details  Supp.Sec.~\ref{sec:homeopert_exp_details}.

\subsubsection{Vector field based perturbation}\label{sec:vfpert_exp_description}
 A smooth perturbation is sampled from a zero-mean Gaussian process with Radial Basis Function (RBF) kernel $k(x, x') = \sigma^2 \exp(-\|x - x'\|^2 / 2\ell^2)$, with a fixed lengthscale $\ell = 0.5$ (random seed $313$). The two vector-field components are sampled independently on a $40\times40$ grid over $[-3,3]^2$ and interpolated continuously.
 %     kernel = ConstantKernel(1.0, (1e-4, 1e1)) * RBF(length_scale, (1e-4, 1e1)) 
 %    #generate random vector field
%    perturb_u = np.random.multivariate_normal(np.zeros(xy.shape[0]), K).reshape(X.shape)
%    perturb_v = np.random.multivariate_normal(np.zeros(xy.shape[0]), K).reshape(X.shape)
 %LC: a rotational term $(-y, x)$ can be added to the perturbation before normalization to inject limit cycle structure.
The perturbation is normalized to have fixed norm $s \in \{0, 0.02, \dots, 0.2\}$ (11 equally spaced values). The initial conditions are sampled randomly in the annulus $0.5$ away from the unit ring.
The resulting perturbed vector field is then numerically integrated resulting in target trajectories $\{x(t)\}_i^s$.
See for further details  Supp.Sec.~\ref{sec:vfpert_exp_details}.

\paragraph{The perturbation norm is a nominal setting, not an effective one.}
The norm $s$ fixes the size of the sampled field, but not how strongly that particular field disturbs the ring: a field of a given norm may lie largely radial, where the ring absorbs it, or largely tangential, where it drives motion along the attractor. The two are not distinguishable by $s$. The effect is large enough to matter. Regenerating the single perturbed-ring target of the benchmark library from its recorded parameters, the mean angular drift over the trajectory window comes out at $13.8^\circ$, whereas the stored target drifts $1.1^\circ$; and drawing a field at $s=0.2$ produced \emph{less} drift ($3.1^\circ$) than one drawn at $s=0.07$ ($13.8^\circ$). The sweep of Figs.~\ref{fig:perturbed_both_trajectories_asy}--\ref{fig:perturbed_both_trajectories_asy_dform} is unaffected, because there a single field is drawn once and then rescaled across $s$, so the levels are comparable to each other by construction; the caveat applies to comparing a nominal $s$ across independently drawn fields. Where the strength of a perturbation matters we therefore report a measured descriptor (angular drift, endpoint radius, largest angular gap) rather than $s$ alone.

\subsection{Two Bounded Line Attractors}\label{sec:2bla}
%BLA+ Bistable
The composite system (see Supp.Sec.\ref{sec:composite}) we use for the experiments is the \emph{Two Bounded Line Attractors} (2 BLAs) system, which is the composite of the \emph{Bistable} and \emph{Bounded Line Attractor} archetypes.
\begin{figure}[htbp]
    \centering
    \includegraphics[width=.5\linewidth]{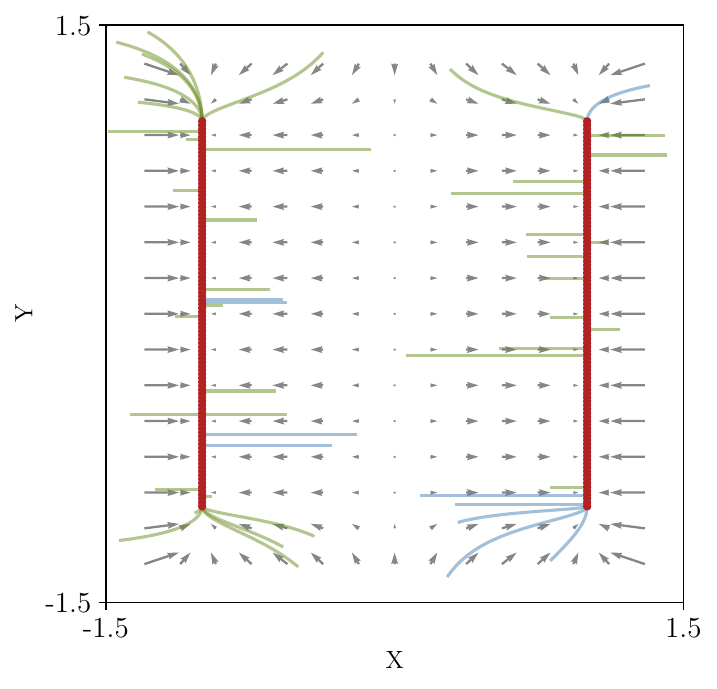}
    \caption{The \emph{Two Bounded Line Attractors} system's vector field and used trajectories to train and evaluate.
     }
    \label{fig:twoblas_vf_trajs}
\end{figure}

\subsection{Compass RNN}
%task
On each trial, the RNN receives a one-dimensional input sequence of angular velocities (positive for rightward, negative for leftward rotation), sampled from a zero-mean Gaussian Process with a radial basis function (RBF) kernel to ensure temporal smoothness. The task is to output a two-dimensional vector \( (\cos \theta(t), \sin \theta(t)) \), where \( \theta(t) \) is the time-integrated angular velocity:
\[
\theta(t) = \int_0^t \omega(s) \, ds \approx \sum_{i=0}^{t} \omega_i \cdot dt.
\]
Initial angles \( \theta(0) \) are random across trials. 
This task serves as a proxy for the internal compass computation performed by head direction cells, and has been used to study the emergence of ring attractors and persistent activity in RNN dynamics.
The trained recurrent neural networks we investigated have hidden sizes of 64, 128, or 256 units with the rectified tanh (ReLU applied to tanh) activation function for their hidden units.

\subsection{Decision-making network (Wong--Wang)}\label{sec:wongwang}
We use the two-variable reduced decision model of \citet{wong2006recurrent}, whose state $\vs = (s_1, s_2)$ holds the NMDA gating fractions of two competing populations:
\begin{align}
    \vy    &= a\,(J_N \vs + \vx_{\text{bias}}) - b, \qquad
    H(\vy)  = \frac{\vy}{1 - e^{-d \vy}}, \nonumber\\
    \frac{d\vs}{dt} &= -\frac{\vs}{\tau_s} + (1 - \vs)\,\gamma\,H(\vy),
    \label{eq:wongwang}
\end{align}
with $H$ the f--I transfer function applied componentwise.
We use the symmetric, zero-coherence setting: $J_N = \begin{psmallmatrix} 0.2609 & -0.0497 \\ -0.0497 & 0.2609\end{psmallmatrix}$~nA (self-excitation and cross-inhibition), $\vx_{\text{bias}} = (0.3411, 0.3411)$~nA, $a = 270$, $b = 108$~Hz, $d = 0.154$~s, $\gamma = 6.41\times10^{-4}$ and $\tau_s = 100$~ms.
Cross-inhibition with self-excitation makes the system \emph{bistable}, with two attracting states separated by a saddle, which is the ground-truth archetype for this target.
Trajectories are integrated by forward Euler at $dt = 1$~ms for $1500$~ms, subsampled to $T = 60$ frames, centred, rescaled by the maximum absolute value, and finally perturbed by additive observation noise of standard deviation $\sigma_{\text{obs}}$.
Initial conditions are drawn \emph{near the separatrix}, $s \sim \mathcal{U}[0.2, 0.8]$ and $\epsilon \sim \mathcal{U}[-0.06, 0.06]$ with $\vs_0 = (s + \epsilon,\, s - \epsilon)$, so that each trajectory begins close to the saddle and resolves into the well selected by $\mathrm{sign}(\epsilon)$; this makes the bistable signature explicit rather than leaving it to chance.
Three variants appear in this work and differ \emph{only} in $\sigma_{\text{obs}}$, sharing the same seed and therefore the same underlying deterministic trajectories: $\sigma_{\text{obs}} = 0.02$ (\emph{Wong--Wang}), $0.002$ (low-noise), and $0$ (the noiseless twin used for the all-deterministic variant of Fig.~\ref{fig:archetype2target}).
Note that this noise is applied to the finished trajectories and is therefore observation noise, not a diffusion term in the sense of \eqref{eq:sde}.

\subsection{Perturbed and bistable line attractors}\label{sec:pert_bla}
Both targets are built from the same two-neuron rectified-linear network,
\begin{equation}
    \frac{d\vh}{dt} = -\vh + \mathrm{relu}\!\left(W\vh + \vb\right), \qquad
    W = \begin{psmallmatrix} 0 & -1 \\ -1 & 0 \end{psmallmatrix}, \quad \vb = (C, C),
    \label{eq:relu_bla}
\end{equation}
with $C = 100$.
Every $\vh$ with $h_1 + h_2 = C$ and both rectifier arguments positive is a fixed point, so the attracting set is a \emph{segment} along the anti-diagonal centred at $(C/2, C/2)$, with contraction transverse to it and the bounds supplied by the rectifier: an exact bounded line attractor.

\emph{Perturbed BLA} (\texttt{pert\_bla}) adds independent Gaussian noise of standard deviation $0.08$ to every entry of $W$ and $\vb$ (seed $0$).
This destroys the exact degeneracy while leaving a slow manifold, and is the line-attractor analogue of the approximate ring attractors of Supp.Sec.~\ref{sec:homeopert_exp_description}.

\emph{Double-well LA} (\texttt{bla\_bistable}) instead leaves the weights exact and adds a deterministic cubic drift \emph{along} the line, so the line of fixed points collapses to two stable points with a saddle between them. We avoid calling this a ``bistable BLA'': the drift destroys the line attractor rather than coexisting with it, so the name would suggest a system that is simultaneously a line attractor and bistable.
With $\ve_L = (1, -1)/\sqrt{2}$ the unit vector along the fixed-point segment and $\alpha_L = (\vh - C/2)\cdot\ve_L$ the coordinate along it,
\begin{equation}
    \frac{d\vh}{dt} = \underbrace{-\vh + \mathrm{relu}(W\vh + \vb)}_{\text{exact BLA}} \; + \; k\,\alpha_L\left(1 - (\alpha_L/L)^2\right)\ve_L ,
    \label{eq:bla_bistable}
\end{equation}
a pitchfork normal form acting only along the manifold, with $k = 1.8$ and $L = 28$.
The continuum of fixed points collapses to two attracting points at $\alpha_L = \pm L$ separated by a saddle at the centre, so the ground-truth archetype is \emph{bistable}, not a line attractor.
This distinguishes it from the two-bounded-line-attractor target of Supp.Sec.~\ref{sec:2bla}, which remains line-attracting.

\subsection{Parameterized sweep families}\label{sec:sweep_targets}
Four one-parameter families are used for the quantitative sweeps of Supp.Sec.~\ref{sec:decomp_sweeps} and~\ref{sec:omega_matrix}; each is generated with $50$ or $60$ trajectories.

\paragraph{Rotation, \texttt{nspiral\_w}$\{\omega\}$ and \texttt{omx\_w}$\{\omega\}$.}
The planar linear family $\dot\vx = A(\omega)\vx$ with $A(\omega) = \begin{psmallmatrix} -a & -\omega \\ \omega & -a \end{psmallmatrix}$ and $a = 1$, so the eigenvalues are $-a \pm i\omega$ and, in polar form, $\dot r = -a r$ with $\dot\theta = \omega$.
The two names differ only in sampling: \texttt{nspiral} sweeps $\omega \in \{0, 0.2, \dots, 2.0\}$ over $T_{\max} = 5$ and supplies the node~$\to$~spiral column, while \texttt{omx} sweeps $\omega \in \{0, 0.2, \dots, 1.0\}$ over $T_{\max} = 2.5$ from initial conditions on an annulus ($r_0 \sim \mathcal{U}[0.8, 1.4]$) and supplies both axes of the rotation matrix.
Keeping the data off the fixed point matters because the conjugacy $\theta \mapsto \theta + (\omega/a)\ln r$ is singular as $r \to 0$.

\paragraph{On-cycle speed, \texttt{lc\_speed\_s}$\{s\}$.}
A limit cycle whose angular speed is made non-uniform,
\begin{equation}
    \dot r = r(1 - r), \qquad \dot\theta = -c(s)\left(1 + s\cos\theta\right), \qquad s \in \{0, 0.1, \dots, 0.8\}.
    \label{eq:lc_speed}
\end{equation}
The factor $c(s) = \omega_0/\sqrt{1 - s^2}$ with $\omega_0 = 1$ holds the \emph{period} fixed at $2\pi/\omega_0$ for every $s$, since $\int_0^{2\pi} \frac{d\theta}{1 + s\cos\theta} = \frac{2\pi}{\sqrt{1-s^2}}$.
This is essential: the period is a conjugacy invariant that a purely spatial diffeomorphism cannot change, so leaving it free would confound a change of on-cycle parameterization with a genuine change of dynamics.
The sign is negative to match the chirality of the limit-cycle archetype.

\paragraph{Process noise, \texttt{ra\_ns}$\{\sigma\}$.}
The ring attractor of Supp.Sec.~\ref{sec:ra} driven by additive Wiener noise in the sense of \eqref{eq:sde}, with drift $\dot r = r(1-r)$, $\dot\theta = 0$, integrated as an It\^o SDE by Euler--Maruyama at step $0.01$ and subsampled, for $\sigma \in \{0, 0.05, \dots, 0.4\}$.
Because the target \emph{is} the ring attractor and only the noise varies, no deformation is required at any $\sigma$ and the exact complexity is identically zero, which is what makes this column a control.

\newpage
\section{Architectures}\label{sec:architectures}
This section describes the two flow-based parameterizations of the diffeomorphism $\homeo$ that we use, both invertible by construction. The Neural ODE realizes $\homeo$ as the time-$T$ flow of a learned vector field, inverted by integrating backward in time; the invertible residual network stacks Lipschitz-constrained residual blocks whose inverse is recovered by fixed-point iteration.

\subsection{Neural Ordinary Differential Equations (Neural ODEs)}\label{sec:node}
Neural Ordinary Differential Equations (Neural ODEs) \citep{chen2018neural,massaroli2020nodes} are being increasingly adopted in computational and systems neuroscience, showing improved performance compared to current approaches \citep{kim2021inferring,geenjaar2023learning,sedler2023expressive,elgazzar2024universal,rubanova2019latent,coelho2024enhancing}.
Beyond their empirical success, Neural ODEs are universal approximators of dynamical systems: recent results establish $\varepsilon$-$\delta$ closeness over the infinite-time horizon $[0,\infty)$ for multistable systems, limit cycles, and normally hyperbolic continuous attractors \citep{Sagodi2026b}. %, which is precisely the class of behaviors we ask the archetype-to-target diffeomorphism to represent.
A Neural ODE is a continuous-time model where the evolution of a system is governed by an ordinary differential equation (ODE) with a neural network defining the dynamics.
Specifically, let \( x(t) \in \mathbb{R}^d \) be the state of the system at time \( t \). 
The evolution of \( x(t) \) is described by the differential equation:
\begin{equation}
    \frac{d}{dt} x(t) = f_\Theta(x(t)),
\end{equation}
where \( f_\Theta(x(t), t) \) is a Multi-Layer Perceptor (MLP) parameterized by \( \Theta \).
The solution to this equation is solved with \texttt{torchdiffeq}'s \texttt{odeint}, with initial condition \( x(0) = x_0 \), for $t\in[0,1]$\citep{torchdiffeq}.
The inverse is acquired by integrating in reverse time i.e. $t\in[1,0]$.

%MLP
We are using single hidden layer MLPs as our parameterized vector field with a ReLU activation function, see further Supp.Sec.~\ref{sec:fnn} and Supp.Sec.~\ref{sec:experiment_details}.

\subsection{Invertible ResNet}
The layers of a ResNet can be considered as Euler-discretization of the integration of a flow of a diffeomorphism\citep{rousseau2020residual}
The basic unit of an Invertible ResNet is the \textit{Invertible ResNet Block}, which operates by splitting the input \( x \in \mathbb{R}^d \) into two parts, \( x_1 \) and \( x_2 \), and applying a transformation only to \( x_1 \). The output is:
\[
 \begin{bmatrix} x_1 + f(x_1) \\ x_2 \end{bmatrix},
\]
where \( f(x_1) = W_2 \cdot \text{ELU}(W_1 \cdot x_1 + b_1) + b_2 \) is a fully connected transformation. %other activation functions?
The block is initialized with identity mapping to ensure invertibility, where \( W_2 = 0 \) and \( W_1 \) is initialized using Kaiming uniform.

The Invertible ResNet consists of stacking multiple Invertible ResNet Blocks. For input \( x \), the forward pass is:
\[
y = \text{InvertibleResNet}(x) = \text{Block}_L(\dots \text{Block}_1(x)),
\]
and the inverse is computed by reversing the block order:
\[
x = \text{InvertibleResNet}^{-1}(y) = \text{Block}_1^{-1}(\dots \text{Block}_L^{-1}(y)).
\]
implemented through 
This structure guarantees that the network is bijective, preserving the volume and making it suitable for density estimation and transformations in generative models.

%\subsection{Normalizing flow}

%\subsection{Naive MLP}
%\begin{tikzpicture}[
%    node distance=0.5cm,
%    neuron/.style={circle, draw, minimum size=0.6cm},
%    layer label/.style={font=\small, anchor=west},
%    every node/.append style={font=\small}
%  ]
%
%  % Input Layer (2 units)
%  \node[layer label] at (-.5, 0.3) {Input};
%  \foreach \i in {1,2} {
%    \node[neuron] (i\i) at (0, -\i) {};
%  }
%
%  % Hidden Layer 1 (8 units, condensed)
%  \node[layer label] at (1.5, 0.3) {$\tanh$};
%  \foreach \i/\y in {1/0.75,3/2.25} {
%    \node[neuron] (h1\i) at (2, -\y) {};
%  }
%  \node at (2, -1.5) {\vdots};
%
%  % Hidden Layer 2 (16 units, condensed)
%  \node[layer label] at (4., 0.3) {$\tanh$};
%  \foreach \i/\y in {1/0.75,3/2.25} {
%    \node[neuron] (h2\i) at (4.5, -\y) {};
%  }
%  \node at (4.5, -1.5) {\vdots};
%
%  % Output Layer (2 units)
%  \node[layer label] at (6.5, 0.3) {Output};
%  \foreach \i in {1,2} {
%    \node[neuron] (o\i) at (6.5, -\i) {};
%  }
%
%  % Connections (sparse for clarity)
%  \foreach \i in {1,2} {
%    \foreach \j in {1,3} {
%      \draw[->] (i\i) -- (h1\j);
%    }
%  }
%
%  \foreach \i in {1,3} {
%    \foreach \j in {1,3} {
%      \draw[->] (h1\i) -- (h2\j);
%    }
%  }
%
%  \foreach \i in {1,3} {
%    \foreach \j in {1,2} {
%      \draw[->] (h2\i) -- (o\j);
%    }
%  }
%
%\end{tikzpicture}

\newpage
\section{Training}
This section gives the optimization protocol shared by all fits: the optimizer and its settings, the division of trajectories into training and test sets, and the training cost.

We used the Adam optimizer \citep{kingma2014adam} with its default parameters ($\beta_1 = 0.9$, $\beta_2 = 0.999$, $\epsilon = 10^{-8}$) for all training runs.

%\subsection{Training parameters}
%
%\subsection{Cross-validation}
%\subsubsection{Validation data}
%\subsubsection{Regularization}

%\subsection{Data centering and normalization}
%Given a batch of trajectories \( \mathbf{X} \in \mathbb{R}^{B \times T \times d} \), where \( B \) is the number of trajectories, \( T \) is the number of time steps, and \( d \) is the dimensionality, the data was normalized as follows :
%
%\[
%\begin{aligned}
%\boldsymbol{\mu} &= \frac{1}{BT} \sum_{b=1}^B \sum_{t=1}^T \mathbf{X}_{b,t} \in \mathbb{R}^d, \\
%\boldsymbol{\sigma} &= \sqrt{ \frac{1}{BT} \sum_{b=1}^B \sum_{t=1}^T \left( \mathbf{X}_{b,t} - \boldsymbol{\mu} \right)^2 } \in \mathbb{R}^d, \\
%\tilde{\mathbf{X}}_{b,t} &= \frac{\mathbf{X}_{b,t} - \boldsymbol{\mu}}{\boldsymbol{\sigma}} \quad \text{for all } b = 1,\dots,B,\; t = 1,\dots,T.
%\end{aligned}
%\]
%
%The model was trained on $\tilde{\mathbf{X}}_{b,t}$.

\subsection{Training and test split}
An 80\%--20\% split was used to divide the data into training and test sets, respectively.

\subsection{Training time}
Training a single network with NODE with an MLP with a hidden layer of size 64 for 200 epochs and calculate Jacobians took around 3 minutes on an Intel Core i7 CPU and occupied 5 percent of an 8GB RAM.

\newpage
\section{Analysis}
This section defines the quantities we report from a fitted diffeomorphism and the procedure used to compare archetype and target trajectories: the trajectory-closeness (distance) score, the simplicity (complexity) score, and the invariant-manifold analysis used to visualize the recovered attractor.

\subsection{Evaluating and comparing}
The archetype trajectories data  for DSA was generated with the same initial values as the target trajectories for matching target-archetype pairs (e.g. ring attractor and noisy ring attractor).
Noisy archetypes were generated with $\sigma=0.025$ in \eqref{eq:sde}.

\subsubsection{Trajectory closeness score}\label{sec:similarity_score}
The scores shown in Fig.~\ref{fig:ring_pert_fig}A and B, as well as Fig.~\ref{fig:archetype2target}, were computed as follows.
First, the \emph{trajectory distance} between the predicted and target dynamics was quantified using the mean squared error (MSE) on the validation dataset (referred to as "distance test" in Fig.~\ref{fig:ring_pert_fig}).
For each archetype system, we normalized the trajectory distances across all targets by dividing them by the maximum distance observed for that archetype.
This normalization yields distances in the range $[0, 1]$.
The \emph{trajectory closeness score} for each target was then computed as $1$ minus the normalized distance, also resulting in a value between $0$ and $1$.
We reserve \emph{dissimilarity} and \emph{similarity} for the combined, multi-valued comparison (the trajectory distance together with the deformation complexity, Supp.Sec.~\ref{sec:simplicity_score}), of which the closeness score is only the trajectory-fit component.

The DSA closeness score was calculated through the same method by first normalizing target-wise and then subtracting the normalized score from $1$.

\subsubsection{Simplicity score}\label{sec:simplicity_score}
We calculate the complexity of the diffeomorphisms under consideration by calculating the deformation as captured by the norm of the difference between the Jacobian of the diffeomorphism and the identity as described in Supp.Sec.~\ref{sec:complexity} with the Frobenius norm (see \eqref{eq:jac_id}).
The complexity is the deformation counterpart of the trajectory distance: exactly as the trajectory \emph{closeness} score is obtained from the trajectory distance (Supp.Sec.~\ref{sec:similarity_score}), the \emph{simplicity} score is obtained from the complexity by the same per-archetype normalisation and inversion. The full \emph{(dis)similarity} of two systems is therefore multi-valued, the pair (trajectory distance, complexity), equivalently (closeness, simplicity), so that two systems count as \emph{similar} only when they are jointly close and simple; this two-dimensionality is exactly what the Pareto comparison in the main text captures.
For each archetype system, we normalized the complexity scores across all targets by dividing them by the maximum complexity observed for that archetype. 
This normalization yields complexity scores in the range $[0, 1]$. 
The \emph{simplicity score} for each target was then computed as $1$ minus the normalized complexity score, also resulting in a value between $0$ and $1$.

\subsubsection{Invariant manifold}\label{sec:inv_man}
For the existing archetypes the attractive invariant manifold is defined by design.
For example, for the ring attractor and the limit cycle archetypes the invariant manifold is the unit ring 
\[
\left\{ (x, y) \in \mathbb{R}^2 \;\middle|\; x^2 + y^2 = 1 \right\}.
\]
The invariant manifold of the model was then transformed by the fitted diffeomorphism to get the mapped source ring in Fig.~\ref{fig:ring_pert_fig}, Fig.~    \ref{fig:deformed_both_trajectories_asy}, Fig.\ref{fig:perturbed_both_trajectories_asy}, Fig.~\ref{fig:avi_rnn_recttanh} and the invariant manifolds in Supp.Sec~\ref{sec:comparison_ts}.

%\paragraph{Hausdorff distance}
%\[
%d_H(A, B) = \max\left\{ \sup_{a \in A} \inf_{b \in B} \|a - b\|, \sup_{b \in B} \inf_{a \in A} \|b - a\| \right\}
%\]
%

\newpage
\section{Numerical experiments}\label{sec:experiment_details}
This section details the setup and parameters of the numerical experiments conducted to evaluate our framework. We report on three main experimental scenarios: 
(i) a deformation applied to the ring attractor via a diffeomorphism, 
(ii) a vector field perturbation applied to the ring attractor, 
and (iii) classification of effective dynamical behavior across a diverse set of target systems.
 For each experiment, we specify the simulation settings and model architectures.

 5 different networks with the ring attractor archetype were trained for Fig.~\ref{fig:ring_pert_fig}A and B.
 A single example network was used to show the trajectories of for Fig.~\ref{fig:ring_pert_fig}A and B.

\subsection{Homeomorphism-based perturbation experiment}\label{sec:homeopert_exp_details}
See also Supp.Sec.~\ref{sec:homeopert_exp_description}.

\begin{table}[h]
\centering
\caption{Experimental parameters for the deformed ring attractor experiment.}

\label{tab:perthomeo_params}
\begin{tabular}{p{5cm}p{6cm}}
\toprule
\textbf{Category} & \textbf{Parameters} \\
\midrule
%\multirow{4}{*}{General} 
%    & Random seed: \texttt{313} \\
%   % & Save directory: \texttt{homeopert\_ring} \\
%    %& Dynamical system motif: \texttt{ring} \\
%\midrule
\multirow{4}{*}{ } 
    & Type: Neural ODE \\
Target Deformation:    & Hidden layers: \texttt{[64]} \\
Random diffeomorphism    & Activation: \texttt{ReLU} \\
    & MLP weight initialization  $\sim\mathcal{N}(0.02,0.5)$ \\
\midrule
\multirow{5}{*}{Simulation for target trajectories} %add AnalyticalDS?
    & Time step for integration: $dt=0.2$ \\
    & $T_{\text{max}} = 2$ \\
    %& Noise standard deviation: \texttt{0.0} \\
    %& Ring vector field :\\
    %& Initial condition mode: \texttt{random} \\
    & Number of trajectories: 50 \\
    %& Margin: \texttt{0.5} \\
    %& Random seed: \texttt{42} \\
    & Train/test split ratio: 0.8-0.2 \\
\midrule
\multirow{3}{*}{Learned diffeomorphism} 
    & Type: Neural ODE \\
    & Hidden layers: \texttt{[128]} \\
    & Activation: \texttt{ReLU} \\
\midrule
\multirow{1}{*}{Source}
    & Time step for integration: $dt=0.2$ \\
 \midrule
\multirow{3}{*}{Training} 
    & Learning rate:  0.01 \\
    & Epochs: 200 \\
    & Batch size: 32 \\
    %& Early stopping patience: \texttt{1000} \\
    %& Inverse formulation: \texttt{True} \\
    %& Annealing: dynamic: \texttt{False} \\
    % & Std schedule: \texttt{0.0} $\rightarrow$ \texttt{0.0} \\
\bottomrule
\end{tabular}
\end{table}

\subsection{Vector field-based perturbation experiment}\label{sec:vfpert_exp_details}
See also Supp.Sec.~\ref{sec:vfpert_exp_description}.

\begin{table}[h]
\centering
\caption{Experimental parameters for the vector field perturbed ring attractor experiment.}
\label{tab:vfpert_params}
\begin{tabular}{p{5cm}p{6cm}}
\toprule
\textbf{Category} & \textbf{Parameters} \\
\multirow{4}{*}{Perturbation}
    & Kernel: RBF, lengthscale $\ell = 0.5$ \\
    & Random seed: 313 \\
    & Grid: $40\times40$ over $[-3,3]^2$ \\
    & $\|p\|$ range $[0,0.2]$ (11 values, steps of 0.02) \\
\midrule
\multirow{4}{*}{Simulation for target trajectories}
    & Time step for integration: $dt=0.05$ \\
    &  $T_{\text{max}} = 5$ \\
    & Number of trajectories: 50 \\
    & Train/test split ratio: 0.8-0.2 \\
\midrule
\multirow{4}{*}{Learned diffeomorphism}
    & Type: Neural ODE (affine + NODE) \\
    & Hidden layers: \texttt{[128]} \\
    & Activation: \texttt{Tanh} \\
    & ODE solver: fixed-step RK4 (step 0.1) \\
\midrule
\multirow{1}{*}{Source}
    & Time step for integration: $dt=0.05$ \\
 \midrule
\multirow{3}{*}{Training}
    & Learning rate:  0.005 \\
    & Epochs: 200 \\
    & Gradient clipping: 1.0 \\
\bottomrule
\end{tabular}
\end{table}

\subsection{Classifying the effective behavior of target systems}\label{sec:manytargets_exp_details}
\begin{table}[h!]
\centering
\caption{Simulation parameters for each dynamical system.}
\begin{tabular}{lcccc}
\toprule
\textbf{System} & $dt$ &  $T_{\text{max}}$ & \textbf{Number of initial points} & \textbf{Range for initial points} \\
\midrule
Ring attractor      	&    0.2         &       2        &             50                    &      $r\pm0.5$ (annulus)  \\
Van der Pol          	&    0.1         &       5        &             50                    &       $(-2,2)\times (-2,2)$ \\
Li\'enard sigmoid    	&    0.1         &      15       &             50                    &        $(-1.5,1.5)\times (-1.5,1.5)$\\
Sel'kov                 	&    0.1         &       5        &             50                    &       $(0,3)\times (0,3)$ \\
2BLAs                  	&    0.1         &       5        &             50                    &       $(0,3)\times (0,3)$ \\
%Repressilator\citep{elowitz2000synthetic}
\bottomrule
\end{tabular}
\label{tab:system_params}
\end{table}

Initial points were sampled uniformly from the Range, except for the Ring attractor where it was sampled uniformly from theannular region (ring-shaped region) between two concentric circles of radii $0.5$ and $1.5$, centered at the origin, i.e., 
\[
\left\{ (x, y) \in \mathbb{R}^2 \;\middle|\; 0.5^2 \leq x^2 + y^2 \leq 1.5^2 \right\}.
\]

\begin{table}[h]
\centering
\caption{Experimental parameters for training diffeomorphism-based dynamical system fitting.}
\label{tab:training_params}
\begin{tabular}{p{5cm}p{6cm}}
\toprule
\textbf{Category} & \textbf{Parameters} \\
\midrule
\multirow{3}{*}{Learned diffeomorphism} 
    & Type: Neural ODE \\
    & Hidden layer sizes: \texttt{[64]} \\
    & Activation function: \texttt{ReLU} \\
\midrule
\multirow{2}{*}{Simulation for source trajectories} 
    &  $T_{\text{max}} = 5$ \\
    & Time step for integration: $dt = T_{\text{max}}/n$ \\
\midrule
\multirow{3}{*}{Training} 
    & Learning rate: 0.01 \\
    & Number of epochs: 200 \\
    & Batch size: 32 \\
\bottomrule
\end{tabular}
\end{table}

Here $n$ is the number of time points per trajectory (see Tab.~\ref{tab:system_params}, where $n= T_{\text{max}}^{\text{target}}/dt^{\text{target}}$).

\subsection{Architecture ablations}\label{sec:ablations}
The comparison our method reports depends on a learned diffeomorphism, so a natural concern is whether the recovered scores are an artifact of the particular network used to parameterize it. We therefore ablate the three choices that define this network: the activation function, the width and depth, and the diffeomorphism family (i-ResNet versus Neural ODE). Each configuration is fit to three targets that span the archetypes used in the paper, a ring attractor and two limit cycles (Van der Pol and Sel'kov), each paired with its correct archetype, over $2$ random seeds and $200$ epochs. We disable the complexity penalty ($\lambda=0$) so that the reported complexity (Frobenius Jacobian score, Supp.Sec.~\ref{sec:complexity}) reflects the deformation the fit selects on its own rather than a regularized one. Table~\ref{tab:ablations} reports the mean dissimilarity (validation MSE, Supp.Sec.~\ref{sec:similarity_score}) and complexity over targets and seeds.

\begin{table}[h!]
\centering
\caption{Architecture ablation. Mean dissimilarity (validation MSE) and complexity (Frobenius Jacobian score, Supp.Sec.~\ref{sec:complexity}) over three targets (ring attractor, Van der Pol, Sel'kov), each fit with its correct archetype and averaged over $2$ seeds, with the complexity penalty disabled ($\lambda=0$). Lower is better in both columns.}
\label{tab:ablations}
\begin{tabular}{lllcc}
\toprule
\textbf{Diffeomorphism} & \textbf{Activation} & \textbf{Width} & \textbf{Dissimilarity} & \textbf{Complexity} \\
\midrule
i-ResNet    & ReLU & $[32]$      & 0.040 & 0.40 \\
i-ResNet    & ReLU & $[64]$      & 0.037 & 0.42 \\
i-ResNet    & ReLU & $[64,64]$   & 0.021 & 0.50 \\
i-ResNet    & tanh & $[32]$      & 0.052 & 0.42 \\
i-ResNet    & tanh & $[64]$      & 0.057 & 0.40 \\
i-ResNet    & tanh & $[64,64]$   & 0.046 & 0.42 \\
i-ResNet    & ELU  & $[32]$      & 0.042 & 0.43 \\
i-ResNet    & ELU  & $[64]$      & 0.042 & 0.44 \\
i-ResNet    & ELU  & $[64,64]$   & 0.026 & 0.45 \\
\midrule
Neural ODE  & ELU  & $[64]$      & 0.033 & 0.13 \\
Neural ODE  & ReLU & $[64]$      & 0.026 & 0.76 \\
Neural ODE  & tanh & $[64]$      & 0.028 & 1.01 \\
\bottomrule
\end{tabular}
\end{table}

\paragraph{The similarity score is robust to the network.} Across all twelve configurations the dissimilarity stays within a narrow band, from $0.021$ to $0.057$, a spread smaller than the gap between individual targets (below). On the i-ResNet the activation has only a mild effect, with ReLU giving the lowest mean dissimilarity ($0.033$, averaged over widths), ELU close behind ($0.037$), and tanh highest ($0.051$). Changing the diffeomorphism family from an i-ResNet to a Neural ODE leaves the dissimilarity essentially unchanged (for example $0.037$ against $0.026$ at width $[64]$ with ReLU). The comparison the method produces therefore does not hinge on any of these choices.

\paragraph{Complexity is more sensitive to the diffeomorphism family.} While fit quality is stable, complexity is not equally invariant. At matched dissimilarity, a Neural ODE with the smooth ELU nonlinearity recovers a warp much closer to the identity (complexity $0.13$) than the i-ResNet (about $0.44$), whereas a Neural ODE with ReLU or tanh inflates the estimated deformation ($0.76$ and $1.01$), driven almost entirely by the stiff Sel'kov target (complexity $1.9$ and $2.8$ there). The i-ResNet gives an intermediate and stable complexity across activations. The practical consequence is that dissimilarity can be compared across diffeomorphism families, whereas complexity should be read relative to a fixed family.

\paragraph{Fit difficulty is set by the target, not the network.} The residual error is dominated by the target geometry rather than the architecture. Averaged over the i-ResNet configurations, the ring attractor is recovered almost exactly (dissimilarity $\approx 10^{-4}$), the Van der Pol cycle to $0.022$, and the stiff Sel'kov relaxation oscillator only to $0.104$, an ordering that holds for every configuration. Increasing depth to $[64,64]$ lowers the mean i-ResNet dissimilarity from about $0.045$ to $0.031$, concentrated on the hardest target, at the cost of a small increase in complexity, while width alone ($[32]$ against $[64]$) has no appreciable effect. All fits use the multi-step trajectory loss of Supp.Sec.~\ref{sec:loss}; a systematic comparison of one-step against multi-step prediction horizons is left to future work.

\newpage
\section{Compared methods}\label{sec:compared_methods}
This section describes the baseline methods our comparison uses and how each is configured on our systems. We begin with Dynamical Similarity Analysis and the Koopman-operator background it builds on, then the remaining baselines, stating for each the quantity it computes and the settings we used.

%%%%%%%%DSA
 \subsection{Dynamical Similarity Analysis\citep{ostrow2023beyond}}
%For DSA\citep{ostrow2023beyond}, we used $n_{\text{delays}}=3$, which is expected to suffice given that the attractors we consider are one-dimensional.
%According to Takens' theorem, a delay embedding with a sufficiently large number of delays (often 3 or more) is typically adequate for reconstructing the state space of low-dimensional systems, such as one-dimensional attractors.
%The rank of the DMD matrix used in the reduced-rank regression was set to 12 to capture the dominant dynamical modes while avoiding overfitting.

\subsubsection{Introduction to DSA}

\paragraph{Dynamic Mode Decomposition (DMD)}
 The objective of data-driven dynamical systems analysis is to approximate the infinite-dimensional Koopman operator, denoted by $\mathcal{K}$, which advances the measurements of a state $g(\mathbf{x}_k)$ forward in time such that $g(\mathbf{x}_{k+1}) = \mathcal{K}g(\mathbf{x}_k)$
While the underlying system $\dot{\mathbf{x}} = f(\mathbf{x})$ may be highly non-linear, the Koopman operator is linear, albeit acting on a specific (potentially infinite) space of observables.
 The following methods represent distinct strategies for finding a finite-dimensional approximation, $\mathbf{A}$, of this operator.

Standard Dynamic Mode Decomposition (DMD) seeks the optimal linear operator that best advances snapshots of the system state forward in time.
Given two data matrices $\mathbf{X} = [\mathbf{x}_1, \dots, \mathbf{x}_{m-1}]$ and $\mathbf{Y} = [\mathbf{x}_2, \dots, \mathbf{x}_m]$, DMD computes the operator $\mathbf{A}$ that minimizes the Frobenius norm reconstruction error $\|\mathbf{Y} - \mathbf{A}\mathbf{X}\|_F$.
The solution is analytically given by $\mathbf{A} = \mathbf{Y}\mathbf{X}^\dagger$, where $\dagger$ denotes the Moore-Penrose pseudoinverse.
In practice, this is computed via the Singular Value Decomposition (SVD) of $\mathbf{X}$ to ensure numerical stability and rank truncation.
The eigenvalues of $\mathbf{A}$ characterize the temporal behavior (growth, decay, and oscillation) of the dynamical modes, effectively performing a linear regression in the original state space.

\paragraph{Hankel Alternative View of Koopman (HAVOK)}
When the system is non-linear or only partially observed, the state variables $\mathbf{x}$ alone may not form a sufficient basis for a linear model. 
The Hankel Alternative View of Koopman (HAVOK)\citep{brunton2017chaos,villarreal2023approximate} addresses this by utilizing Takens' embedding theorem\citep{takens2006detecting}. 
Instead of analyzing raw snapshots, HAVOK constructs a Hankel matrix $\mathbf{H}$ containing time-delayed copies of the data, stacking measurement windows into columns.
An SVD is performed on this matrix, $\mathbf{H} = \mathbf{U}\Sigma\mathbf{V}^T$. The dynamics are then modeled on the time-delay coordinates (the columns of $\mathbf{V}$), which provide a linear embedding of the underlying non-linear attractor. 
The resulting linear operator $\mathbf{A}$ describes the evolution of these delay coordinates.

\subsubsection{DSA sanity checks}\label{sec:dsa_sanity}
We tested how well DSA and DMD perform on simple systems (especially focusing on the ring attractor).
We investigated the dependence of the quality of the fit and eigenspectrum on the number of trajectories used to fit the DMD model.

\begin{figure}[htbp]
    \centering
    \includegraphics[width=\linewidth]{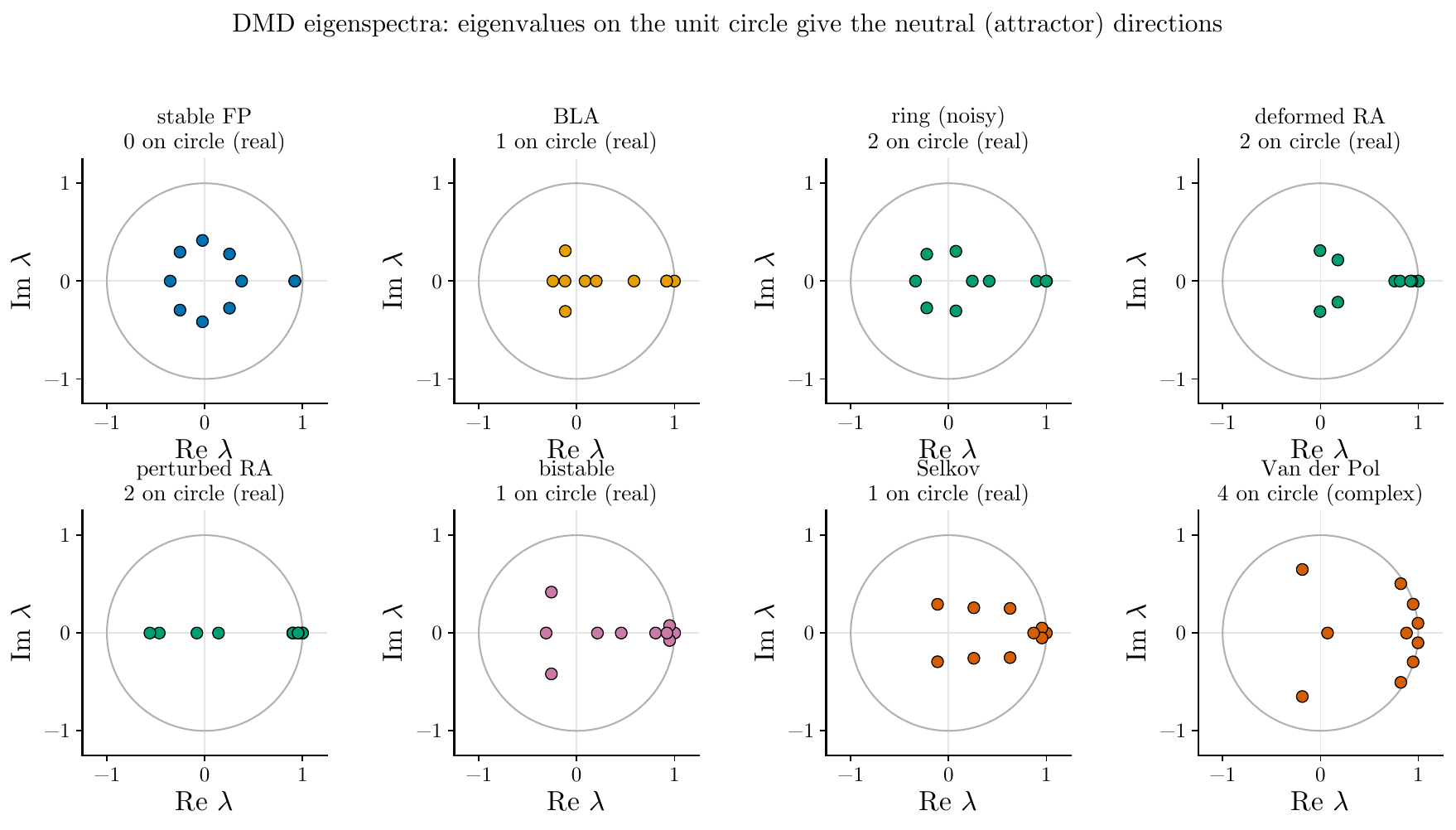}
    \caption{\textbf{DMD eigenspectra of the target library} ($5$ delays, rank $10$; eigenvalues of the reduced operator in the complex plane, grey circle is $|\lambda| = 1$).
    The count of eigenvalues on the unit circle gives the dimension of the neutrally-stable set and their argument distinguishes a continuous attractor from an oscillator: none for the stable fixed point, one for the bounded line attractor, two real ones for every member of the ring family, and a complex pair for Van der Pol.
    The topological class is therefore already encoded in the operators DSA compares (Supp.Sec.~\ref{sec:dsa_sanity}).
    }
    \label{fig:eigenspectra}
\end{figure}

\paragraph{The eigenspectrum already carries the topology.}
Because DSA compares DMD operators, what it can in principle distinguish is fixed by their spectra, so it is worth asking what those spectra contain before asking how the comparison performs.
Fig.~\ref{fig:eigenspectra} shows them for the $2$-D library ($5$ delays, rank $10$).
The classification is visible directly, and in the simplest possible form: the number of eigenvalues on the unit circle is the dimension of the neutrally-stable set, and whether they are real or complex separates a continuous attractor from an oscillator.
A stable fixed point has \emph{none} ($\max|\lambda| = 0.921$, everything decays); the bounded line attractor has \emph{one}; all three members of the ring family have \emph{two}, and both are real ($\arg\lambda = 0$), i.e.\ a two-dimensional neutral subspace with no rotation; Van der Pol instead has a complex \emph{pair} on the circle at $\arg\lambda = \pm 5.7^\circ$, the rotation frequency of its cycle, while Sel'kov's pair sits just inside the circle, spiralling onto its cycle.

This makes the accuracy DSA achieves on the same systems the more informative number.
Taking the best-scoring archetype per target, DSA classifies $7$ of $14$ correctly: it identifies the deformed and vector-field-perturbed rings, Van der Pol, Li\'enard, the bounded line attractor and both bistable targets, but assigns the \emph{noisy} ring and the noisy Van der Pol to the line attractor, Sel'kov to the plane, and the stable fixed point to the limit cycle.
The information needed to separate these cases is present in the operators DSA itself fits, as Fig.~\ref{fig:eigenspectra} shows; what the Frobenius-norm similarity between aligned operators does not do is give it decisive weight, since the eigenvalues that determine the topology lie near $|\lambda| = 1$ and contribute little to a norm dominated by the decaying bulk.
We report this as an observation about where the discriminative signal sits rather than a claim about the method's design.

\subsubsection{Robust DSA}
%\section{Robust Hyperparameter Selection for DSA}
\label{app:robust_dsa}

Dynamical Similarity Analysis (DSA) \citep{ostrow2023beyond} fits a linear dynamical system to each dataset via Hankel-DMD and measures similarity between the resulting operators.
The quality of the fit---and therefore the validity of the trajectory closeness score---depends critically on two hyperparameters: the number of time-lags $d$ used to construct the Hankel embedding, and the rank $r$ of the reduced-order DMD operator.
Standard model-selection approaches (e.g.\ minimising AIC or applying the one-standard-error rule to held-out MASE) tend to select high ranks because they optimise mean predictive accuracy without penalising cross-split instability.
Here we describe a \emph{minimal-variance} selection criterion that instead identifies the $(d, r)$ pair whose predictive behaviour is the most stable across data splits, yielding more conservative and reproducible fits.

Robust DSA replaces the pure accuracy criterion with a stability-first condition: it identifies a ``pool'' of acceptable models that achieve a mean MASE within a specified tolerance of the global minimum. 
Within this pool, it selects the (delay, rank) pair that exhibits the lowest coefficient of variation ($\text{CV} = \text{SE} / |\text{mean}|$). 

% ------------------------------------------------------------------
%\subsection{Notation}
%\label{app:robust_dsa:notation}
% ------------------------------------------------------------------

Let $\mathcal{D} = \{x^{(k)}\}_{k=1}^{K}$ denote a dataset of $K$ trajectories, each of length $T$ and living in $\mathbb{R}^N$.
For a given delay $d$ and rank $r$, let $\hat{A}_{d,r}(\mathcal{D}_\text{train}) \in \mathbb{R}^{r \times r}$ denote the DMD operator fitted on a training split $\mathcal{D}_\text{train} \subset \mathcal{D}$, and let $\text{MASE}(d, r, s)$ denote the mean absolute scaled error evaluated on the corresponding test split for data split $s \in \{1,\ldots,S\}$.

% ------------------------------------------------------------------
%\subsection{Cross-split Variance Estimation}
%\label{app:robust_dsa:cv}
% ------------------------------------------------------------------

We repeat the train/test split $S$ times (we use $S=3$ for the deformation and perturbation experiments, $S=10$ for the rank sanity checks) with independent random seeds.
For each $(d, r)$ pair we compute the mean and standard deviation of the MASE across splits:
\begin{align}
    \overline{m}(d,r) &= \frac{1}{S}\sum_{s=1}^{S} \text{MASE}(d,r,s), \\
    \widehat{\sigma}(d,r) &= \sqrt{\frac{1}{S-1}\sum_{s=1}^{S}\bigl(\text{MASE}(d,r,s) - \overline{m}(d,r)\bigr)^2}.
\end{align}
The \emph{coefficient of variation} (CV) is the scale-free ratio
\begin{equation}
    \text{CV}(d,r) = \frac{\widehat{\sigma}(d,r)}{|\overline{m}(d,r)|},
\end{equation}
and measures relative prediction instability: high CV indicates that the fitted operator changes substantially depending on the random data split, suggesting the model is over-parameterised or fitting to transient noise for that $(d,r)$ combination.

% ------------------------------------------------------------------
%\subsection{Minimal-Variance Selection Criterion}
%\label{app:robust_dsa:criterion}
% ------------------------------------------------------------------

Let $\overline{m}^* = \min_{d,r}\overline{m}(d,r)$ denote the global minimum mean MASE achieved during the grid search. 
We define a \emph{candidate pool} $\mathcal{P}$ of acceptable hyperparameters as those whose mean error falls within a tolerance factor $\eta$ of the global best:
\begin{equation}
    \mathcal{P} = \{(d,r) : \overline{m}(d,r) \leq \eta \cdot \overline{m}^*\}.
\end{equation}
By default, we set $\eta = 2.0$, ensuring that the pool only contains models that are reasonably accurate. 
The selected hyperparameters are simply those that minimise the cross-split instability within the candidate pool:
\begin{equation}
    (d^*, r^*) = \operatorname*{arg\,min}_{(d,r)\,\in\,\mathcal{P}}\; \text{CV}(d,r).
    \label{eq:robust_selection}
\end{equation}

\paragraph{Sensitivity Check.}
To ensure the selected $(d^*, r^*)$ is genuinely robust and not merely an artifact of the tolerance boundary, we perform a sensitivity check using a predefined set of strict CV thresholds, $\mathcal{T}_{\text{CV}} = \{0.05, 0.10, 0.15, 0.20, 0.25\}$. 
If the pool contains fewer than two candidates, or if the ``best'' mean model changes drastically depending on which strict CV threshold is applied, the selection is flagged as \emph{threshold sensitive}, indicating that the dataset may lack a definitive stability plateau.

% ------------------------------------------------------------------
%\subsection{Numerical Stability Upper Bound}
%\label{app:robust_dsa:ub}
% ------------------------------------------------------------------

Independently of the variance criterion, we impose a hard upper bound on the rank derived from the dimensions of the Hankel embedding.
The Hankel matrix constructed from $K$ training trajectories of length $T$ with $d$ delays has $dN$ rows and $K(T-d)$ columns, so its rank cannot exceed
\begin{equation}
    r_\text{max}(d) = \min\bigl(dN,\; K(T - d)\bigr).
    \label{eq:rank_ub}
\end{equation}
Any candidate rank $r > r_\text{max}(d)$ is excluded before the variance criterion is evaluated, preventing numerically degenerate DMD fits.
At small training fractions $K \ll K_\text{total}$, \eqref{eq:rank_ub} enforces an automatic rank reduction, which is a desirable regularisation effect.

% ------------------------------------------------------------------
%\subsection{Pairwise Selection for DSA}
%\label{app:robust_dsa:pair}
% ------------------------------------------------------------------

DSA requires a \emph{shared} $(d, r)$ for both datasets being compared, so that the two DMD operators live in a common subspace.
Given two datasets $\mathcal{D}^X$ and $\mathcal{D}^Y$, we apply the minimal-variance criterion independently to each, obtaining $(d^*_X, r^*_X)$ and $(d^*_Y, r^*_Y)$, and then take
\begin{equation}
    r^*_\text{pair} = \max(r^*_X, r^*_Y),
    \qquad
    d^*_\text{pair} = \begin{cases} d^*_X & \text{if } r^*_X \geq r^*_Y \\ d^*_Y & \text{otherwise.} \end{cases}
    \label{eq:pair_selection}
\end{equation}
Taking the maximum rank is conservative: it ensures neither dataset is forced into a subspace too small to capture its dynamics.
In practice we observe that $r^*_X \approx r^*_Y$ for the dynamical systems studied here, so \eqref{eq:pair_selection} rarely inflates the rank substantially above either individual optimum.

% ------------------------------------------------------------------
\paragraph{Algorithm Summary}
%\label{app:robust_dsa:alg}
% ------------------------------------------------------------------

\begin{algorithm}[H]
\caption{Robust hyperparameter selection for DSA}
\label{alg:robust_dsa}
\begin{algorithmic}[1]
\State \textbf{Input:} Datasets $\mathcal{D}^X, \mathcal{D}^Y$; delay grid $\mathcal{G}_d$; rank grid $\mathcal{G}_r$; $S, \eta, \mathcal{T}_{\text{CV}}$
\For{$\mathcal{D} \in \{\mathcal{D}^X, \mathcal{D}^Y\}$}
    \For{$s = 1, \ldots, S$}
        \State Split $\mathcal{D}$ into $\mathcal{D}^s_{\text{train}}, \mathcal{D}^s_{\text{test}}$ randomly
        \For{$(d, r) \in \mathcal{G}_d \times \mathcal{G}_r$ with $r \leq r_{\max}(d)$} \Comment{\eqref{eq:rank_ub}}
            \State Fit $\hat{A}_{d,r}(\mathcal{D}^s_{\text{train}})$; compute $\text{MASE}(d,r,s)$
        \EndFor
    \EndFor
    \State Compute $\overline{m}(d,r)$, $\widehat{\sigma}(d,r)$, $\text{CV}(d,r)$ for all $(d,r)$
    \State Build candidate pool $\mathcal{P} = \{(d,r) : \overline{m}(d,r) \leq \eta \cdot \min \overline{m}\}$
    \State $(d^*, r^*) \leftarrow \operatorname*{arg\,min}_{(d,r)\,\in\,\mathcal{P}}\; \text{CV}(d,r)$ \Comment{\eqref{eq:robust_selection}}
    \State Evaluate sensitivity check using thresholds $\mathcal{T}_{\text{CV}}$
\EndFor
\State $(d^*_{\text{pair}}, r^*_{\text{pair}}) \leftarrow$ \eqref{eq:pair_selection}
\State Fit DSA with $(d^*_{\text{pair}}, r^*_{\text{pair}})$ and return trajectory closeness score
\end{algorithmic}
\end{algorithm}

% ------------------------------------------------------------------
%\subsection{Hyperparameters}
%\label{app:robust_dsa:hyperparams}
%% ------------------------------------------------------------------

Table~\ref{tab:robust_dsa_hparams} summarises the hyperparameters of Robust DSA used in all experiments.

\begin{table}[H]
\centering
\caption{Hyperparameters of the Robust DSA rank selection procedure.}
\label{tab:robust_dsa_hparams}
\begin{tabular}{llll}
\toprule
Symbol & Description & Value & Rationale \\
\midrule
$S$ & \# of train/test splits & $3$--$10$ & 3 for inner loops, 10 for diagnostics \\
$\eta$ & Accuracy tolerance & $2$ & Allows up to $2\times$ the minimum MASE \\
$\mathcal{T}_{\text{CV}}$ & Sensitivity thresholds & $\{0.05, \dots, 0.25\}$ & Assesses stability plateau robustness \\
$\mathcal{G}_d$ & Delay grid & $\{1,2,3,5,10,15\}$ & Covers sub-cycle to multi-cycle lags \\
$\mathcal{G}_r$ & Rank grid & $\{2, 5,10,15,20,25,30\}$ & Below $r_\text{max}$ for all experiments \\
\bottomrule
\end{tabular}
\end{table}

 \subsection{SPE\citep{friedman2025characterizing}}
Smooth Prototype Equivalences (SPE)\citep{friedman2025characterizing} classify a system by aligning its vector field to each of a small library of canonical \emph{prototype} systems. For each prototype $\vg$, SPE learns a smooth invertible map that pushes the prototype vector field onto the target vector field, regularised to stay close to volume-preserving and centred, and reports the residual alignment error as a dissimilarity; the best-aligning prototype labels the target. The prototype library is two-dimensional and consists of a fixed point, a ring attractor, a limit cycle, a bistable system, and a line attractor. SPE acts on the vector field rather than on trajectories directly, so for trajectory data the field is first estimated by finite differences; it therefore requires an estimable vector field and, with a fixed low-dimensional prototype set, cannot represent higher-dimensional or higher-topology targets.

We modified the ring attractor prototype to better match the attraction dynamics of the target system.
In the original SPE formulation, the limit cycle is defined by:
\begin{equation}\label{eq:ra_spe}
\dot{r} = r(a - r^2), \quad \dot{\theta} = \omega,
\end{equation}
which exhibits cubic attraction towards the ring of radius \( a \).
In contrast, our ring attractor archetype (\eqref{eq:ra}) features linear radial attraction.
To facilitate a more faithful comparison and improve alignment with the target dynamics, we modify \eqref{eq:ra_spe} to adopt the same linear form.

\subsubsection{Hyperparameter tuning}\label{sec:spe_hyperparameter}
SPE aligns each prototype vector field $\vg$ to the target's (finite-difference estimated) vector field by learning a normalising-flow diffeomorphism $H$ with $n_{\text{layers}}$ coupling layers and $K$ basis functions per layer. $H$ is trained by Adam for $N_{\text{it}}$ iterations to minimise the pushed-forward vector-field alignment error, regularised by a centering term (weight $\lambda_c$, which fixes the aligned system's origin) and a determinant term (weight $\lambda_d$, which penalises the log-Jacobian so that $H$ stays close to volume-preserving); the residual alignment error is the reported dissimilarity. We use the linear-attraction ring prototype $\dot r = r(a-r),\ \dot\theta = \omega$ with $a=1,\ \omega=0$ (\eqref{eq:ra_spe} modified to the linear radial attraction of our archetype, \eqref{eq:ra}), so SPE and DAA compare against the same ring.

We fix $N_{\text{it}}=1500$, learning rate $10^{-3}$, weight decay $10^{-3}$, and $\lambda_c=0$, and choose the diffeomorphism capacity and volume regularisation once, by sweeping $K\in\{4,8\}$, $n_{\text{layers}}\in\{3,4\}$, and $\lambda_d\in\{10^{-3},10^{-2}\}$ and retaining the single configuration that maximises archetype-classification accuracy on the toy benchmark (ties broken by the lowest matched-prototype dissimilarity). This frozen configuration is then used for all SPE experiments, so no per-target tuning enters the comparison. The selected configuration is $K=4$, $n_{\text{layers}}=3$, and $\lambda_d=10^{-3}$ (benchmark archetype-classification accuracy $0.55$).

%\subsection{SPE results}\label{sec:spe_results}

 %%%%%%%%%%%DFORM
 \subsection{DFORM}
DFORM~\citep{chen2024dform, chen2026comparing} also compares two systems through their vector fields, but instead of a residual it learns an orbital diffeomorphism that aligns the two fields and reports a \emph{topological similarity}: the minimal bidirectional cosine similarity between the pushed-forward vector fields, taken as the best over $n_{\text{rep}}$ random restarts and lying in $[-1,1]$, where $1$ indicates identical flow directions. Like SPE, DFORM is vector-field-based; for trajectory data we estimate the field by finite differences and fit a small neural surrogate $\hat{\vf}$ so it can be evaluated wherever DFORM samples. Because it compares vector-field \emph{directions}, DFORM cannot separate a limit cycle from a stable spiral, which have nearly identical direction fields; this is its principal failure mode on our benchmark.

\subsubsection{Hyperparameter tuning}
We use the reference DFORM~\citep{chen2024dform,chen2026comparing} implementation (\url{https://github.com/rq-Chen/DFORM_stable}). DFORM trains an \emph{orbital} diffeomorphism that aligns the two vector fields in two stages (a linear alignment followed by a nonlinear refinement) and scores topological similarity as the minimal bidirectional cosine similarity between the pushed-forward fields, taken as the best over $n_{\text{rep}}$ random restarts. For trajectory data the target field is estimated by finite differences and represented by a small multilayer-perceptron surrogate $\hat{\vf}$ (two hidden layers, $\tanh$ activations) trained for $500$ epochs with Adam (learning rate $3\times10^{-3}$); the ring archetype is supplied analytically. The alignment samples $n_{\text{batch}}$ phase-space points for the nonlinear stage and $n_{\text{batch}}/2$ for the linear stage, over $n_{\text{rep}}$ random restarts. The configuration used for all DFORM experiments reported here is $n_{\text{batch}}=120$ (hence $60$ linear-stage points) and $n_{\text{rep}}=4$, giving a benchmark archetype-classification accuracy of $0.75$.

Two details of how DFORM is applied here are worth stating, because both affect the numbers and neither is a property of the method itself. First, DFORM scores a pair of vector fields by averaging a similarity over each model's own sampling measure, so the two models must sample comparable regions of state space for the comparison to be about dynamics rather than about scale. Sampling the archetype uniformly from a fixed box while the target is sampled from its data cloud makes the score depend on the target's spatial extent: on one system, rescaled, we measured $0.75$ at ring radius $1.4$ against $0.48$ at radius $1.0$, and on an unperturbed unit ring --- where the archetype should map to itself --- the fixed box returns a mapped ring of radius $1.29$ with a $25\%$ variation in that radius, against a truth of $1$ and $0$. We therefore sample the archetype from its own relaxed trajectories, which returns $1.03$ and $3\%$ on the same test and raises benchmark accuracy from $0.58$ to $0.75$. Second, a single restart is not stable: repeated runs on identical inputs gave mapped ring radii between $0.98$ and $1.35$, and some fits returned no model at all, so we use $n_{\text{rep}}=4$ and keep the best-scoring restart.

%\subsection{DFORM results}\label{sec:dform_results}

\newpage
\section{Additional results}
This section collects the supporting figures referenced from the main text: the full method-comparison matrices, the high-dimensional and complex-topology comparisons, and example fitted trajectories with their recovered invariant manifolds.

\subsection{Trajectory closeness comparisons}

\subsubsection{How the squares are scaled}\label{sec:squares_scaling}
Every squares figure encodes two quantities per archetype--target pair, and they are scaled differently because they answer different questions.

\paragraph{Trajectory closeness (yellow).}
The distance is shown either on a fixed absolute scale, where a full square means a distance at or below $10^{-2}$ and an empty one a distance at or above $1$, or \emph{column-normalized}, where within each target the sizes are stretched to that target's own range of archetype fits. The column-normalized form asks ``which archetype explains this target best'', and is the appropriate comparison when targets differ in ambient dimension or intrinsic difficulty; the absolute form asks ``is this fit good in absolute terms'' and is the one to read across columns. Both are on a $\log_{10}$ scale.

\paragraph{Complexity (green).}
Complexity is always shown on a \emph{fixed absolute} scale, spanning $10^{-2}$ to $2$, so a green square of a given size means the same deformation everywhere in the figure. It is deliberately not column-normalized. Normalizing complexity within a column would set the zero point from whichever archetype happens to be worst there, and that archetype is typically one that barely fits: its map is near-identity only because it never had to deform anything, so a genuine fit that had to undo a real deformation would be scored against a non-fit and rendered as almost nothing. The scale is also floored, so a reported complexity is drawn as a small square rather than vanishing; an absent green square means ``not reported'', never ``most complex''.

\paragraph{When complexity is reported.}
A complexity is only interpretable where the fit is close enough that the diffeomorphism is describing the target rather than the optimiser. We therefore report it only for pairs whose distance is within a factor of three of the best distance in that column, and leave the remaining squares yellow. The need for this is not hypothetical, and a fixed distance cutoff would not achieve it: fitting the BLA archetype to the Selkov oscillator attains a distance of $0.29$ --- unremarkable in absolute terms, but $3.6$ times the best available for that target --- while paying a complexity of $33$, roughly twenty times the largest complexity among accepted fits. Applying the criterion removes that pair and three like it, and brings the largest reported complexity from $33.4$ down to $1.67$, which is what makes a fixed scale usable at all. Such pairs are informative in their own right, since they show that a topologically mismatched archetype can be driven towards the data if the deformation is allowed to grow without bound (Supp.Sec.~\ref{sec:finite_time_conj}); they are excluded from the complexity encoding, not from the analysis.

\begin{figure}[ht]
    \centering
    \includegraphics[width=\linewidth,height=0.92\textheight,keepaspectratio]{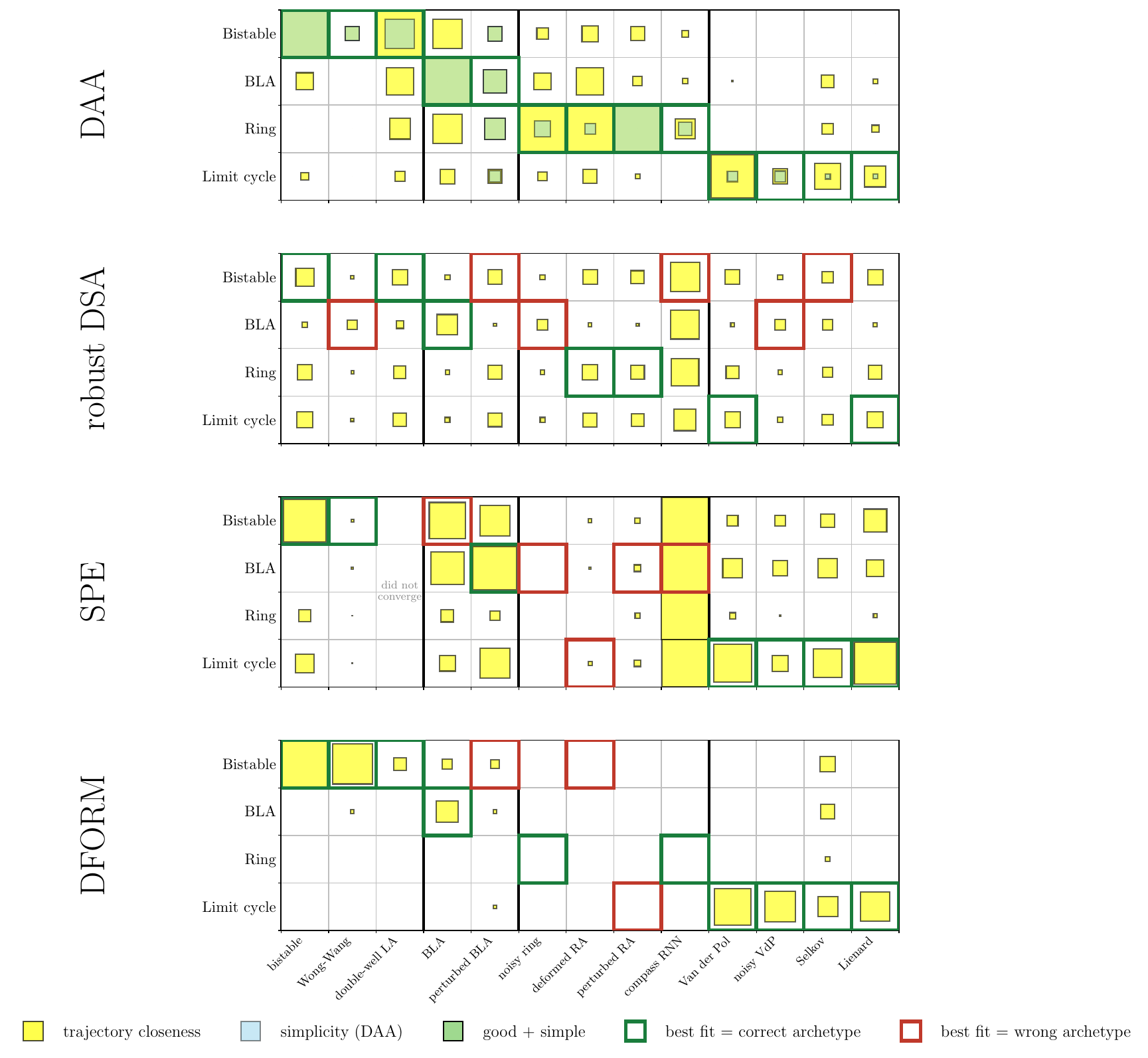}
    \caption{Archetype--target fit for the four compared methods (DAA, robust DSA, SPE, DFORM) on the toy benchmark. Rows are the five archetypes, in the same order in every panel (Stable FP, Ring, Limit cycle, BLA, Bistable); columns are the target systems. Each yellow square's area encodes \emph{trajectory closeness} on an absolute scale (a good fit fills the cell and an imperfect one is proportionally smaller, so size reflects fit quality and not merely the ranking); within each column the largest square is the best-fitting archetype.
    For DAA the blue overlay encodes \emph{simplicity} (small deformation / low complexity), so a green square marks a fit that is both good and simple. For each target (column), that best-fitting archetype is outlined in green when it is of the correct dynamical class and in red otherwise, i.e.\ whether the method selects the right archetype.
    %Archetype names are unified across methods: the single-fixed-point archetype is labelled \emph{Stable FP}, the bounded line attractor \emph{BLA}.
    }
    \label{fig:squares_compared_methods}
\end{figure}

\begin{figure}[ht]
    \centering
    \includegraphics[width=0.75\linewidth]{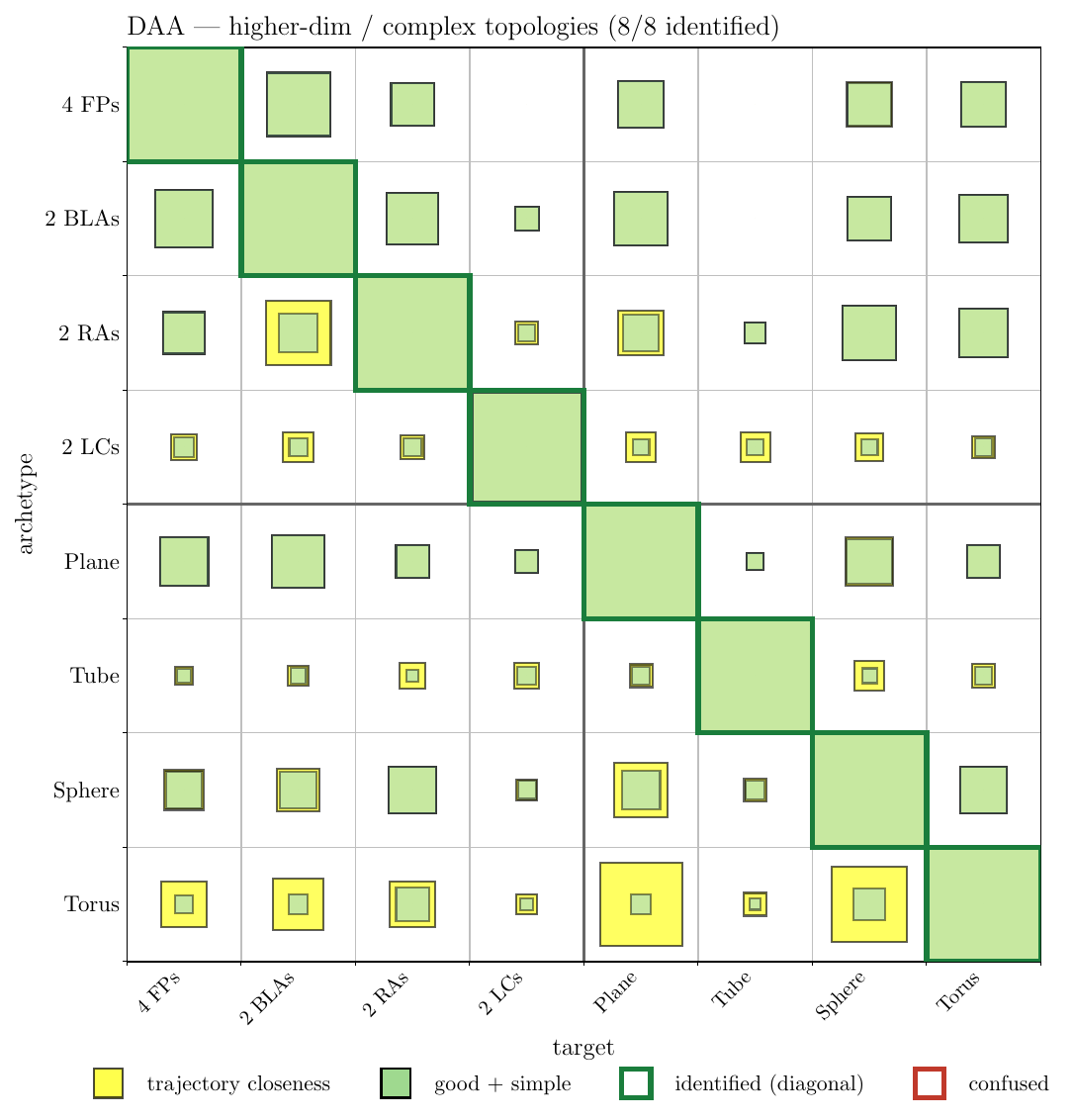}
    \caption{DAA on higher-dimensional and complex-topology structures, all generated at a common ambient dimension. Rows are archetypes and columns targets in the same order; the lines separate disconnected sets (four fixed points, two bounded line attractors, two ring attractors, two limit cycles) from connected manifolds (plane, tube, sphere, torus). Square area encodes similarity on an absolute scale (imperfect fits are proportionally smaller), and for each target the best-fitting archetype is outlined in green when it lies on the diagonal (correctly identified) and red otherwise. Because the DAA diffeomorphism is a homeomorphism it cannot bridge distinct topologies, so the diagonal dominates ($7/8$ identified); the only confusion is sphere versus torus, the two closed two-manifolds. DSA, SPE and DFORM cannot represent these structures (respectively a topology-blind delay embedding, a two-dimensional prototype library, and a direction-only vector-field comparison), so only DAA is shown.}
    \label{fig:squares_highdim}
\end{figure}

\subsection{Fitted trajectories and invariant manifolds}
In this section we compare the archetype fits to the targets. 
\subsection{Ring attractor experiments}
\subsubsection{Deformed ring attractor}
\begin{figure}[htbp]
    \centering
    \includegraphics[width=\linewidth]{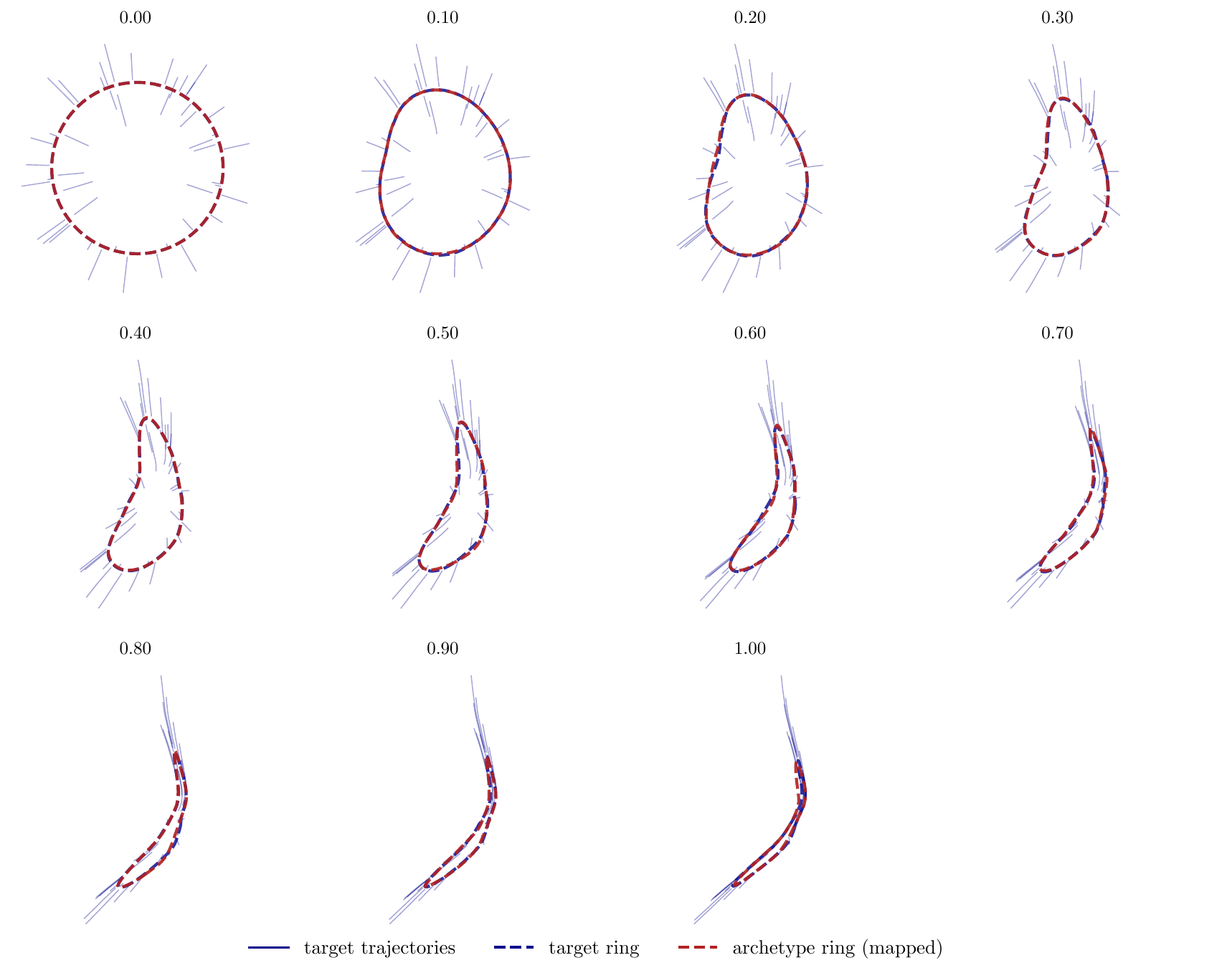}
    \caption{Our method applied to the intermediate interpolation values for the deformed ring attractor as shown in Fig.~\ref{fig:ring_pert_fig}A and described in
	Supp.Sec.~\ref{sec:homeopert_exp_description}.
	Each panel shows the target trajectories, the target ring, and the archetype's ring carried into the target's coordinates by the fitted diffeomorphism; these are the same fits reported in Fig.~\ref{fig:ring_pert_fig}.
	The invariant manifold is recovered across the whole sweep, including at $s=1$, where the ring has been folded into a strongly bent curve.
	The interpolation parameter is shown above each plot.
     }
    \label{fig:deformed_both_trajectories_asy}
\end{figure}

\begin{figure}[htbp]
    \centering
    \includegraphics[width=\linewidth]{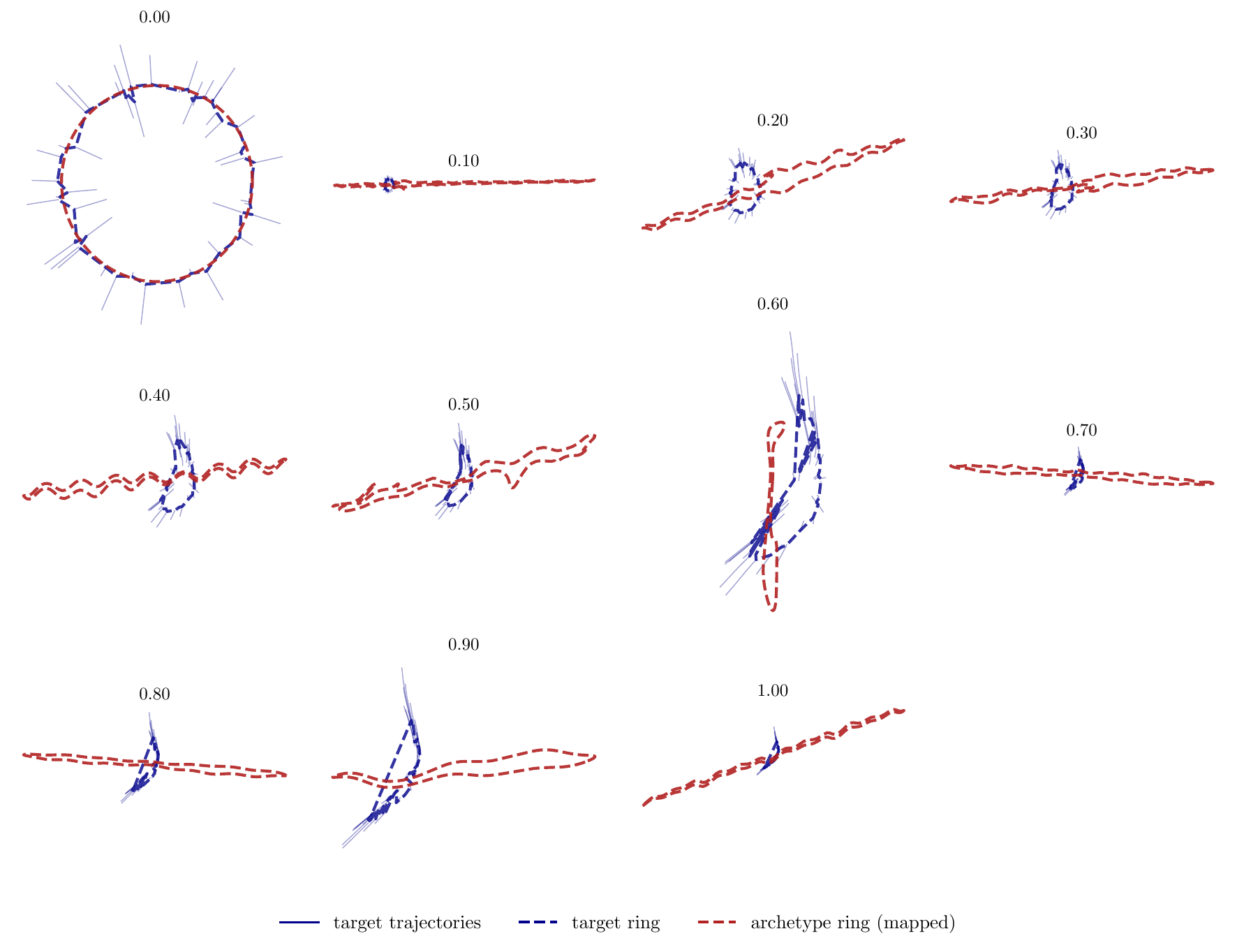}
    \caption{The Smooth Prototype Equivalences method applied to the deformed ring attractor as described in Supp.Sec.~\ref{sec:homeopert_exp_description}.
    Each panel shows the target trajectories, the target ring, and the archetype's ring carried into the target's coordinates by the fitted map, so that the panels are directly comparable to those of the other methods.
     }
    \label{fig:deformed_both_trajectories_asy_spe}
\end{figure}

\begin{figure}[htbp]
    \centering
    \includegraphics[width=\linewidth]{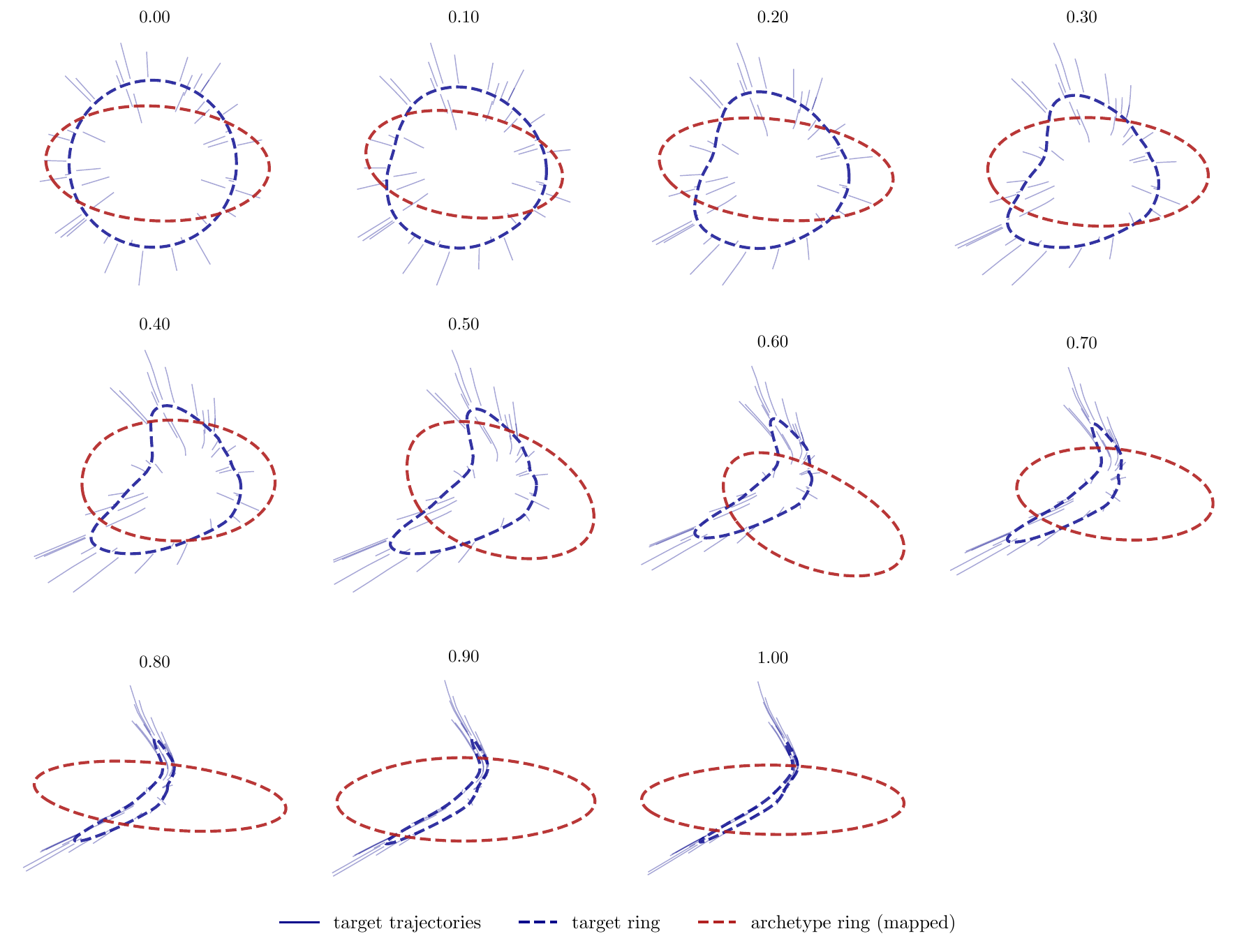}
    \caption{DFORM~\citep{chen2024dform} applied to the deformed ring attractor as described in Supp.Sec.~\ref{sec:homeopert_exp_description}, drawn with the same three curves as Fig.~\ref{fig:deformed_both_trajectories_asy_spe}.
    The archetype ring is carried into the target's coordinates by the inverse deformation, which is the direction in which DFORM's map takes the archetype to the target.
     }
    \label{fig:deformed_both_trajectories_asy_dform}
\end{figure}

\FloatBarrier
\subsubsection{Perturbed ring attractor}
\begin{figure}[htbp]
    \centering
    \includegraphics[width=\linewidth]{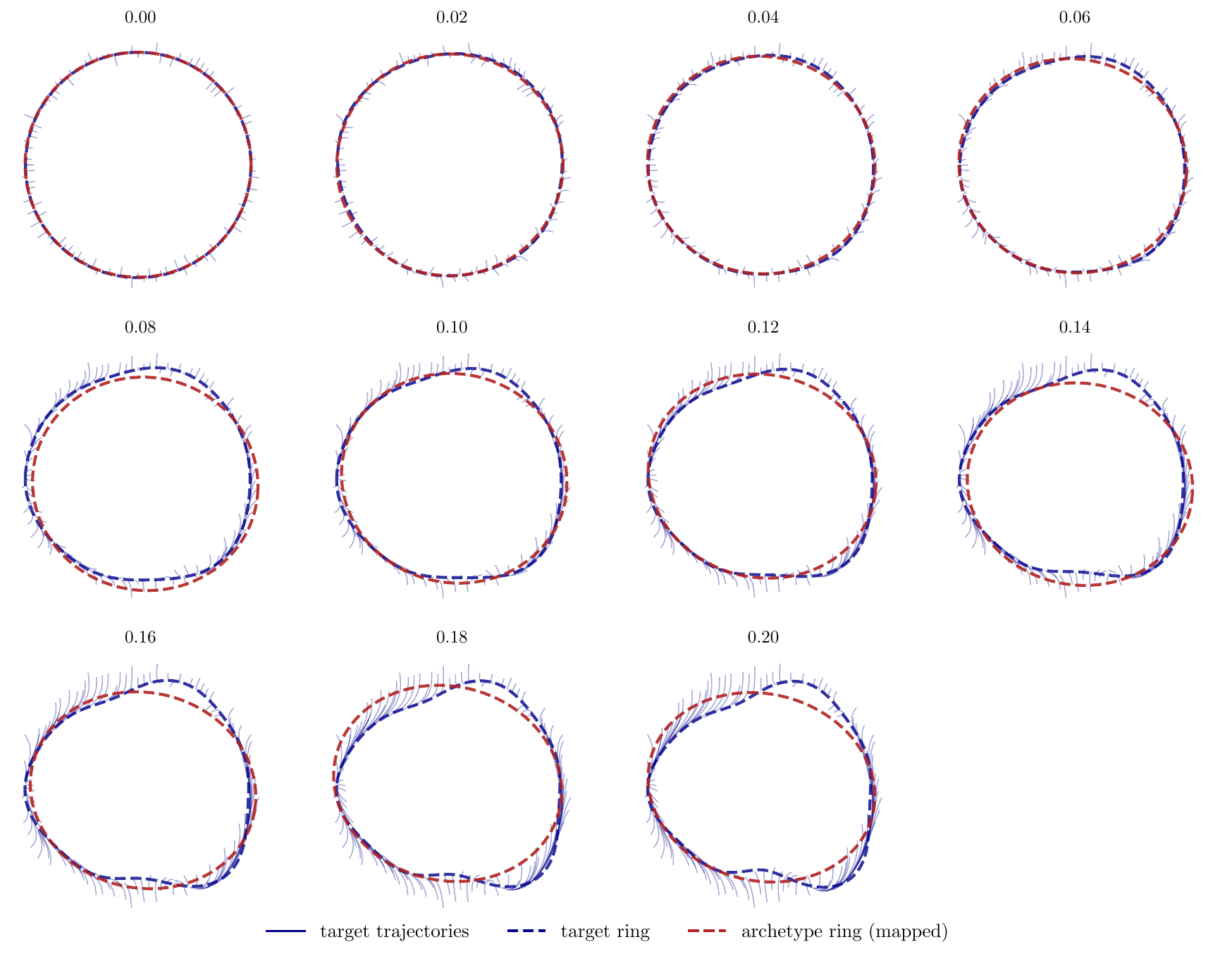}
    \caption{Our method applied to the perturbed ring attractor as shown in Fig.~\ref{fig:ring_pert_fig}B and described in
	Supp.Sec.~\ref{sec:vfpert_exp_description}, at every perturbation norm of that experiment.
	Each panel shows the target trajectories, the target ring, and the archetype's ring carried into the target's coordinates by the fitted diffeomorphism; these are the same fits reported in Fig.~\ref{fig:ring_pert_fig}.
	The invariant manifold is recovered across the sweep, with the error growing gradually and no abrupt change.
	The norm of the perturbation is shown above each subplot.
     }
    \label{fig:perturbed_both_trajectories_asy}
\end{figure}

\begin{figure}[htbp]
    \centering
    \includegraphics[width=\linewidth]{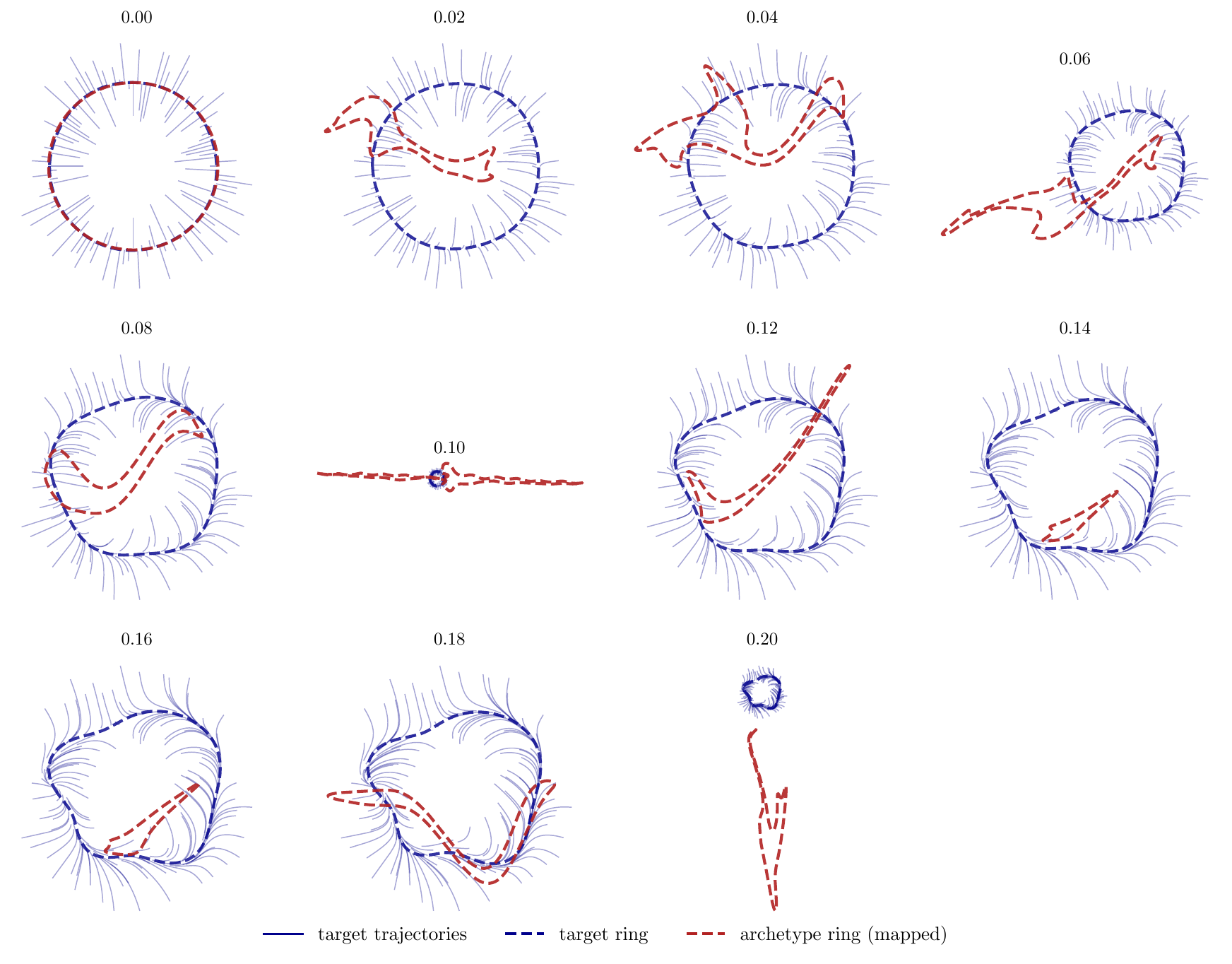}
    \caption{The Smooth Prototype Equivalences (SPE\citep{friedman2025characterizing}) method applied to the perturbed ring attractor of Supp.Sec.~\ref{sec:vfpert_exp_description}, at the same perturbation norms as Fig.~\ref{fig:perturbed_both_trajectories_asy}.
    The norm of the perturbation is shown above each subplot.
    SPE and DFORM both estimate the target vector field by finite differences, so for these two methods we re-integrate the \emph{same} perturbed systems (identical vector field, random seed and perturbation norms) from initial conditions spread across the annulus, exactly as Fig.~\ref{fig:ring_pert_fig} does for its displayed trajectories.
    This matters because the stored targets start on the ring itself, where the vector field vanishes: fitted to those samples the map is left almost unconstrained, and even the unperturbed system --- whose answer is a perfect circle --- is recovered with a $58\%$ variation in the mapped radius.
    Each level is fitted from several random initialisations and the best-scoring fit is shown, as a single initialisation is not stable across this sweep.
     }
    \label{fig:perturbed_both_trajectories_asy_spe}
\end{figure}

\begin{figure}[htbp]
    \centering
    \includegraphics[width=\linewidth]{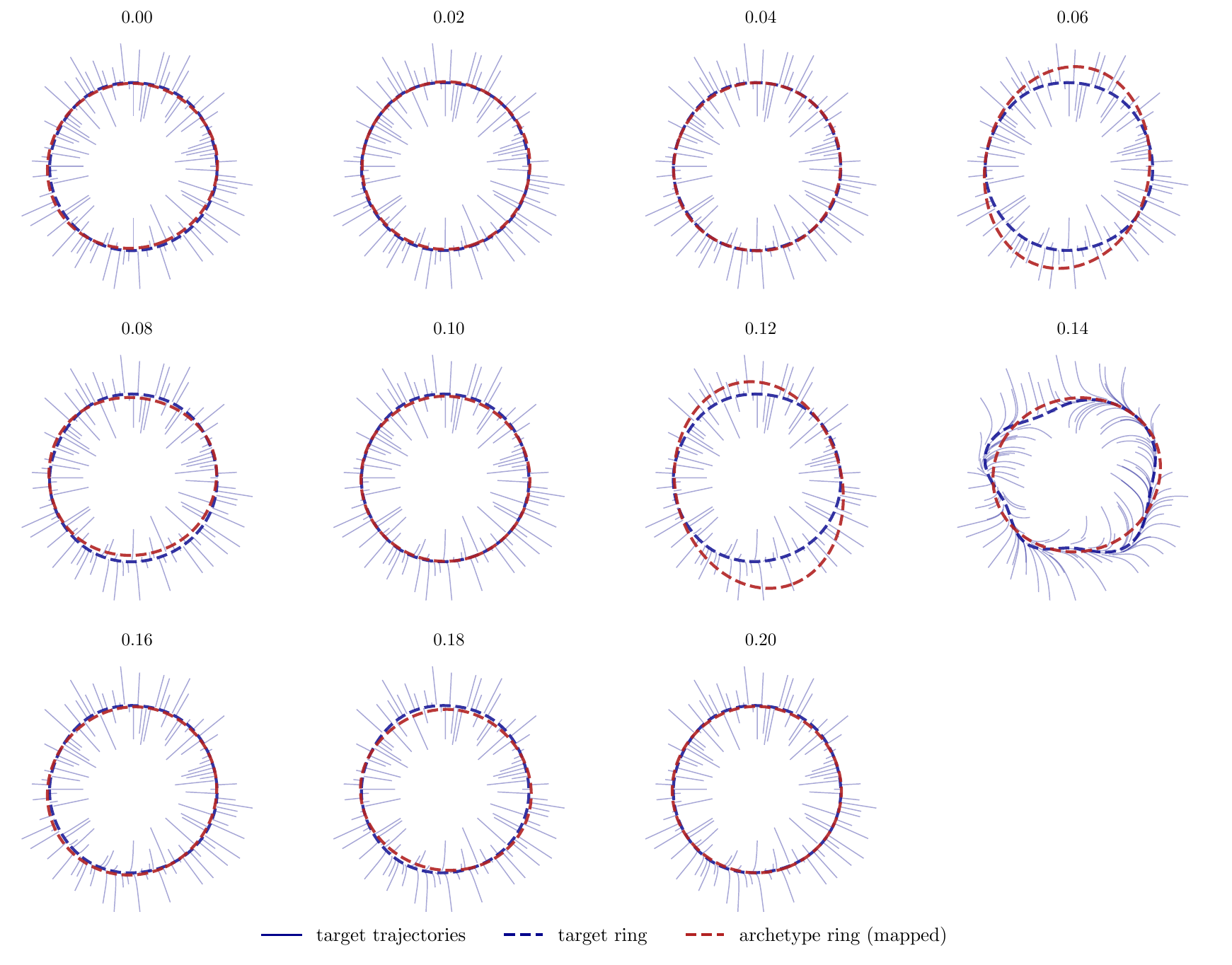}
    \caption{DFORM~\citep{chen2024dform} applied to the perturbed ring attractor of Supp.Sec.~\ref{sec:vfpert_exp_description}, drawn with the same three curves as Fig.~\ref{fig:perturbed_both_trajectories_asy_spe} and on the same re-integrated trajectories described there.
    The norm of the perturbation is shown above each subplot.
    Two departures from the DFORM configuration used for the scores reported elsewhere in this paper are needed to report the method fairly here, and both are stated for transparency.
    First, DFORM compares two vector fields by averaging a similarity over each model's own sampling measure; our archetype wrapper samples a fixed box while the target samples its data cloud, which makes the comparison depend on the target's spatial scale, so here the archetype is sampled from its own relaxed trajectories instead.
    Second, the fit has substantial run-to-run variance at a single restart --- on the unperturbed system, whose answer is a unit circle, repeated runs on identical inputs gave mapped radii between $0.98$ and $1.35$ --- so each level is fitted with a restart budget and the best-scoring fit is shown.
    With both in place the unperturbed system is recovered with a mapped radius of $0.99$ and a $1.6\%$ variation in that radius.
     }
    \label{fig:perturbed_both_trajectories_asy_dform}
\end{figure}

\paragraph{The abrupt failure in Fig.~\ref{fig:perturbed_both_trajectories_asy_spe} is a property of the fit, not of the system.}
Comparing Figs.~\ref{fig:perturbed_both_trajectories_asy} and~\ref{fig:perturbed_both_trajectories_asy_spe} on the \emph{same} targets, the SPE fit reproduces the ring at the smallest perturbations and then fails abruptly rather than degrading: at $\|p\| = 0$ the mapped ring coincides with the target to within $0.5\%$ in radius, and by $\|p\| = 0.04$ it is an order of magnitude too large.
It would be natural to read that as the ring being destroyed at a critical perturbation, and this is worth checking explicitly because it is not what happens.
Measured on the target trajectories themselves, every descriptor varies smoothly and monotonically across the sweep, with nothing distinguishing the level at which the fit fails: the mean final radius drifts from $1.000$ to $0.941$, its coefficient of variation from $0.003$ to $0.142$, the largest angular gap between endpoints from $3.0^\circ$ to $9.5^\circ$, and the mean on-manifold speed from $2.8\times10^{-3}$ to $8.8\times10^{-2}$, the last growing essentially linearly in $\|p\|$.
That is the expected picture: a generic perturbation turns the ring into a slow manifold carrying a finite number of fixed points, and the trajectories drift along it at a rate proportional to the perturbation, clumping progressively rather than discontinuously.
Our own fit on the identical systems (Fig.~\ref{fig:perturbed_both_trajectories_asy}) degrades in the same graceful way, still tracking the target at the largest perturbation with a locally growing error and no jump, and so does DFORM (Fig.~\ref{fig:perturbed_both_trajectories_asy_dform}) once it is given a restart budget.

The discontinuity therefore lies in the fitting, not the dynamics: the deformation SPE must represent grows smoothly, and beyond some magnitude the fit stops degrading and instead fails outright.
This is the same signature we encounter elsewhere in this work whenever a required deformation grows past what a fitting procedure will find, notably the non-convergence of the largest-$|\Delta\omega|$ cells under a cold-started complexity penalty (Supp.Sec.~\ref{sec:omega_matrix}).
The natural explanation would be an unlucky initialisation, and we tested it directly by refitting the collapsed levels from twelve random initialisations: it is not that.
At $\|p\| = 0.04$ every one of the twelve returns a mapped ring of radius $8.6$--$10.1$ against a target of $0.995$, so no initialisation in that sample recovers the ring.
The failure is also invisible to SPE's own fit score, which is what a practitioner would use to choose among restarts: across those twelve fits the correlation between the score and the mapped ring's radial variation is $-0.12$ at $\|p\| = 0.10$ and $+0.37$ at $\|p\| = 0.04$, so the best-scoring fit is not the most faithful one.
This is the opposite of what we find for DFORM, whose score does track fidelity closely enough that selecting on it turns an unstable procedure into a stable one, and it is why a restart budget repairs one baseline and not the other.

\FloatBarrier
\subsubsection{Combined comparison across deformation and perturbation}\label{sec:fig2_merged}
\begin{figure}[htbp]
    \centering
    \includegraphics[width=\linewidth]{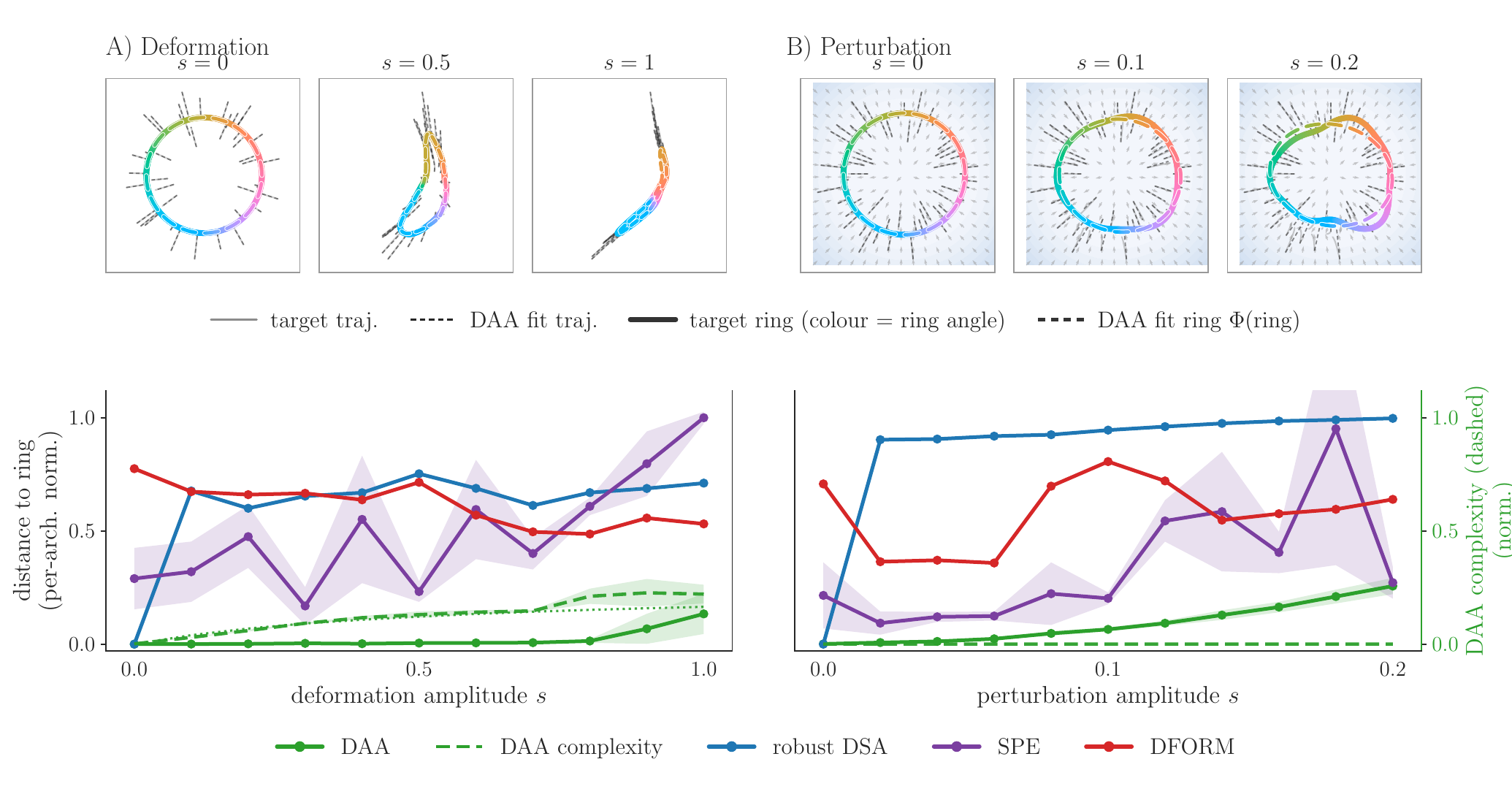}
    \caption{Merged comparison figure. (\textbf{a}) Homeomorphism deformation and (\textbf{b}) vector-field perturbation of the ring attractor, showing target (dark blue) versus DAA fit (dark red) trajectories and invariant rings; the perturbation panels also draw the perturbed vector field.
    (\textbf{c},\textbf{d}) Trajectory distance to the ring archetype for DAA, robust DSA, SPE, and DFORM (per-archetype normalized), with the DAA deformation complexity as a dashed line.
    A version showing trajectory closeness and simplicity ($1-$normalized value) is available, as is a variant drawing the equivalent (pullback) vector field on the deformation panels.}
    \label{fig:fig2_merged}
\end{figure}

Figure~\ref{fig:fig2_merged} summarizes both distortion types on a single, shared per-archetype axis.
Both halves use the same Neural ODE diffeomorphism: the deformation half reuses the well-trained fits (mean over five random seeds), and the perturbation half is refit with the same architecture on the vector-field-perturbed sweep of Supp.Sec.~\ref{sec:vfpert_exp_description} (norm $\|p\|\in\{0,0.02,\dots,0.2\}$).
The top row shows, for three amplitudes each, the target trajectories and rings (dark blue) against the DAA fit (dark red); the perturbation panels also show the perturbed vector field.
The bottom row compares the trajectory distance of DAA against robust DSA\citep{ostrow2023beyond}, SPE\citep{friedman2025characterizing}, and DFORM\citep{chen2024dform}, together with the DAA deformation complexity (dashed), all normalized per archetype.
For the diffeomorphic deformation (left) the DAA distance stays near zero while the complexity grows, whereas the compared methods immediately report a large dissimilarity; for the vector-field perturbation (right) the DAA distance grows slowly as the fixed-point continuum breaks.

\clearpage  
\subsection{Task trained RNNs: Internal compass model}
\begin{figure}[htbp]
    \centering
    \includegraphics[width=0.95\linewidth]{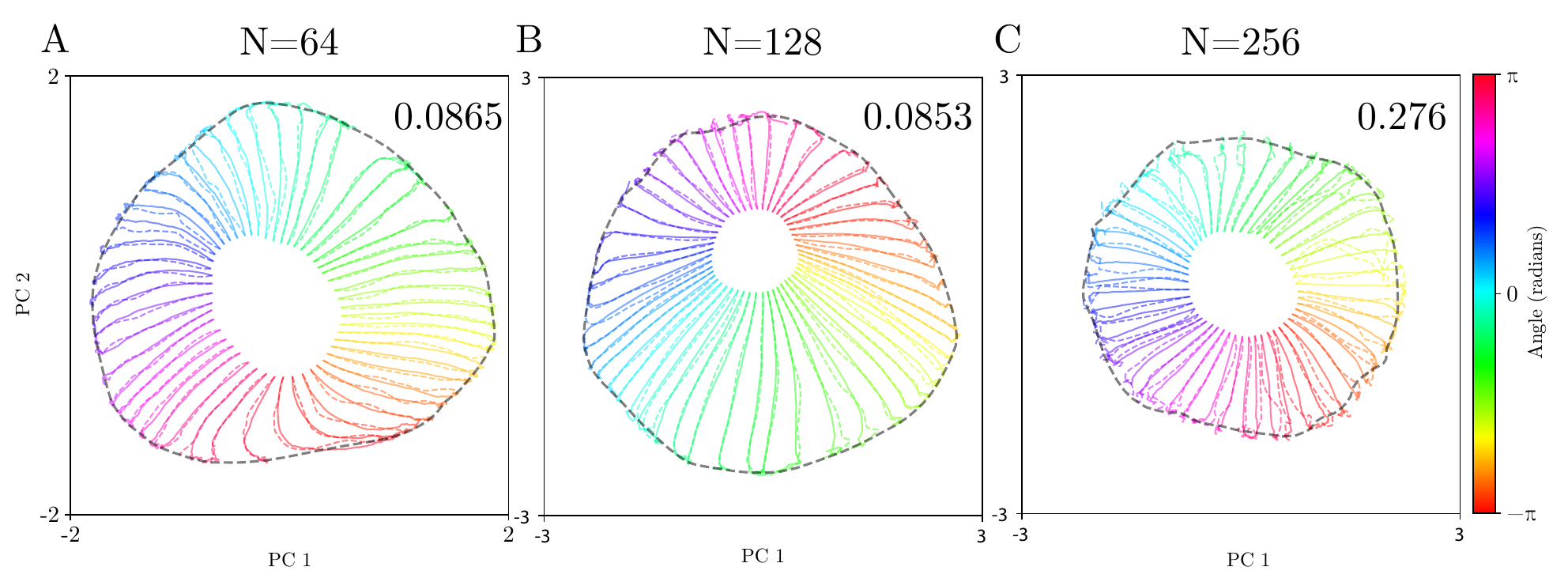}
    \caption{The ring attractor archetype fit to RNNs of different sizes ($N=64,128,256$ hidden units) trained on an angular velocity integration task from \citet{Sagodi2024a}.
    The diffeomorphism brings the canonical ring attractor near the invariant slow manifold (composed out of a finite number of fixed points) even though topological equivalency is lacking.
    Trajectories colored by memory content (angular memory).
    The MSE similarity is reported in the upper right corner for each network.
    }
    \label{fig:avi_rnn_recttanh}
\end{figure}

To study low-dimensional attractors embedded in high-dimensional neural dynamics, we fit our canonical archetype library to RNNs trained on an angular velocity integration task, designed to model the internal compass in animals \citep{Sagodi2024a}. 
The ring attractor archetype captures both the transient trajectories and the invariant manifold across network sizes (Fig.~\ref{fig:avi_rnn_recttanh}).

\clearpage  
\subsection{Comparison between the target and the mapped source trajectories}\label{sec:comparison_ts}
Each figure in this section takes one archetype and fits it to every target in turn.
A panel overlays the target trajectories with the archetype trajectories carried into the target's coordinates by the fitted map, together with the image of the archetype's invariant manifold.
The archetype trajectories are started from the target's own initial conditions, pulled back through the map, so that the two families of curves are directly comparable rather than two unrelated bundles.
Panels are framed on the extent of the target: when a fit maps the archetype's invariant manifold far outside the data the curve is simply clipped, which is itself the signal that the archetype has landed nowhere near the target.

We show the same construction for all three methods that produce a map between state spaces, so that the comparison is one of fits rather than of plotting conventions.
DSA is not included here because it aligns DMD operators and never constructs a map on state space, so it has no curve to draw in these coordinates; the corresponding comparison for DSA is the eigenspectra of Fig.~\ref{fig:eigenspectra}.

\subsubsection{Dynamical archetype analysis}
\begin{figure}[htbp]
    \centering
    \includegraphics[width=0.95\linewidth]{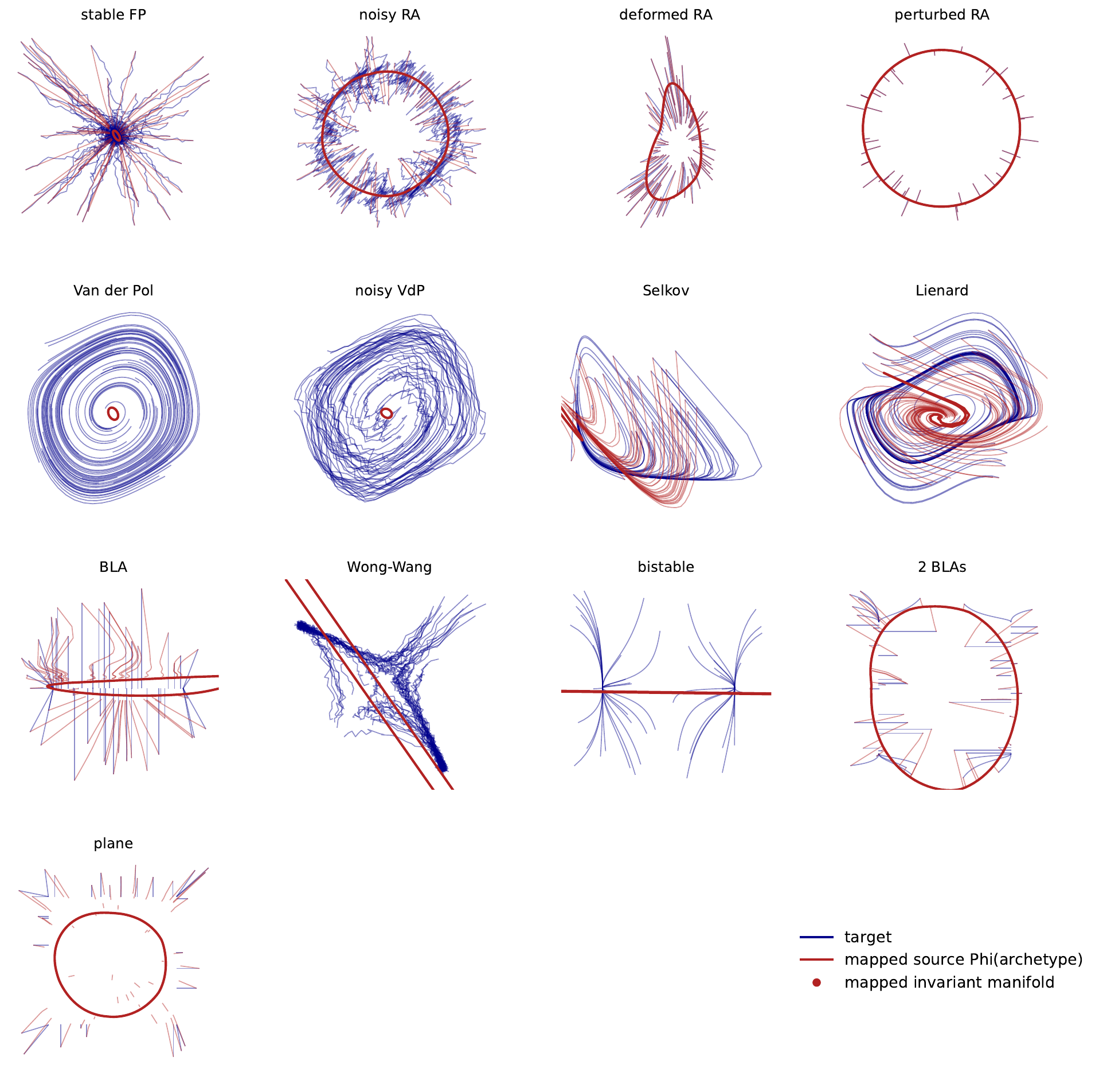}
    \caption{The \emph{ring attractor} archetype fitted to the different target systems.
    }
    \label{fig:ring_traj_invman}
\end{figure}

\begin{figure}[htbp]
    \centering
    \includegraphics[width=0.95\linewidth]{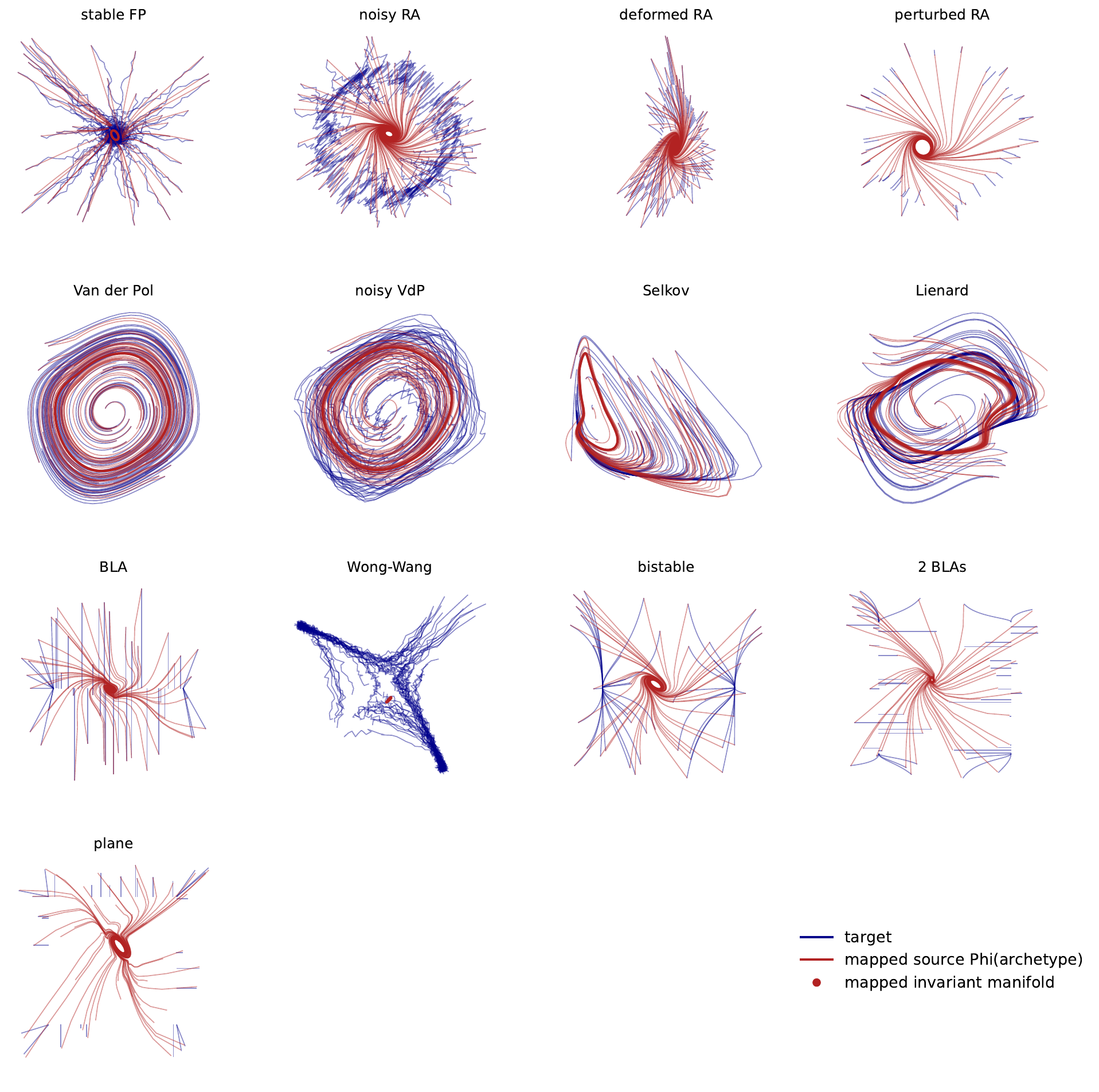}
    \caption{The \emph{limit cycle} archetype fitted to the different target systems.
    }
    \label{fig:lc_traj_invman}
\end{figure}

\begin{figure}[htbp]
    \centering
    \includegraphics[width=0.95\linewidth]{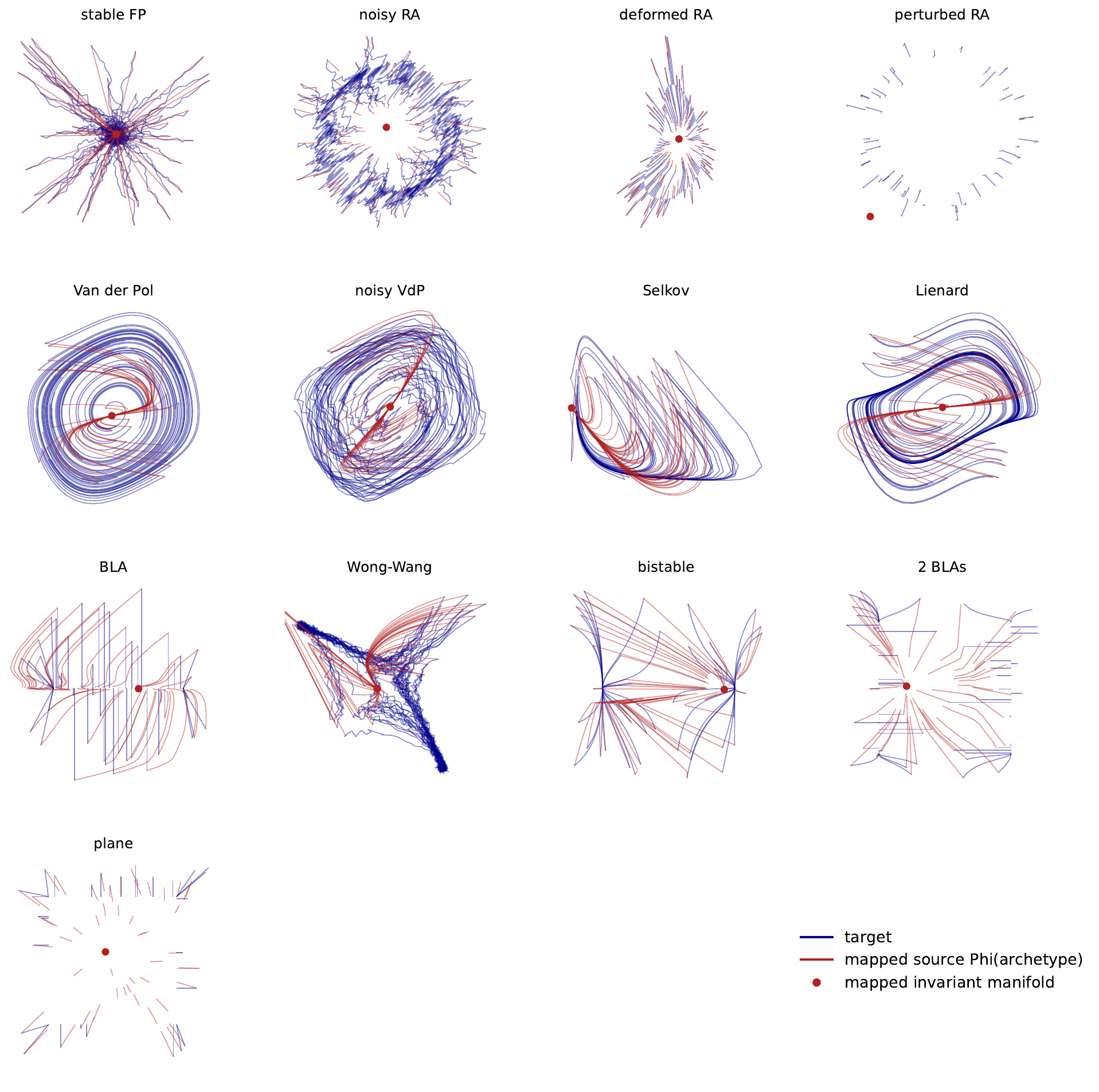}
    \caption{The \emph{single fixed point} archetype fitted to the different target systems.
    }
    \label{fig:lds_traj_invman}
\end{figure}

\begin{figure}[htbp]
    \centering
    \includegraphics[width=0.95\linewidth]{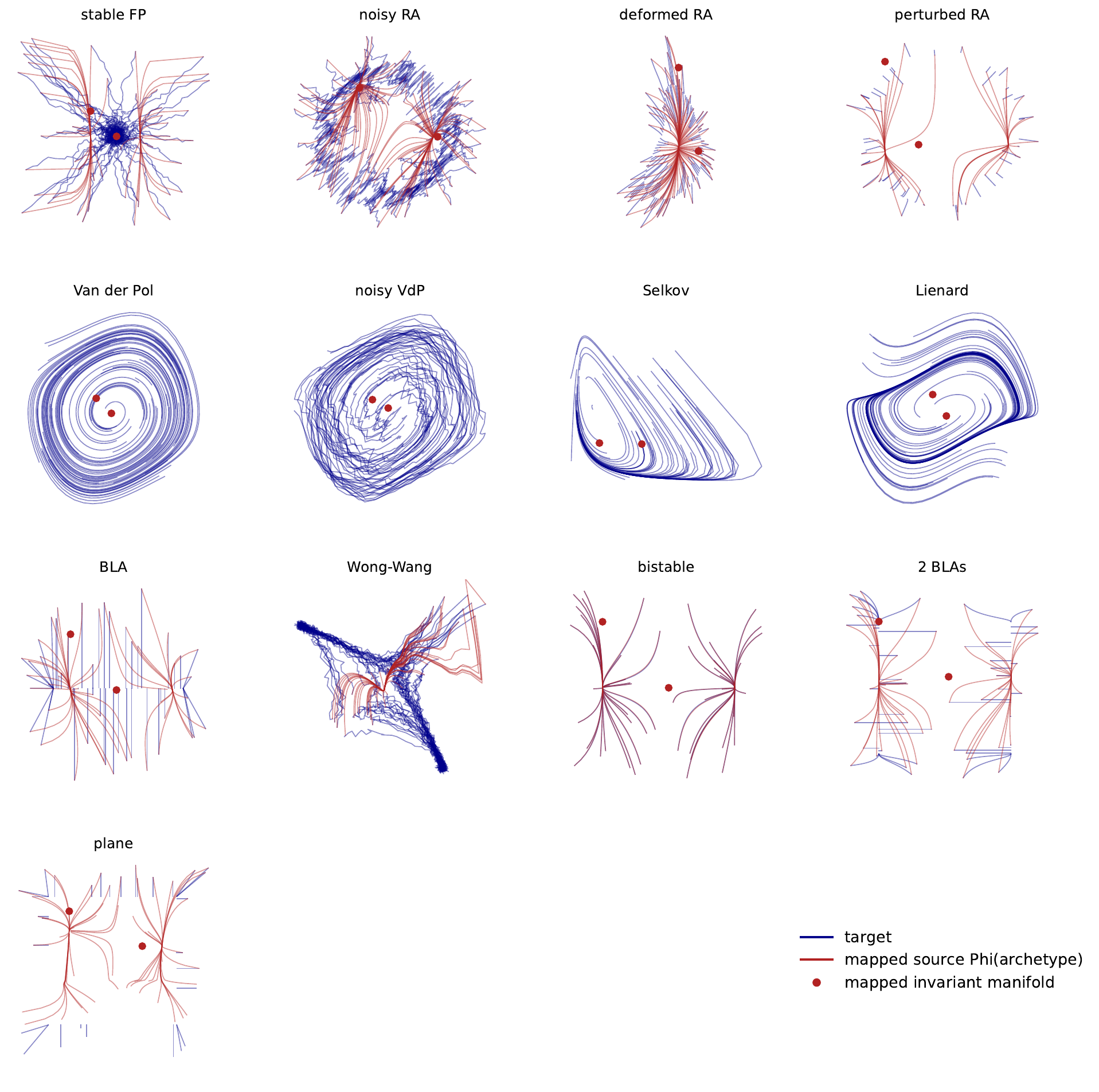}
    \caption{The \emph{bistable} archetype fitted to the different target systems.
    }
    \label{fig:bistable_traj_invman}
\end{figure}

\begin{figure}[htbp]
    \centering
    \includegraphics[width=0.95\linewidth]{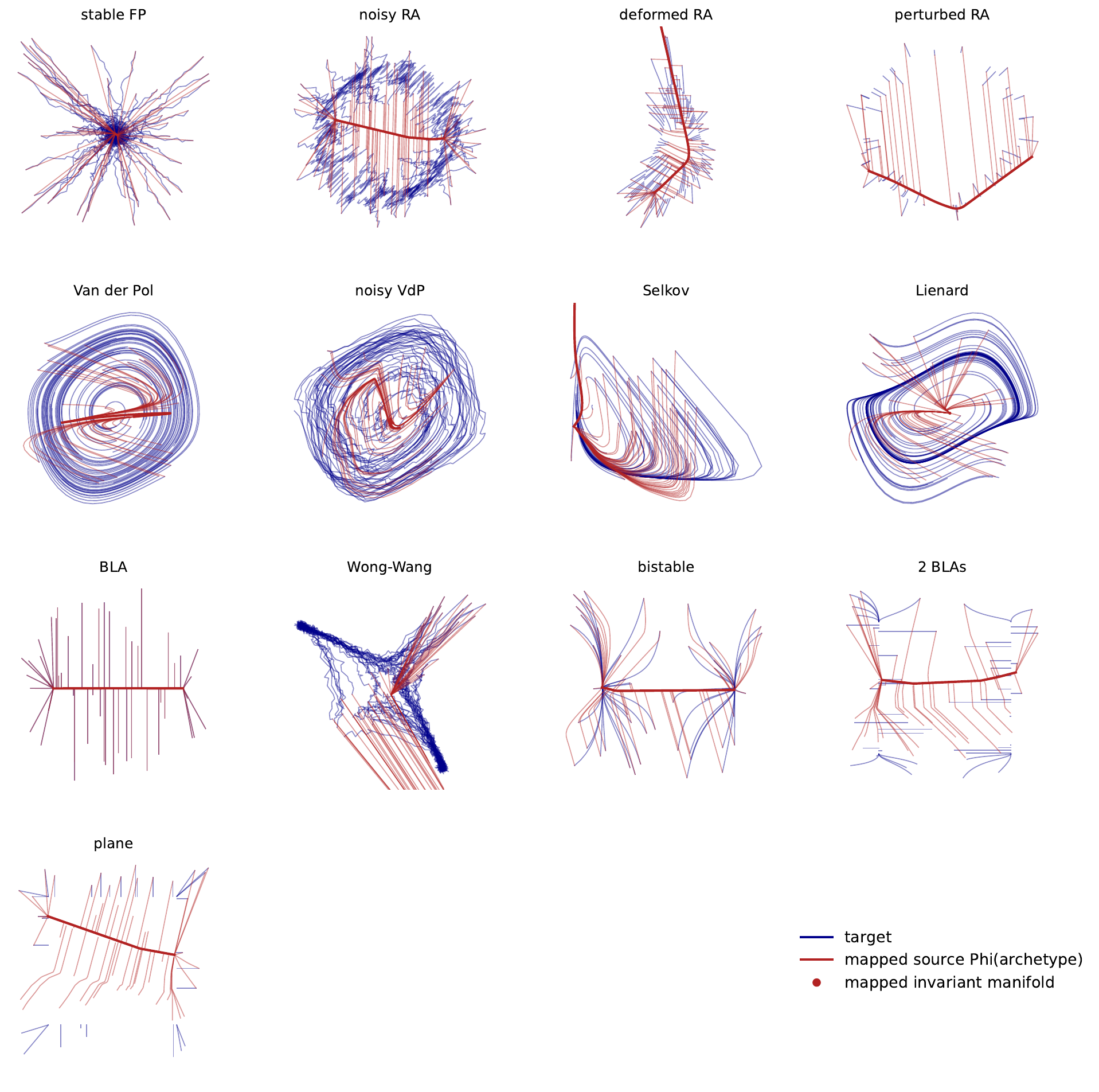}
    \caption{The \emph{Bounded Line Attractor} archetype fitted to the different target systems.
    }
    \label{fig:bla_traj_invman}
\end{figure}

\FloatBarrier
\subsubsection{Smooth prototype equivalences}
The same five archetypes fitted by SPE~\citep{friedman2025characterizing}, using the normalizing-flow diffeomorphism that SPE learns for each prototype--target pair.
SPE was not scored on the \emph{2 BLAs} and \emph{plane} targets, so those panels show the target alone and are marked accordingly.
The mapped trajectories are visibly more erratic here than for the other two methods.
This appears to be a matter of scale rather than of the fit failing outright: on the pairs we measured, SPE's map carries the target's initial conditions to radii of order $10$--$10^2$ in prototype coordinates, where the prototype's own invariant manifold sits at radius $1$, so the archetype trajectories traverse a large region before relaxing and their images magnify any local irregularity of the map.
The maps themselves are accurate inverses, with a round-trip error of order $10^{-7}$ relative to the data scale.

\begin{figure}[htbp]
    \centering
    \includegraphics[width=0.95\linewidth]{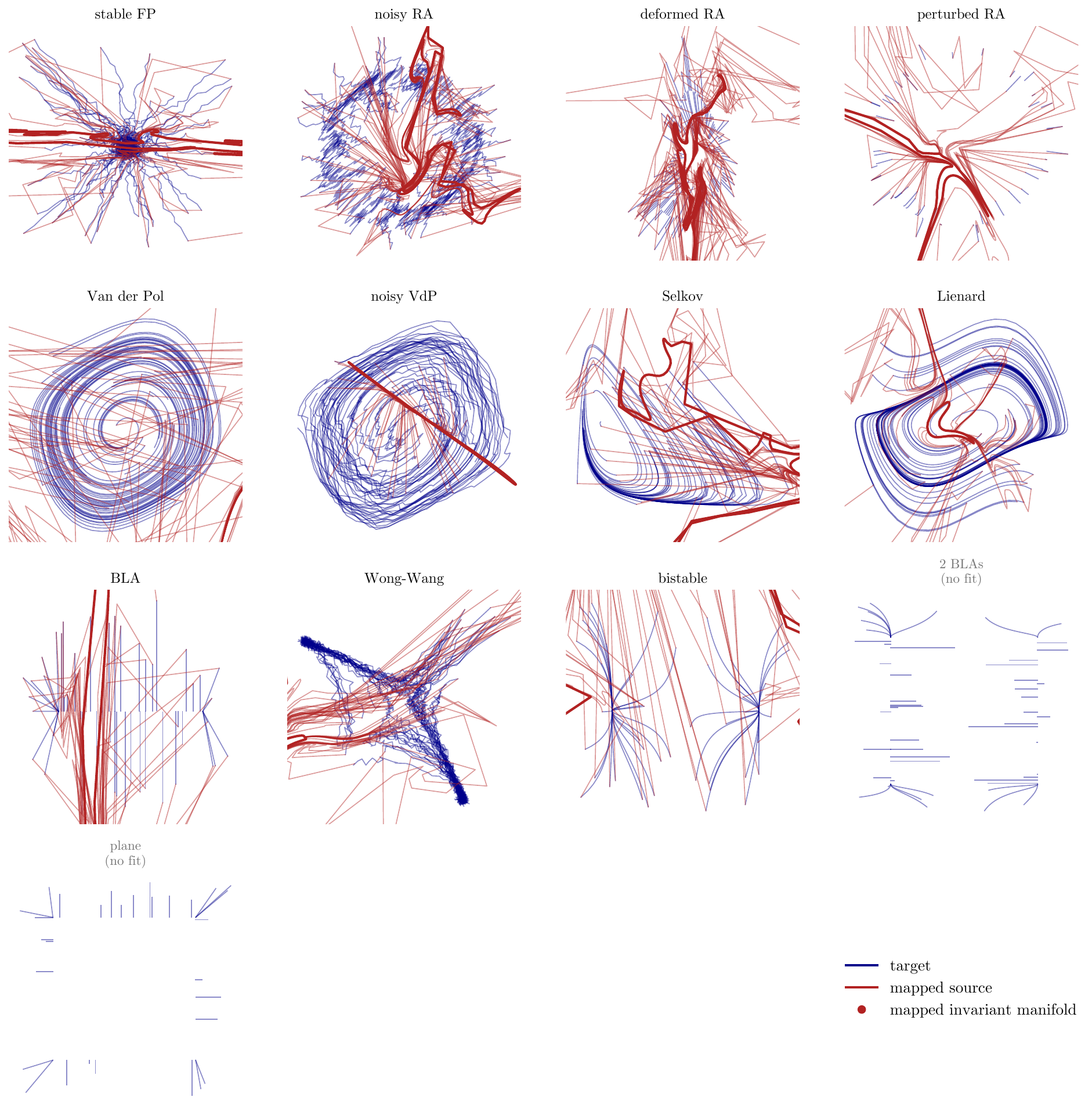}
    \caption{SPE: the \emph{ring attractor} prototype fitted to the different target systems.
    }
    \label{fig:ring_traj_invman_spe}
\end{figure}

\begin{figure}[htbp]
    \centering
    \includegraphics[width=0.95\linewidth]{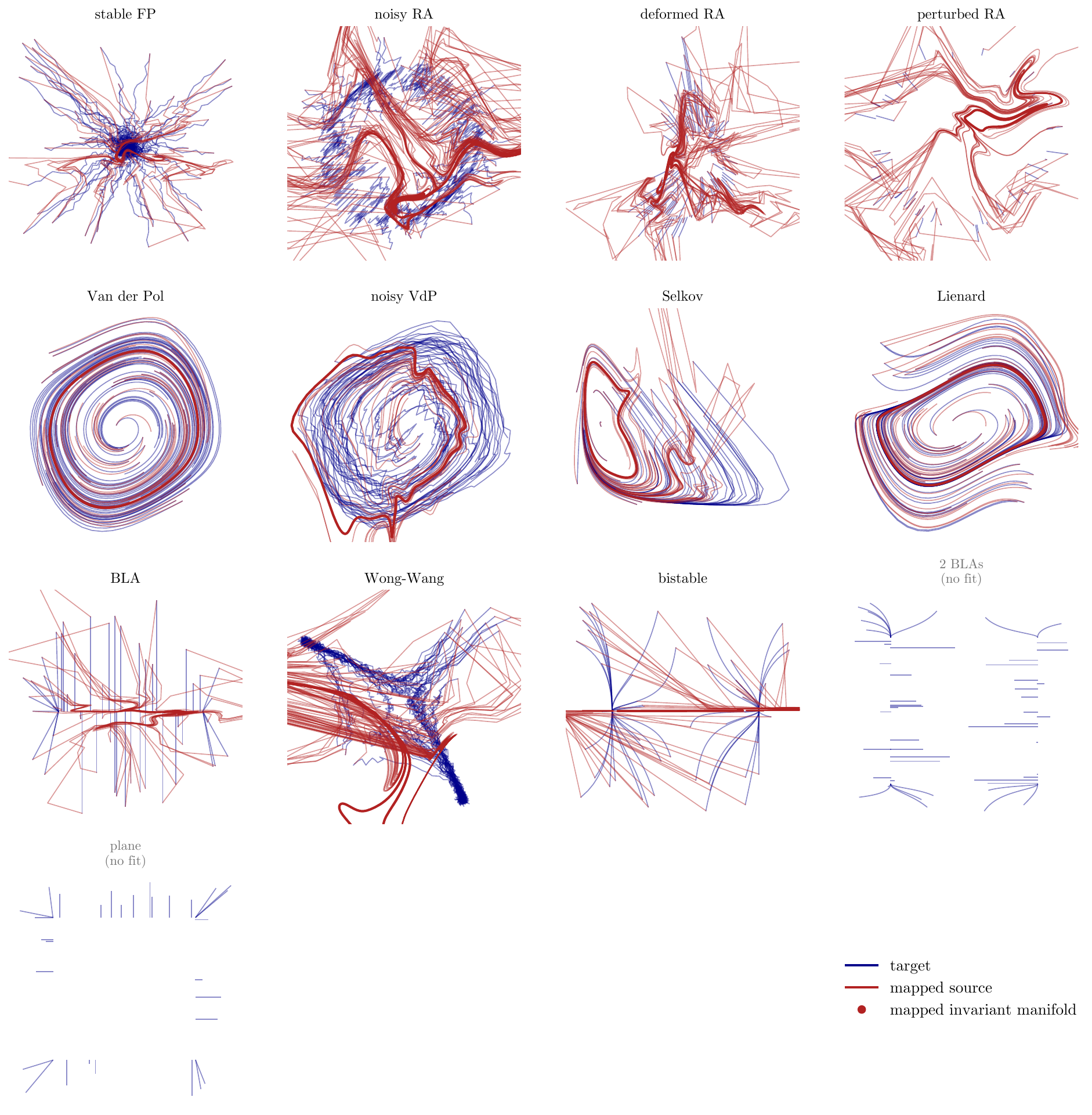}
    \caption{SPE: the \emph{limit cycle} prototype fitted to the different target systems.
    }
    \label{fig:lc_traj_invman_spe}
\end{figure}

\begin{figure}[htbp]
    \centering
    \includegraphics[width=0.95\linewidth]{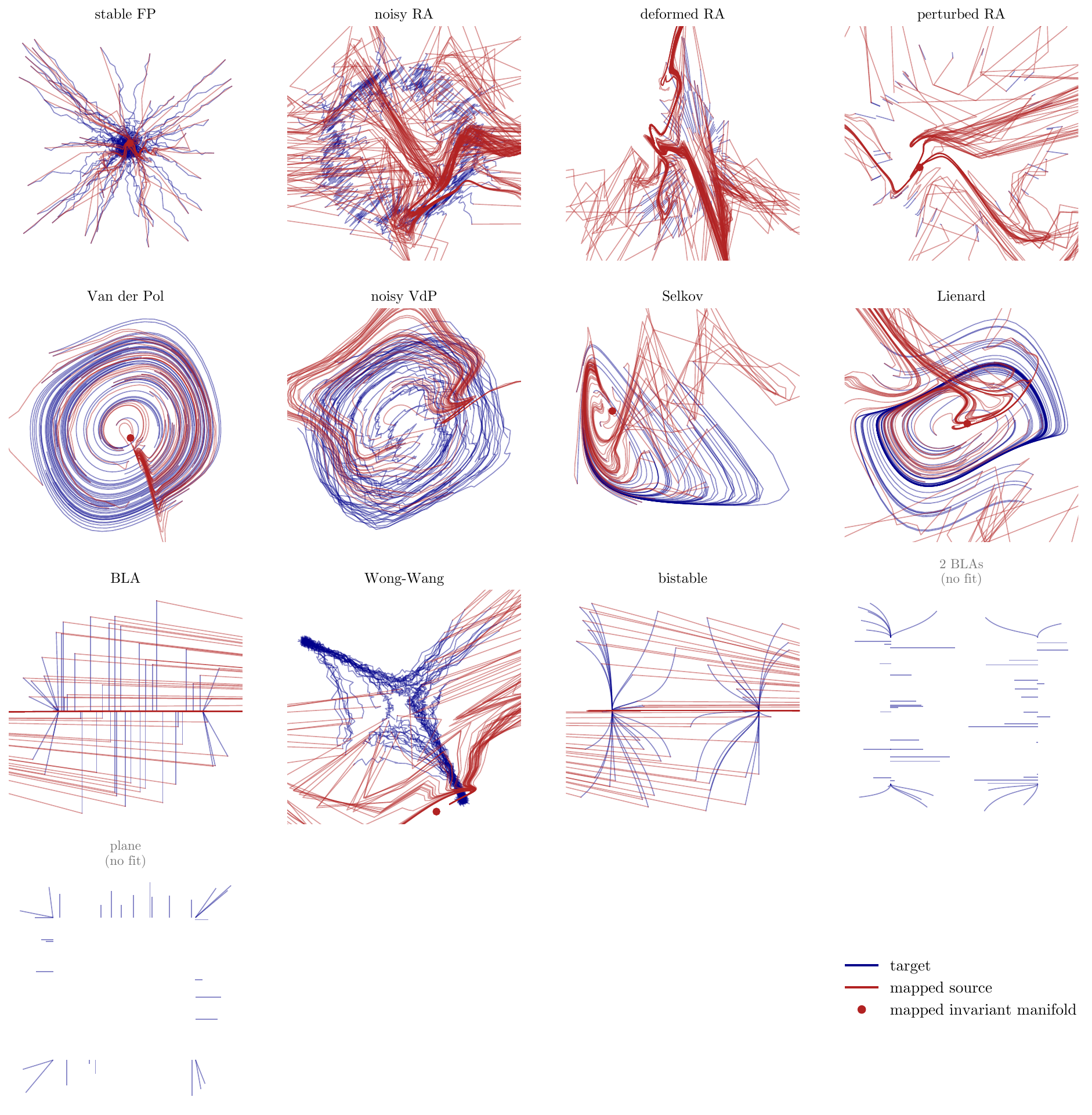}
    \caption{SPE: the \emph{single fixed point} prototype fitted to the different target systems.
    }
    \label{fig:lds_traj_invman_spe}
\end{figure}

\begin{figure}[htbp]
    \centering
    \includegraphics[width=0.95\linewidth]{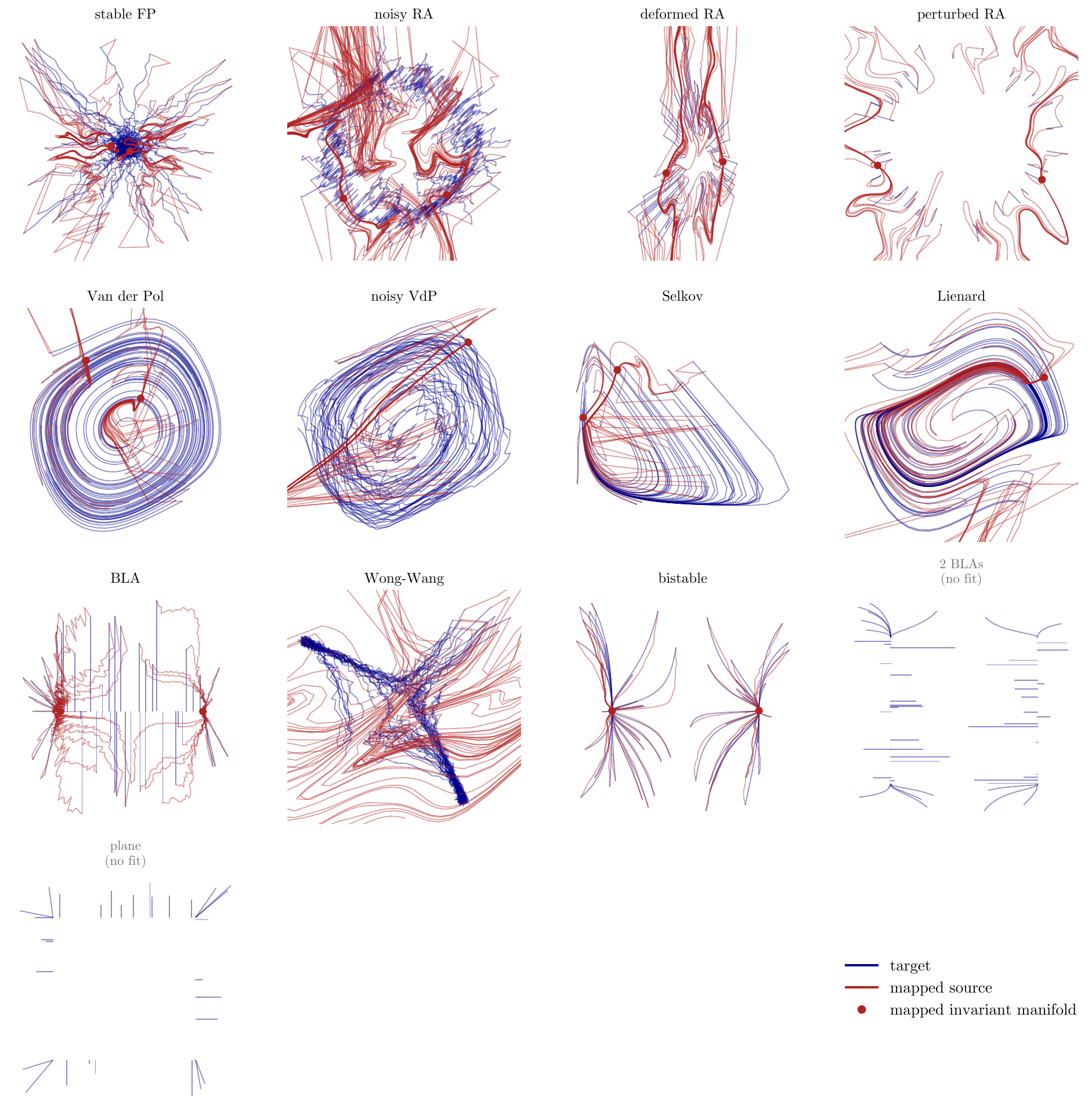}
    \caption{SPE: the \emph{bistable} prototype fitted to the different target systems.
    }
    \label{fig:bistable_traj_invman_spe}
\end{figure}

\begin{figure}[htbp]
    \centering
    \includegraphics[width=0.95\linewidth]{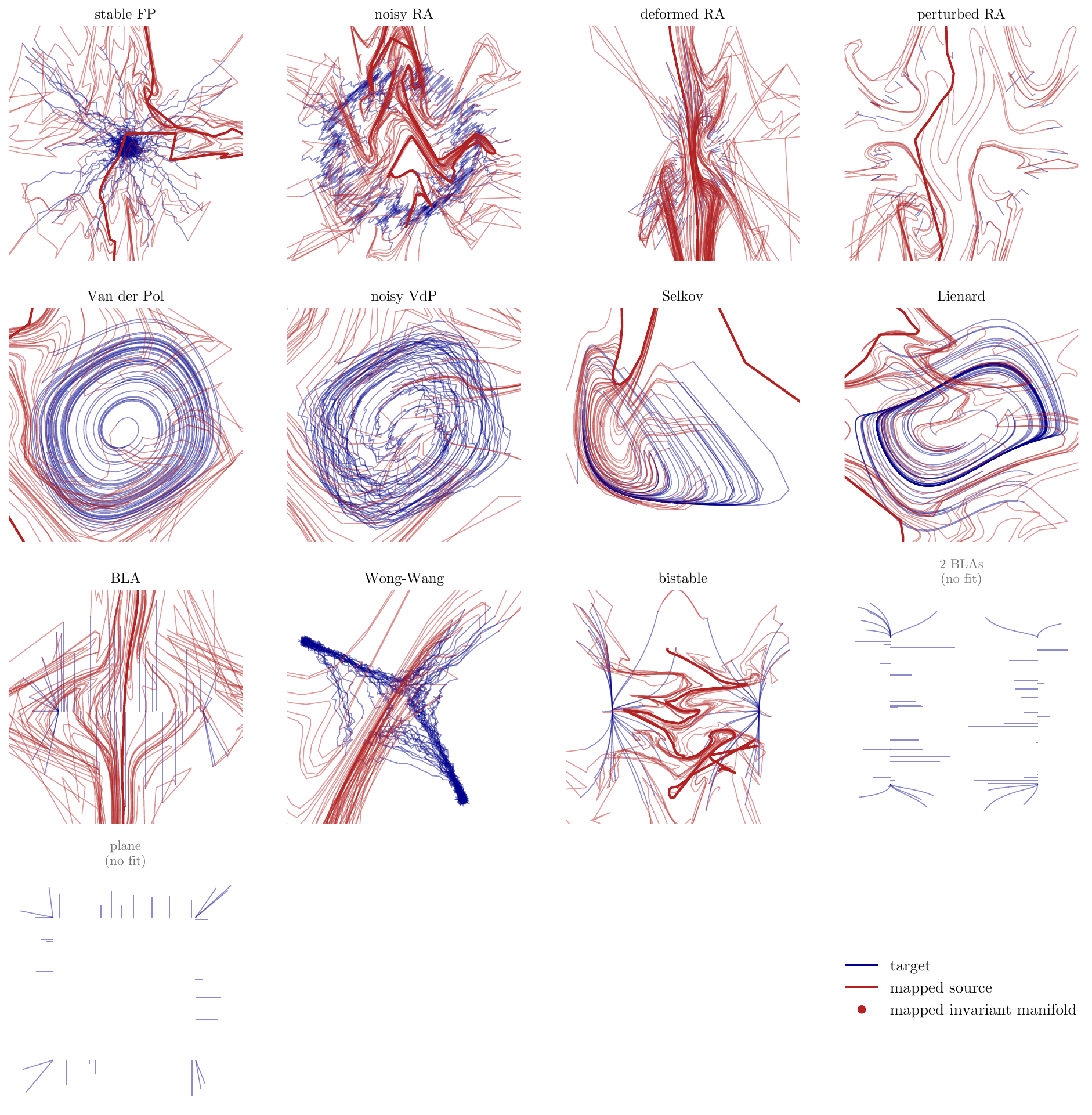}
    \caption{SPE: the \emph{Bounded Line Attractor} prototype fitted to the different target systems.
    }
    \label{fig:bla_traj_invman_spe}
\end{figure}

\FloatBarrier
\subsubsection{DFORM}
The same five archetypes fitted by DFORM~\citep{chen2024dform}.
The archetype is carried into the target's coordinates by the inverse deformation, which is the direction in which DFORM's learned map takes the archetype to the target.

\begin{figure}[htbp]
    \centering
    \includegraphics[width=0.95\linewidth]{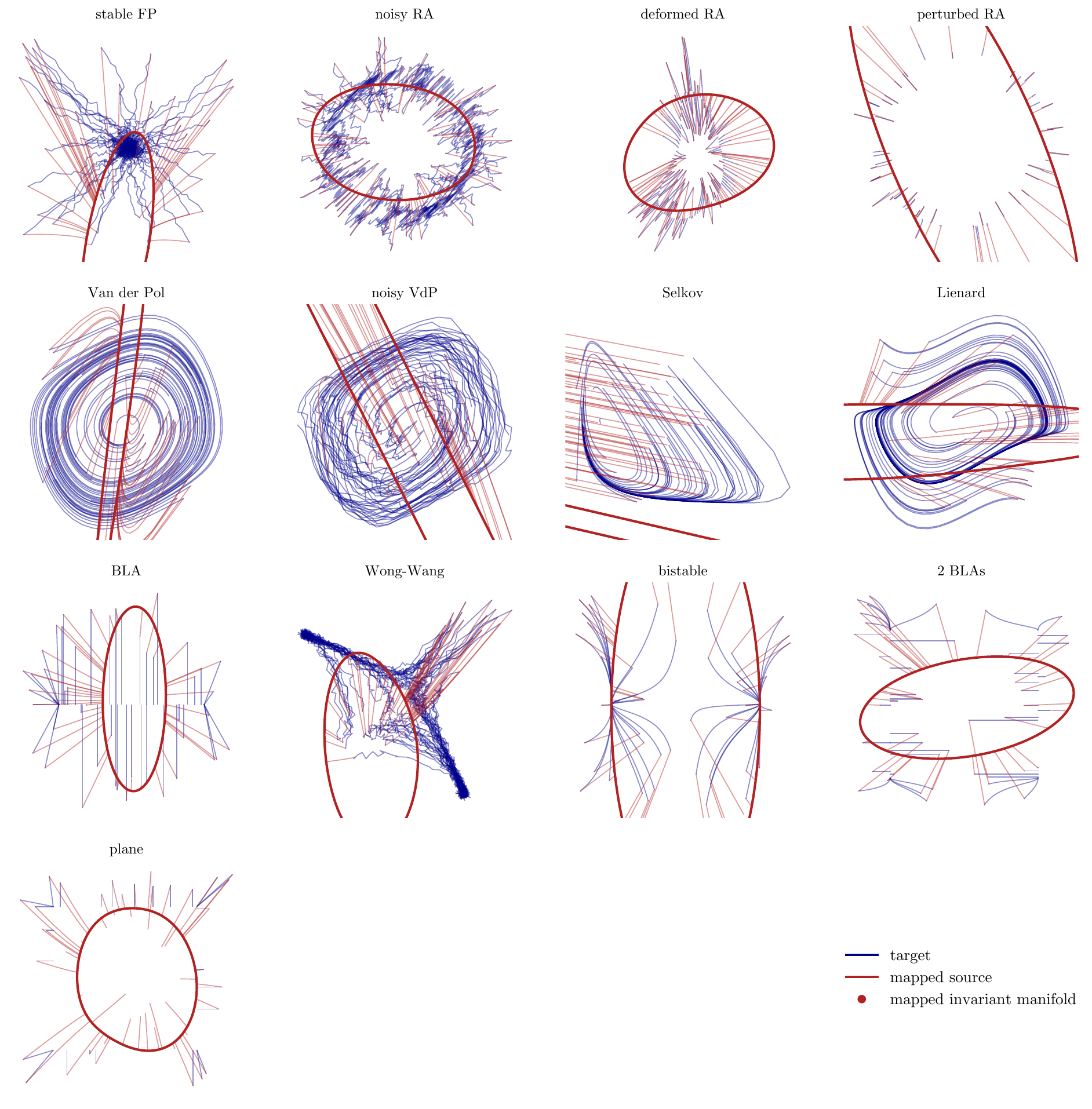}
    \caption{DFORM: the \emph{ring attractor} archetype fitted to the different target systems.
    Because the map is a diffeomorphism, the image of the ring remains a closed curve even for targets with no ring, such as the stable fixed point: a diffeomorphism cannot remove the archetype's topology, and the mismatch appears instead as the ring sitting away from the data.
    }
    \label{fig:ring_traj_invman_dform}
\end{figure}

\begin{figure}[htbp]
    \centering
    \includegraphics[width=0.95\linewidth]{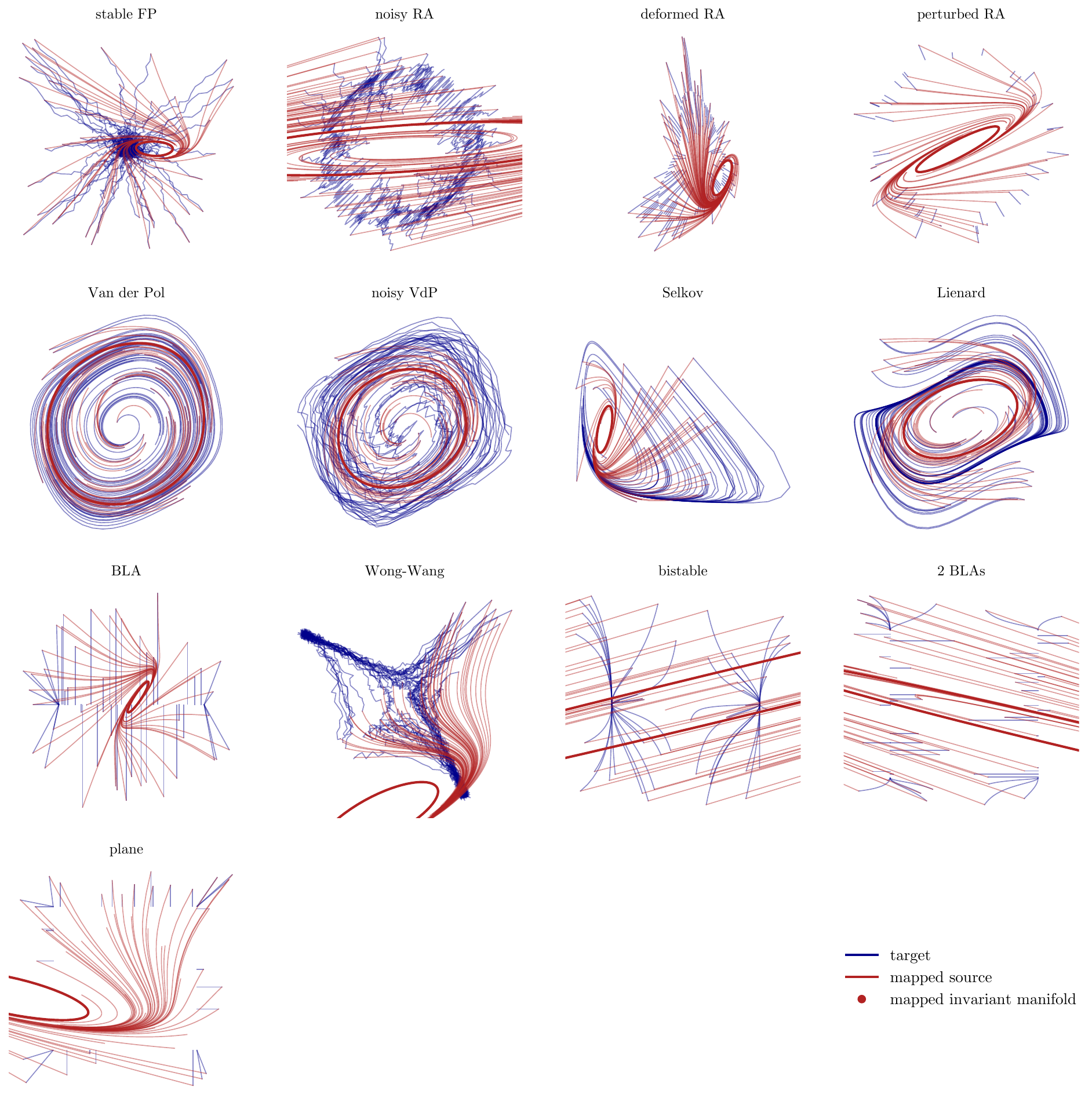}
    \caption{DFORM: the \emph{limit cycle} archetype fitted to the different target systems.
    }
    \label{fig:lc_traj_invman_dform}
\end{figure}

\begin{figure}[htbp]
    \centering
    \includegraphics[width=0.95\linewidth]{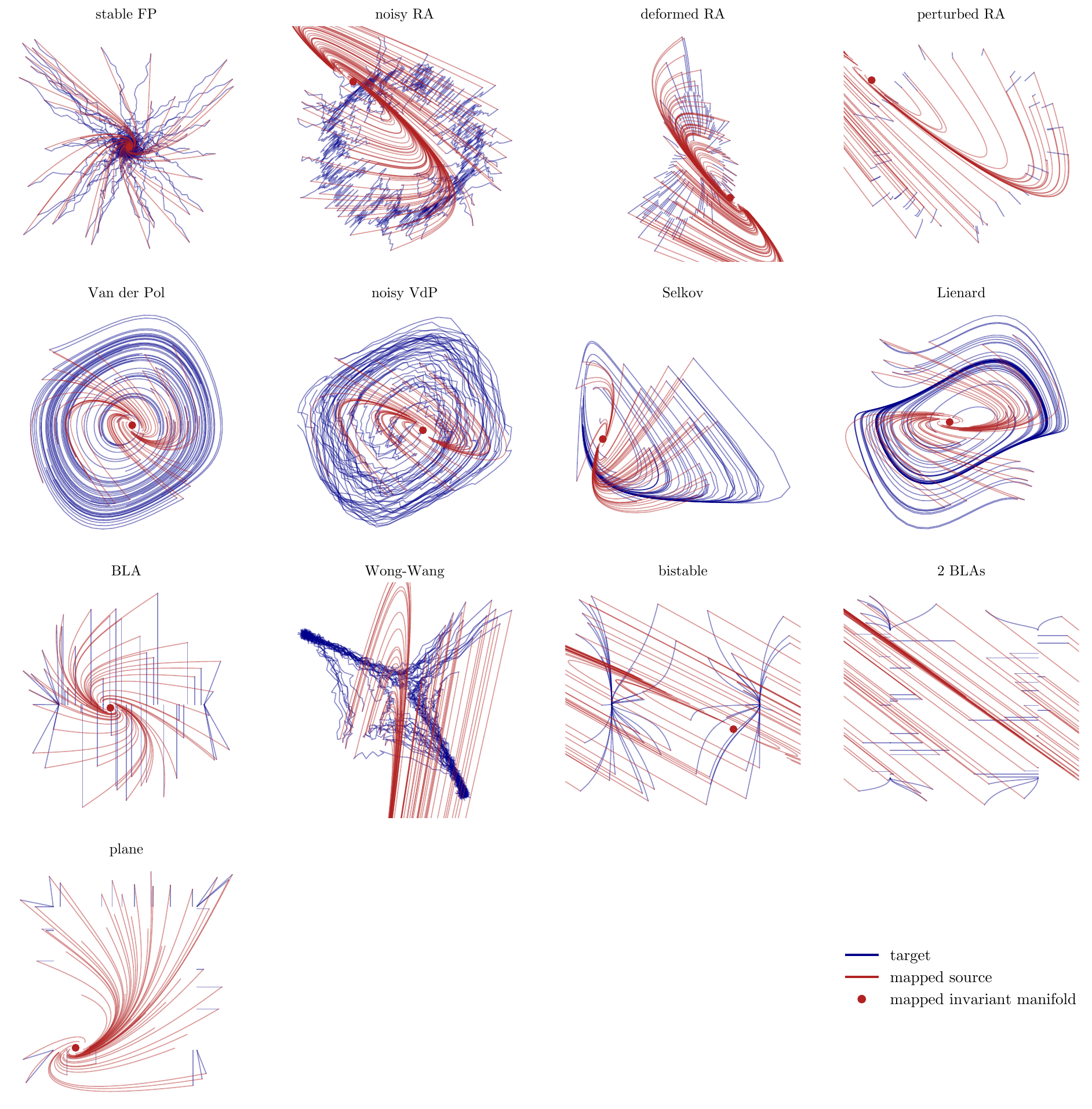}
    \caption{DFORM: the \emph{single fixed point} archetype fitted to the different target systems.
    }
    \label{fig:lds_traj_invman_dform}
\end{figure}

\begin{figure}[htbp]
    \centering
    \includegraphics[width=0.95\linewidth]{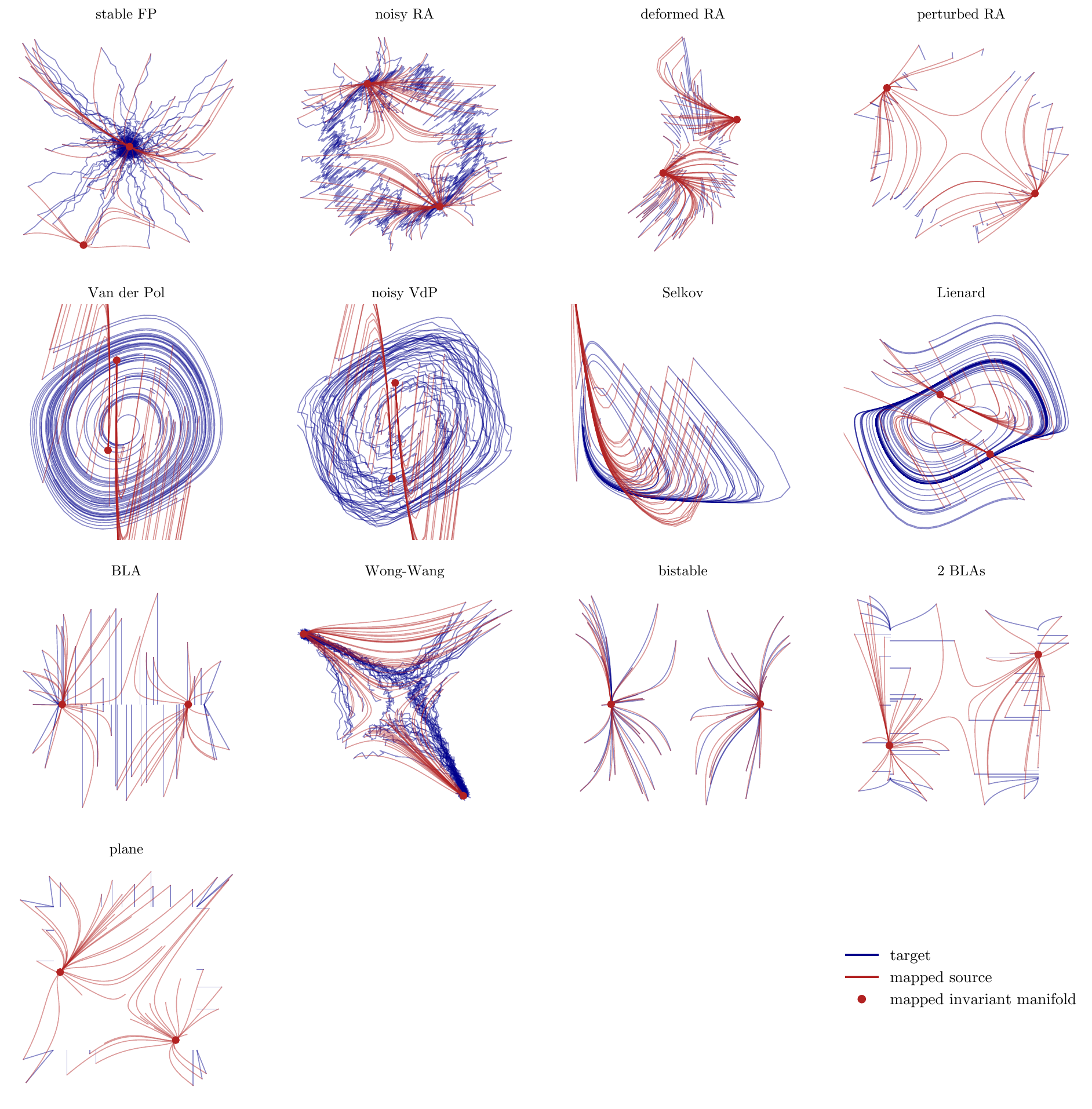}
    \caption{DFORM: the \emph{bistable} archetype fitted to the different target systems.
    }
    \label{fig:bistable_traj_invman_dform}
\end{figure}

\begin{figure}[htbp]
    \centering
    \includegraphics[width=0.95\linewidth]{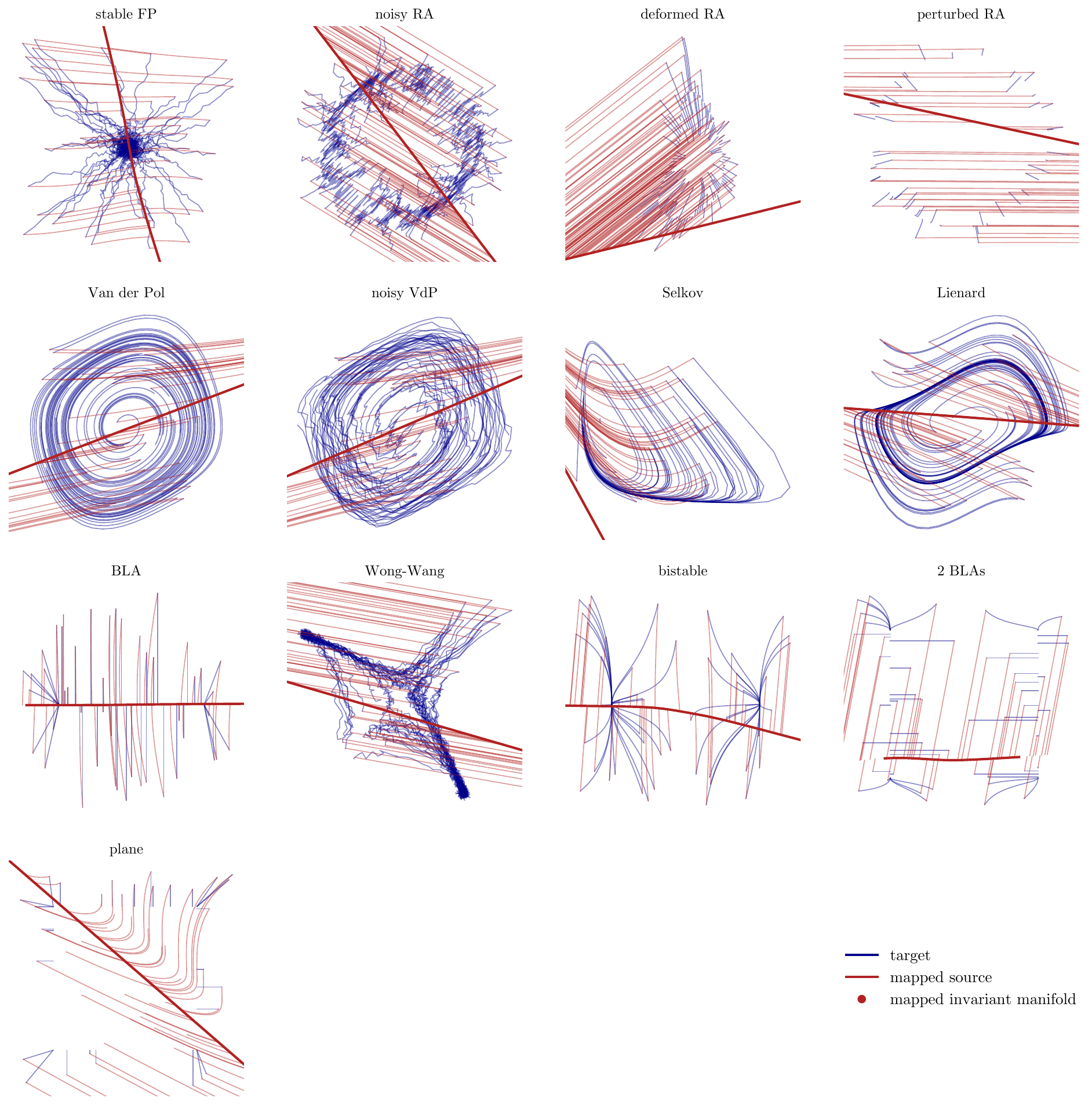}
    \caption{DFORM: the \emph{Bounded Line Attractor} archetype fitted to the different target systems.
    }
    \label{fig:bla_traj_invman_dform}
\end{figure}

\makeatletter
\if@submission
\makeatother
\newpage
\section*{NeurIPS Paper Checklist}

\begin{enumerate}

\item {\bf Claims}
    \item[] Question: Do the main claims made in the abstract and introduction accurately reflect the paper's contributions and scope?
    \item[] Answer: \answerYes{} % Replace by \answerYes{}, \answerNo{}, or \answerNA{}.
    \item[] Justification: The abstract and introduction states the claims of the paper and these claims are demonstrated through theory and numerical experiments.

\item {\bf Limitations}
    \item[] Question: Does the paper discuss the limitations of the work performed by the authors?
    \item[] Answer: \answerYes{} % Replace by \answerYes{}, \answerNo{}, or \answerNA{}.
    \item[] Justification: The limitations are described in the Discussion section of the paper.

\item {\bf Theory assumptions and proofs}
    \item[] Question: For each theoretical result, does the paper provide the full set of assumptions and a complete (and correct) proof?
    \item[] Answer: \answerYes{} % Replace by \answerYes{}, \answerNo{}, or \answerNA{}.
    \item[] Justification: Yes, the theoretical bounds are derived stepwise from results from the literature.

    \item {\bf Experimental result reproducibility}
    \item[] Question: Does the paper fully disclose all the information needed to reproduce the main experimental results of the paper to the extent that it affects the main claims and/or conclusions of the paper (regardless of whether the code and data are provided or not)?
    \item[] Answer: \answerYes{} % Replace by \answerYes{}, \answerNo{}, or \answerNA{}.
    \item[] Justification: All experimental setting are described in full detail in the Supplementary.

\item {\bf Open access to data and code}
    \item[] Question: Does the paper provide open access to the data and code, with sufficient instructions to faithfully reproduce the main experimental results, as described in supplemental material?
    \item[] Answer: \answerYes{} % Replace by \answerYes{}, \answerNo{}, or \answerNA{}.
    \item[] Justification: All experimental methods are described in full detail in the Supplementary. 
    \item[] Guidelines:
    \begin{itemize}
        \item The answer NA means that paper does not include experiments requiring code.
        \item Please see the NeurIPS code and data submission guidelines (\url{https://nips.cc/public/guides/CodeSubmissionPolicy}) for more details.
        \item While we encourage the release of code and data, we understand that this might not be possible, so ``No'' is an acceptable answer. Papers cannot be rejected simply for not including code, unless this is central to the contribution (e.g., for a new open-source benchmark).
        \item The instructions should contain the exact command and environment needed to run to reproduce the results. See the NeurIPS code and data submission guidelines (\url{https://nips.cc/public/guides/CodeSubmissionPolicy}) for more details.
        \item The authors should provide instructions on data access and preparation, including how to access the raw data, preprocessed data, intermediate data, and generated data, etc.
        \item The authors should provide scripts to reproduce all experimental results for the new proposed method and baselines. If only a subset of experiments are reproducible, they should state which ones are omitted from the script and why.
        \item At submission time, to preserve anonymity, the authors should release anonymized versions (if applicable).
        \item Providing as much information as possible in supplemental material (appended to the paper) is recommended, but including URLs to data and code is permitted.
    \end{itemize}

\item {\bf Experimental setting/details}
    \item[] Question: Does the paper specify all the training and test details (e.g., data splits, hyperparameters, how they were chosen, type of optimizer, etc.) necessary to understand the results?
    \item[] Answer: \answerYes{} % Replace by \answerYes{}, \answerNo{}, or \answerNA{}.
    \item[] Justification: Yes, these are fully documented in the Supplementary or the Code.

\item {\bf Experiment statistical significance}
    \item[] Question: Does the paper report error bars suitably and correctly defined or other appropriate information about the statistical significance of the experiments?
    \item[] Answer: \answerYes{} % Replace by \answerYes{}, \answerNo{}, or \answerNA{}.
    \item[] Justification: For the relevant experiments, error bars are included.

\item {\bf Experiments compute resources}
    \item[] Question: For each experiment, does the paper provide sufficient information on the computer resources (type of compute workers, memory, time of execution) needed to reproduce the experiments?
    \item[] Answer: \answerYes{} % Replace by \answerYes{}, \answerNo{}, or \answerNA{}.
    \item[] Justification: The used computer resources are detailed in the Supplementary.
    
\item {\bf Code of ethics}
    \item[] Question: Does the research conducted in the paper conform, in every respect, with the NeurIPS Code of Ethics \url{https://neurips.cc/public/EthicsGuidelines}?
    \item[] Answer: \answerYes{} % Replace by \answerYes{}, \answerNo{}, or \answerNA{}.
    \item[] Justification: The authors read and comply with the  NeurIPS Code of Ethics.

\item {\bf Broader impacts}
    \item[] Question: Does the paper discuss both potential positive societal impacts and negative societal impacts of the work performed?
    \item[] Answer: \answerNA{} % Replace by \answerYes{}, \answerNo{}, or \answerNA{}.
    \item[] Justification: The paper presents theoretical results that are not directly impactful to society at large.
    
\item {\bf Safeguards}
    \item[] Question: Does the paper describe safeguards that have been put in place for responsible release of data or models that have a high risk for misuse (e.g., pretrained language models, image generators, or scraped datasets)?
    \item[] Answer: \answerNA{} % Replace by \answerYes{}, \answerNo{}, or \answerNA{}.
    \item[] Justification: The paper does not present models that pose a risk of misuse.

\item {\bf Licenses for existing assets}
    \item[] Question: Are the creators or original owners of assets (e.g., code, data, models), used in the paper, properly credited and are the license and terms of use explicitly mentioned and properly respected?
    \item[] Answer: \answerYes{} % Replace by \answerYes{}, \answerNo{}, or \answerNA{}.
    \item[] Justification: The usage of code is well-documented in the Supplementary.

\item {\bf New assets}
    \item[] Question: Are new assets introduced in the paper well documented and is the documentation provided alongside the assets?
    \item[] Answer: \answerNA{} % Replace by \answerYes{}, \answerNo{}, or \answerNA{}.
    \item[] Justification: The paper does not release new assets.

\item {\bf Crowdsourcing and research with human subjects}
    \item[] Question: For crowdsourcing experiments and research with human subjects, does the paper include the full text of instructions given to participants and screenshots, if applicable, as well as details about compensation (if any)? 
    \item[] Answer: \answerNA{} % Replace by \answerYes{}, \answerNo{}, or \answerNA{}.
    \item[] Justification: The paper does not involve crowdsourcing nor research with human subjects.

\item {\bf Institutional review board (IRB) approvals or equivalent for research with human subjects}
    \item[] Question: Does the paper describe potential risks incurred by study participants, whether such risks were disclosed to the subjects, and whether Institutional Review Board (IRB) approvals (or an equivalent approval/review based on the requirements of your country or institution) were obtained?
    \item[] Answer: \answerNA{} % Replace by \answerYes{}, \answerNo{}, or \answerNA{}.
    \item[] Justification: The paper does not involve crowdsourcing nor research with human subjects..

\item {\bf Declaration of LLM usage}
    \item[] Question: Does the paper describe the usage of LLMs if it is an important, original, or non-standard component of the core methods in this research? Note that if the LLM is used only for writing, editing, or formatting purposes and does not impact the core methodology, scientific rigorousness, or originality of the research, declaration is not required.
    %this research? 
    \item[] Answer: \answerNA{} % Replace by \answerYes{}, \answerNo{}, or \answerNA{}.
    \item[] Justification: The paper's core methodology was not directly dependent on any LLMs.

\end{enumerate}
\fi

\end{document}